\documentclass[10pt,twoside]{article}
\usepackage{repdefs,latexsym,longtable,graphicx,varioref,cite,subfigure,comment,ulem}
\usepackage[usenames,dvipsnames]{xcolor}

\newcommand{\xs}{x^{\ast}}
\newcommand{\al}{\alpha}
\newcommand{\au}{\underline{\alpha}}
\newcommand{\gm}{\gamma}

\newcommand{\rr}{\rho}
\newcommand{\dd}{\delta}

\newcommand{\xb}{\bar x}
\newcommand{\RR}{\mathrm{I\!R\!}}
\newcommand{\PP}{\mathrm{I\!P\!}}
\newcommand{\EE}{\mathrm{I\!E\!}}
\newcommand{\Gm}{\Gamma}
\newcommand{\Lm}{\Lambda}
\newcommand{\sm}{\setminus}
\newcommand{\IL}{\mathcal{IL}}
\newcommand{\IU}{\mathcal{IU}}
\newcommand{\LL}{\mathcal{L}}
\newcommand{\UU}{\mathcal{U}}

\newcommand{\IFF}{\mathcal{IF}}
\newcommand{\astar}{a^{\ast}(\dd,\rr;\zeta)}
\newcommand{\ra}{\rightarrow}
\newcommand{\rh}{\hat{\rho}}

\newcommand{\bNN}{\overline{NSP_{\au}}}
\newcommand{\bNa}{\overline{NSP_{\al}}}
\newcommand{\mN}{\mathcal{N}}
\newcommand{\mW}{\mathcal{W}}
\newcommand{\lm}{\lambda}
\newcommand{\Sg}{\sigma}
\newcommand{\sg}{\sigma}
\newcommand{\sti}{\sigma}
\newcommand{\Ta}{\Theta^1_n}
\newcommand{\Tb}{\Theta^2_n}
\newcommand{\et}{e}
\newcommand{\FF}{\mathcal{F}}

\newcommand{\paperauthor}{Coralia Cartis\hspace*{-0.02cm}\footnotemark[1]\,\,\, and\,
 Andrew Thompson\footnotemark[2]}
\newcommand{\papertitle}{A new and improved quantitative recovery analysis for iterative hard thresholding algorithms in compressed sensing}

\title{\papertitle}
\author{\paperauthor}

\newcommand{\theabstract}{
We present a new recovery analysis for a standard compressed sensing algorithm, Iterative Hard Thresholding (IHT) (Blumensath and Davies, 2008), which considers the fixed points of the algorithm. In the context of arbitrary measurement matrices, we derive a sufficient condition for  convergence of IHT to a fixed point and a necessary condition for the existence of fixed points. These conditions allow us to perform  a sparse signal recovery analysis in the deterministic noiseless case by implying that the original sparse signal is the unique fixed point and limit point of IHT,
and in the case of Gaussian measurement matrices and noise by generating a bound on the approximation error of the IHT limit as a multiple of the noise level. By generalizing the notion of fixed points, we extend our analysis to the variable stepsize Normalised IHT (N-IHT) (Blumensath and Davies, 2010).  For both stepsize schemes, we obtain lower bounds on asymptotic phase transitions in a proportional-dimensional framework, quantifying the sparsity/undersampling trade-off for which recovery is guaranteed.
Exploiting the reasonable average-case assumption that the underlying signal and measurement matrix are independent,
comparison with previous results within this framework shows a substantial quantitative improvement.

}

\begin{document}
\renewcommand{\thefootnote}{\fnsymbol{footnote}}
\footnotetext[1]{Mathematical Institute, University of Oxford,
            Radcliffe Observatory Quarter, Oxford,
           OX2 6GG, UK.\\
       email: coralia.cartis@maths.ox.ac.uk.}
\footnotetext[2]{Department of Mathematics, Duke University,
            Box 90320, Durham, NC 27708-0320, USA.\\
       email: thompson@math.duke.edu.}

\renewcommand{\thefootnote}{\arabic{footnote}}

\maketitle

\nocite{*}
\vspace*{-10mm}

\begin{abstract}
\theabstract
\end{abstract}

\setcounter{page}{1}
\numsection{Introduction}
\label{intro}

Compressed Sensing (CS) seeks to recover sparse or compressible
signals from undersampled linear measurements \cite{Candes, robust, CS}; it asserts that the
number of measurements should be proportional to the information
content of the signal, rather than its dimension. More specifically, one seeks to recover a sparse signal from noisy linear measurements. We refer to a signal which has at most $k$ nonzero entries as being $k$-sparse, and the problem can be stated as follows.\\
\textbf{Sparse recovery from noisy measurements:} \textit{Recover a} $k$\textit{-sparse signal} $x^*\in\Re^N$ \textit{from the linear measurements}
\eqn{sparserecovery}{b=Ax^*+e\in\Re^n,}
\textit{where} $e\in\Re^n$ \textit{is an unknown noise vector and where} $0<2k\le n\le N$\textit{.}\\
Since the introduction of CS in 2004, many algorithms have been
proposed to solve this seemingly (and generally) intractable problem;
see \cite{TroppWright} for a recent survey. A common approach is
to solve, by means of classical or recently-proposed optimization
methods, a (convex or nonconvex) optimization relaxation
that penalizes the lack of sparsity of $x$ by means of $l_p$-norms
with $0<p \leq 1$. Alternatively, greedy methods --- such as (Orthogonal) Matching Pursuit \cite{mallat, omp1, omp2},
SP \cite{sp}, CoSAMP \cite{cosamp}, amongst others --- can be used to tackle
the so-called $l_0$-problem directly, namely,
\eqn{c:l0problem}{
\min_{x\in \smallRe^N} \Psi(x)\eqdef \half\|Ax-b\|^2 \tim{subject to} \|x\|_0\leq k,
}
where $\|\cdot\|_0$ counts the number of nonzero entries of the
argument. Problem \req{c:l0problem} is  nonconvex, with many local
minimizers and in the perfect case of zero noise, with a (unique) global
minimizer at the $k$-sparse vector $x^*$ \cite{Mila} that we are aiming to recover. It is now well known that, under certain conditions, many of these algorithms have stable recovery properties, namely that the error in approximating the original signal is some (usually small) multiple of the noise level, which further implies exact recovery of the original signal $x^*$ in the absence of noise \cite{TroppWright, greedy}.

Here, we focus on a  simple and widely-used greedy technique --
Iterative Hard Thresholding (IHT) \cite{thresh, normalized} -- that generates
feasible steepest descent steps for problem \req{c:l0problem}, obtained by projecting  steps along the negative
gradient direction of $\Psi$ onto the $l_0$-norm constraint by means of the \textit{hard
threshold} operator which simply sets all but the $k$ largest in magnitude coefficients of a
vector to zero.
IHT \cite{thresh} performs gradient projection with constant stepsize, while Normalised Iterative Hard Thresholding (N-IHT) \cite{normalized} employs a variable stepsize scheme.

Early CS theory focussed upon algorithms which solve convex relaxations of \req{c:l0problem} \cite{robust,CS}, and perhaps for this reason IHT algorithms were slow to gain acceptance in the CS community. However, more recently, empirical studies in \cite{large_scale} have shown that, surprisingly, these nonconvex approaches are in fact competitive in practice in terms of sparse recovery properties with more established CS algorithms based on $l_1$-minimization.
In common with many gradient methods proposed for $l_1$-based CS recovery, IHT algorithms also have low computational cost, with the most costly operations being matrix-vector products and hard threshold operations. As we detail below, existing theoretical recovery analyses are unduly pessimistic and
fail to account for this excellent practical behaviour of IHT variants,  and it is the aim of this paper to improve quantitative recovery guarantees of IHT algorithms
by means of a new probabilistic analysis.



Regarding state-of-the-art
theoretical properties, Blumensath and Davies \cite{thresh}
obtained the first convergence result for IHT, proving convergence
to a fixed point of IHT/local minimizer of (\ref{c:l0problem})
provided the spectral norm of the measurement matrix $A$ is less
than one, a somewhat restrictive condition. The same authors
\cite{threshRIP} then proved that stable recovery is guaranteed
provided $A$ satisfies a
\textit{restricted isometry property} (RIP) \cite{Candes}, which
requires the matrix to act as a near isometry on all $k$-sparse
vectors, a now ubiquitous tool in CS recovery analysis. Other
RIP-based recovery conditions were subsequently obtained for IHT in
\cite{gargkandekhar,foucart,HTP}, and for N-IHT in
\cite{normalized}.

Determining whether a given measurement matrix satisfies a restricted isometry property is in itself, however, NP-hard. It has been shown \cite{simple_proof} that certain random matrices, such as Gaussian matrices in which each entry is i.i.d. Gaussian, satisfy the RIP provided
$$n\geq C\cdot k\ln\left(\frac{N}{k}\right).$$
Quantifying the constant $C$ is, however, vital to practitioners who wish to know how aggressively a signal may be undersampled given its dimension and sparsity. Based on the RIP analysis in \cite{threshRIP,gargkandekhar,normalized,foucart,HTP}, quantitative results were obtained for IHT in \cite{greedy,thesis} for the case of Gaussian matrices in an asymptotic framework in which the problem dimensions are assumed to grow proportionally. We will refer to such a framework as the \textit{proportional-growth asymptotic}, defined as follows.

\begin{definition}[\textbf{Proportional-growth asymptotic~\cite{lqphase}}]\label{propdimdef}
We say that a sequence of problem sizes $(k,n,N)$, where $0<k\le n\le N$, grows proportionally if, for some $\dd\in(0,1]$ and $\rr\in(0,1]$,
$$\frac{n}{N}\ra\dd\;\;\;\;\mbox{and}\;\;\;\;\frac{k}{n}\ra\rr\;\;\;\;\mbox{as}\;(k,n,N)\ra\infty.$$
\end{definition}

This framework, advocated by Donoho and others~\cite{neighborliness,precise}, defines a two-dimensional phase space for asymptotic analysis in which the variables $\dd$ and $\rr$ have a simple practical interpretation. The parameter $\dd$ is the ratio by which the signal is undersampled (an \textit{undersampling ratio}), while the ratio $\rr$ indicates how many measurements need to be taken as a multiple of the sparsity (an \textit{oversampling ratio}).

By making use of RIP analysis for Gaussian matrices, first performed in \cite{lqphase} and subsequently improved upon in \cite{BT}, it was shown in \cite{greedy,thesis} that all of the RIP conditions proved to date for IHT algorithms are pessimistic compared to these algorithms' numerically-observed average-case behaviour. This is not altogether surprising, since the RIP gives worst-case guarantees. There is, therefore, a need for improved quantitative recovery guarantees for IHT algorithms which narrow the gap between theoretical guarantees and observed performance. This is in contrast to $l_1$-minimization, for which average-case phase transitions in the proportional-growth asymptotic have been precisely determined for Gaussian matrices \cite{neighborliness}.\footnote{To the best of our knowledge, the average-case analysis techniques of approximate message passing (Donoho, Maleki and Montanari, 2009) cannot be applied
to IHT methods because the hard thresholding operator is not Lipschitz-continuous.}

The main contributions of this paper are as follows:

{\bf 1)\,  We present an entirely new recovery analysis of IHT algorithms.}\quad
In the context of constant stepsize IHT, whereas previous recovery analyses \cite{threshRIP,gargkandekhar,normalized,foucart,HTP} take
 the direct approach of bounding the approximation error from iteration to iteration, we take a two-part approach in which we analyse the fixed points of the algorithm. First, we prove a \textit{stable point condition}, namely a necessary condition for there to be a fixed point on a given support. Second, we give a \textit{convergence condition} which guarantees the convergence of IHT to one of its fixed points. In the case of no noise, this analysis allows us to establish conditions under which, surprisingly,
 IHT converges to its unique fixed point, namely, the original signal $x^*$; noise-dependent recovery results
are also given.
By extending the notion of a fixed point to the (new) concept of an $\au$-\textit{stable point}, we obtain similar recovery results for
the variable-stepsize N-IHT.

{\bf 2)\,  We use average-case assumptions to obtain improved lower bounds on recovery phase transitions for IHT algorithms with Gaussian matrices and Gaussian noise.}\quad
While it is possible to analyse the stable point condition using the RIP \cite{thesis}, we take a different approach. Because the stable point condition has no dependence upon the iterates of the algorithm, it is amenable to analysis for Gaussian matrices under the assumption that the measurement matrix is independent of the signal -- a realistic assumption in CS. We derive precise distributions of this condition's constituent terms, and obtain large deviations bounds on these terms over all possible support sets in the proportional-dimensional asymptotic; in this context, we deduce bounds on some {\it independent RIP} constants that occur naturally in our results. For the convergence condition, we still make use of the RIP, and upper bounds thereon in the proportional-growth asymptotic for Gaussian matrices. However, the RIP condition involved is substantially weaker than any others that have appeared in the literature to date for IHT algorithms. Combining these results, we obtain lower bounds on recovery phase transitions for IHT and N-IHT, namely regions of the phase plane in which stable recovery is guaranteed. In the case of zero noise, we have exact recovery of the original signal. In the case of noise, we derive \textit{stability factors} which bound the approximation error as a multiple of the expectation of the noise. Comparison with state-of-the-art results that have been quantified in the phase transition
framework in \cite{greedy, thesis}, shows a substantial quantitative improvement, both in terms of recovery guarantees as expressed by the height of the
phase transition bounds, and of robustness to noise as expressed by
the size of the stability factors; thus narrowing the gap to observed average-case behaviour. In particular, for the variable-stepsize N-IHT, we obtain
about a factor $10$ improvement in the height of the phase transition bound over best-known results.

We refer to the assumption of independence between signal and measurements as an average-case assumption. The reason for this choice of terminology is that the assumption implies that results hold for given signal instances chosen independently of the measurement matrix, and not for all $k$-sparse signals (as is the case in worst-case RIP analysis). In particular, the independence assumption excludes the (unlikely) scenario in which the worst possible signal is chosen for a given measurement matrix. However, we are not claiming that our analysis is entirely average-case. Though the independence assumption is utilized in analysing the stable point condition, this analysis also involves the use of union bounds, which are often viewed as worst-case techniques. Furthermore, the independence assumption is not used in the convergence analysis, which 
is based entirely upon the worst-case notion of the RIP.
Two comments will help further clarify our contribution. Firstly, our use of union bounds in the stable point analysis is somewhat different from their use in RIP analyses. In our analysis, we fix a signal support set (and coefficients), and obtain probability bounds over all other incorrect support sets. We therefore make use of union bounds within the context of the average-case assumption that the signal is independent of the measurements. Secondly, the RIP conditions we require for convergence are much weaker than existing RIP-based
analyses for IHT algorithms. Thus the improvement we obtain over existing results for IHT algorithms can be attributed jointly to the exploitation of average-case assumptions (in the stable point conditions) --- see the discussion in Section~\ref{largedev} --- and to the weakening of RIP requirements (in the convergence condition) --- see the discussion in Section~\ref{discussion}.
While we achieve significant quantitative improvements over previous analysis~\cite{greedy,thesis}, we emphasize that our results are still lower bounds on the average-case phase transition, and it remains to fully bridge the gap between theory and empirical average-case behaviour.

{\bf 3) We determine a region of phase space within which constant-stepsize IHT is asymptotically guaranteed to have a single fixed point in the case of zero measurement noise.}\quad
Since IHT attempts to solve a nonconvex problem with many local minimizers, it might be natural to expect that the algorithm has a very large number of fixed points. In the noiseless case, our analysis implies a radically new insight: that, within some region of the phase plane (which depends upon the stepsize), IHT has a single fixed point (and hence minimizer), namely the original signal.

{\bf Outline of the paper.}\quad
We begin in Section \ref{algorithms} by describing in more detail the generic IHT algorithm and two stepsize scheme variants, IHT and N-IHT. In Section \ref{general}, we introduce our new recovery analysis, proving our stable point condition, and convergence conditions for both stepsize schemes. Then we focus our attention for the remainder of the paper upon Gaussian matrices: in Section~\ref{phasetrans}, we prove various distributional and large deviations results, and we use these in Section~\ref{phase_proofs} to obtain improved lower bounds on recovery phase transitions in the proportional-growth asymptotic, after which we conclude in Section 6.

{\bf Notation.}\quad
We let $\|\cdot\|$ denote the Euclidean norm.
The support set of the $k$-sparse signal  $x^*$ we aim to recover
will be denoted by
$\mbox{supp}(x^*)=\Lm$ with cardinality $|\Lm|=k$.
Given some index set $\Gm\subseteq\{1,2,\ldots N\}$, we define the complement of $\Gm$ to be $\Gm^C=\{1,2,\ldots
N\}\setminus\Gm$. We
write $x_{\Gm}$ for the restriction of the vector $x$ to the
coefficients indexed by the elements of $\Gm$, and we write $A_{\Gm}$ for the restriction of the matrix $A$ to those columns
indexed by the elements of $\Gm$.

\numsection{Iterative hard thresholding algorithms}
\label{algorithms}

Let us describe in detail the algorithms that are the focus of the analysis in this paper.
Generically, on each hard thresholding iteration $m$, a steepest descent step, possibly with linesearch, is calculated
for the objective $\Psi$ in \req{c:l0problem}, namely, a move is performed from the current iterate
$x^m$ along the negative gradient
of $\Psi$,
\eqn{gradientPsi}{-\nabla \Psi(x^m)=-A^T(Ax^m-b).}
The resulting step is then projected onto the (nonconvex) $l_0$-constraint in \req{c:l0problem} using the so-called {\it hard threshold}
operator ${\cal{H}}_k(\cdot)$ defined as
\[
{\cal{H}}_k(x)\eqdef {\rm arg}\min_{\|z\|_0\leq k} \|z-x\|.
\]
As the name suggests,  ${\cal{H}}_k(\cdot)$ is indeed a thresholding operator, keeping the $k$ largest entries in magnitude of its argument
and setting the rest to zero, namely,
\eqn{Hthreshold}{
{\cal{H}}_k(x) = \left\{
\begin{array}{ll}
x_i &\tim{for}\,\,i\in\Gamma,\\
0 &\tim{for}\,\,i\notin \Gamma,
\end{array}
\right.
\tim{where} \Gamma\eqdef\{\tim{indices of the $k$ largest in magnitude entries of $x$}\}
}
(See \cite[Lemma 1.10]{thesis} for a proof of \req{Hthreshold} given its definition.) To avoid a situation in which the support set  $\Gamma$ is not
uniquely defined, if for instance some of the coefficients are equal in magnitude, then a support set for the identical components
can be selected either randomly or according to some predefined ordering.


The generic IHT algorithm, that includes variants allowing constant or
variable linesearch choices, can be summarized as follows\footnote{The reason we introduce the G-IHT framework
is to allow a more concise presentation of results for the
(fully specified) IHT variants that are of interest.}.
\vspace*{-0.35cm}
\algo{giht}{Generic Iterative Hard Thresholding (G-IHT) algorithm \cite{thresh, normalized}.}{
Given $A$, $b$ and $k$ for problem \req{sparserecovery}, do:
\begin{description}
\item{{\bf Step 0}:} Set $x^0=0$ and $m=0$.\\[1ex]
\hspace*{-0.8cm}While some termination criterion is not satisfied, do:
\item{{\bf Step 1}:} Compute
\eqn{ihtIT}{
x^{m+1}:=H_k\left\{x^m-\al^m A^T(Ax^m-b)\right\},
}
\hspace*{0.6cm}with ${\cal{H}}_k$ defined in \req{Hthreshold} and $\alpha^m>0$
some (pre-defined or computed) stepsize.
\item{{\bf Step 2}:} Set $m=m+1$ and return to Step 1.
\end{description}
}
\vspace*{-0.45cm}
In our analysis, we will consider the possibly infinite sequence of iterates
generated by G-IHT, though in practice a useful termination criterion
such as requiring the residual to be sufficiently small,
would need to be employed.
Two popular stepsize choices will be addressed: {\it constant
stepsize} $\alpha^m=\alpha \in (0,1)$ for all $m$, with the
resulting G-IHT variant being denoted simply as IHT \cite{thresh},
and {\it variable stepsize}
as prescribed in the Normalised IHT  (N-IHT) variant proposed in
\cite{normalized}.

The IHT variant of G-IHT can be summarized as follows.
\vspace*{-0.4cm}
\algo{iht}{Iterative Hard Thresholding (IHT) algorithm \cite{thresh}.}{
Given some $\alpha>0$, on each iteration $m\geq 0$ of G-IHT, do:
\begin{description}
\item In {\bf Step 1}, set $\alpha^m$ in \req{ihtIT} as follows:
\vspace*{-0.25cm}
\eqn{IHTstep}{\alpha^m:=\alpha.}
\vspace*{-0.35cm}
\end{description}}

The N-IHT variant defined below follows \cite{normalized},
having the stepsize $\al^m$  chosen
according to an {\it exact linesearch} \cite{NocedalWright} when the support set of
consecutive iterates stays the same, and using a
{\it shrinking} strategy  when the support set changes so as to ensure
sufficient decrease in the objective of \req{c:l0problem}.
\algo{niht}{Normalised Iterative Hard Thresholding (N-IHT) algorithm \cite{normalized}.}{
Given some $c\in (0,1)$ and $\kappa>1/(1-c)$, on each iteration $m\geq 0$ of G-IHT, do:
\begin{description}
\item In {\bf Step 1}, compute $\alpha^m$ in \req{ihtIT} as follows:
\begin{description}
\item{{\bf Step 1.0}:} Set $\Gm^m:={\rm supp}(x^m)$ and
\eqn{N-IHTstep}{
\alpha^{m}:=
\frac{\|A_{\Gm^m}^T(b-Ax^m)\|^2}{\|A_{\Gm^m}A^T_{\Gm^m}(b-Ax^m)\|^2}.}
Compute $\tilde{x}^{m+1}:=H_k\left\{x^m+\al^{m}
  A^T(b-Ax^m)\right\}$.
If
$\mbox{supp}(\tilde{x}^{m+1})=\Gm^m$,
terminate with $\alpha^m$ given in \req{N-IHTstep}.\\[1ex]
\hspace*{-0.8cm} While $\al^m  \geq
(1-c)\displaystyle\frac{\|\tilde{x}^{m+1}-x^m\|^2}{\|A(\tilde{x}^{m+1}-x^m)\|^2}$, do:

\item{{\bf Step 1.1}:} $\al^m:=\al^m/[\kappa(1-c)]$;

\item{{\bf Step 1.2}:} $\tilde{x}^{m+1}:=H_k\left\{x^m+\al^m
    A^T(b-Ax^m)\right\}$;\\[1ex]
\hspace*{-0.8cm} End.
\end{description}
\end{description}
}

Under the (weakest) assumptions of this paper, we can  ensure that both the exact linesearch and the shrinkage stepsizes in N-IHT
are well-defined, until termination;
see Section \ref{NIHTconvergencesection}.
The shrinkage iteration between Steps 1.1--1.2 of N-IHT can be shown to terminate in finitely many steps \cite{normalized}.

\numsection{Deterministic conditions for a recovery analysis}\label{general}

In this section we derive conditions for IHT algorithms when
applied to general measurement
matrices $A$. Hence we
make the following common assumption for compressed sensing algorithms.

\vspace*{0.1cm}
\ass{genposass}{\quad  The matrix $A$ is in $2k$-general position, namely any
$2k$ of its columns are linearly independent.}
\vspace*{0.1cm}
 The (weak) assumption {\bf A.1} is equivalent to the condition that, for
any $\Gm$ such that $|\Gm|=2k$, the matrix $A^T_{\Gm}A_{\Gm}$
is nonsingular. Thus, whenever {\bf A.1} holds and there is no noise in the system (i.e., $e=0$ in \req{sparserecovery}), we have
$\|A(x^*-x)\|>0$ for any $k$-sparse $x\neq x^*$, and so
$x^*$ is the unique $k$-sparse exact solution to the
linear system $b=Ax^*$. Note also that {\bf A.1}  holds if $A$ is in general position. It is also satisfied with probability $1$ if
$A$ is a Gaussian matrix and $2k\leq n$ \cite{mvd}.


The results derived in this section come in two parts: a necessary condition for the existence of (generalized) fixed points of G-IHT
and a sufficient condition guaranteeing convergence for particular stepsize schemes. In the deterministic noiseless case, the
former condition can be used to guarantee the existence of at most one fixed point, namely, the original signal $x^*$; thus,
provided we also have convergence of the algorithm to \textit{some} such fixed/stable point, signal recovery is ensured.
The below deterministic results also yield similar recovery properties (based on proximity/closeness
of fixed points to the original signal) in the presence of noise and Gaussian measurement matrices,
as we show in later sections.

\subsection{A stable-point condition}
\label{stablesection}

We introduce the concept of an $\au$\textit{-stable point} of
G-IHT, a generalization of fixed points.

\begin{definition}[\textbf{$\underline{\alpha}$-stable points of G-IHT}]\label{stable}
Given $\au>0$ and an index set $\Gm$ with $|\Gm|=k$, we say $\xb\in\RR^N$ is an $\au$-stable point of  G-IHT on $\Gm$ if $\mbox{supp}(\xb) \subseteq  \Gm$ and
\begin{equation}\label{stable1}
\left\{A^T(b-A\xb)\right\}_{\Gm}=0\;\;\;\;\mbox{and}
\end{equation}
\begin{equation}\label{stable2}
\min_{i\in\Gm}|\xb_i|\geq
  \au\displaystyle\max_{j\in\Gm^C}|\left\{A^T(b-A\xb)\right\}_j|.
\end{equation}
\end{definition}

Note that in the noiseless case ($e=0$ in \req{sparserecovery}), the original signal $x^*$
is clearly an $\au$-stable point on ${\rm supp}(x)=\Lm$, for any value of $\au>0$.

In the case of the constant-stepsize IHT algorithm, an $\al$-stable point is nothing other than a \textit{fixed point}
of IHT  (see Blumensath \& Davies~\cite[Lemma
  6]{thresh}) or an $L$-stationary point of \req{c:l0problem} in \cite[\S2.3]{beck}.
(Indeed, if a further IHT iteration is applied at a fixed point $\xb$, there is no change in the support set;
thus the gradient term on the complement of the support of $\xb$ must be suitably small, which is \req{stable2}. Also, the coefficients on
the support of $\xb$ must remain unchanged, and so we require the gradient on the support of $\xb$ to be zero, namely \req{stable1}.)
A generalization of the notion
of a fixed point and $L$-stationary point to stable points is, however, required to allow for variable stepsize
schemes in G-IHT\footnote{When \req{stable1} holds for N-IHT, the exact linesearch stepsize \req{N-IHTstep} is not well-defined with its numerator and denominator both being zero.}; we will be interested in values of $\au$ that lower bound the
stepsize $\al^m$ of G-IHT.

Next
we show that any $\au$-stable point is a
minimum-norm solution on some $k$-subspace.
\begin{lemma}\label{necfp} Let {\bf A.1} hold and
 $\xb$ be an $\au$-stable point of G-IHT on $\Gm$ for some $\au>0$. Then
\begin{equation}\label{exactk}
\xb_{\Gm}=A_{\Gm}^{\dag}b,
\end{equation}
where $A_{\Gm}^{\dag}$ is the Moore-Penrose pseudo-inverse, namely,
\eqn{pseudoinv}{A_{\Gm}^{\dag}\eqdef (A_{\Gm}^T A_{\Gm})^{-1}A_{\Gm}^T.}
\end{lemma}

\proof{
It follows from (\ref{stable1}) that
 $A^T_{\Gm}(b-A_{\Gm}\xb_{\Gm})=0$. By {\bf A.1},
 the pseudoinverse $A^{\dag}_{\Gm}$ is well-defined and we may rearrange to
 give (\ref{exactk}).}

While the previous lemma tells us that any stable point is necessarily a minimum-norm
 solution on some $k$-subspace, the converse may not hold. Next, we give a more useful necessary condition for there to exist
a stable point on a given support set. We will use the latter condition  in a sufficient sense later on, to guarantee that under certain conditions,
all G-IHT stable points are {\it close} to the underlying signal,
which in the noiseless case reduces to G-IHT having \textit{at most
 one} stable point, namely, the original signal.


\begin{theorem}[\textbf{Stable point condition; noise case}]\label{singlefp}
Consider problem \req{sparserecovery} and let $\Lm={\rm supp}(x^*)$.  Suppose {\bf A.1}
 holds\footnote{Assumption {\bf A.1} may in fact be weakened in Theorem~\ref{singlefp} to requiring the matrix $A$ to be in $k$-general position.} and suppose there exists an $\au$-stable point of G-IHT on some $\Gm$ such that $\Gm\neq\Lm$. Then
\begin{equation}\label{stablecond_noise}
\left\|A_{\Gm}^{\dag}A_{\Lm\sm\Gm}x^*_{\Lm\sm\Gm}\right\|+\left\|A_{\Gm}^{\dag}e\right\|\geq\au\left\{\left\|A_{\Lm\sm\Gm}^T(I-A_{\Gm}A^{\dag}_{\Gm})A_{\Lm\sm\Gm}x^*_{\Lm\sm\Gm}\right\|-\left\|A_{\Lm\sm\Gm}^T(I-A_{\Gm}A^{\dag}_{\Gm})e\right\|\right\},
\end{equation}
where $A_{\Gm}^{\dag}$ is defined in \req{pseudoinv}.
\end{theorem}

\proof{Assume $\xb$ is an $\au$-stable point on $\Gm$.
Since $\Gm\sm\Lm\subseteq\Gm$ and
$\Lm\sm\Gm\subseteq\Gm^C$,  where $\Lm={\rm supp}(x^*)$,
(\ref{stable2}) implies that
\begin{equation}\label{fpminmaxl}\displaystyle\min_{i\in\Gm\sm\Lm}|\xb_i|\geq
  \au\displaystyle\max_{j\in\Lm\sm\Gm}|\left\{A^T(b-A\xb)\right\}_j|.
\end{equation}
Definition \ref{stable} implies that $|\Gm|=|\Lm|$, and so $|\Gamma\sm\Lm|=|\Lm\sm\Gm|$.
This, properties of the Euclidean norm and \req{fpminmaxl} provide
\eqn{2norm}{
\|\xb_{\Gm\sm\Lm}\|^2\geq |\Gamma\sm\Lm|\left\{\min_{i\in\Gm\sm\Lm}|\xb_i|\right\}^2
\geq |\Lm\sm\Gamma|\left\{\au\max_{j\in\Lm\sm\Gm}|\left\{A^T(b-A\xb)\right\}_j|\right\}^2
\geq \au^2\|A^T_{\Lm\sm\Gm}(b-A\xb)\|^2.
}
Problem \req{sparserecovery}  and $x_{\Lm^C}^*=0$ imply
\eqn{b_expand}{b=Ax^*+e=A_{\Gm}x^*_{\Gm}+A_{\Lm\sm\Gm}x^*_{\Lm\sm\Gm}+e.}
This and Lemma~\ref{necfp} now provide, under {\bf A.1},
$$\xb_{\Gm}=A_{\Gm}^{\dag}b=x^*_{\Gm}+A_{\Gm}^{\dag}A_{\Lm\sm\Gm}x^*_{\Lm\sm\Gm}+A_{\Gm}^{\dag} e,$$
where in the last equality, we used $A_{\Gm}^{\dag}A_{\Gm}=I$. Therefore, since $x^*_{\Gm\sm\Lm}=0$, we deduce
$$
\xb_{\Gm\sm\Lm}=\left(A_{\Gm}^{\dag}A_{\Lm\sm\Gm}x^*_{\Lm\sm\Gm}+A_{\Gm}^{\dag} e\right)_{\Gm\sm\Lm}
$$
and so,
\begin{equation}\label{lhsbound}
\|\xb_{\Gm\sm\Lm}\|\le\left\|\left(A^{\dag}_{\Gm}A_{\Lm\sm\Gm}x^*_{\Lm\sm\Gm}\right)_{\Gm\sm\Lm}\right\|+\left\|\left(A^{\dag}_{\Gm}e\right)_{\Gm\sm\Lm}\right\| \le\left\|A^{\dag}_{\Gm}A_{\Lm\sm\Gm}x^*_{\Lm\sm\Gm}\right\|+\left\|A^{\dag}_{\Gm}e\right\|,
\end{equation}
which upper bounds the left-hand side of (\ref{2norm}).
Under {\bf A.1}, we may next use Lemma~\ref{necfp} and (\ref{b_expand}) to express the right-hand side of (\ref{2norm}) independently of $\xb$, as
$$
\begin{array}{l}
A^T_{\Lm\sm\Gm}(b-A\xb)=A^T_{\Lm\sm\Gm}(I-A_{\Gm}A^{\dag}_{\Gm})b=A_{\Lm\sm\Gm}^T(I-A_{\Gm}A^{\dag}_{\Gm})(A_{\Lm\sm\Gm}x^*_{\Lm\sm\Gm}+e),
\end{array}
$$
where in the last equality, we used $A_{\Gm}^{\dag}A_{\Gm}=I$. We therefore may deduce
\begin{equation}\label{rhsbound}
\left\|A^T_{\Lm\sm\Gm}(b-A\xb)\right\|\geq\left\|A_{\Lm\sm\Gm}^T(I-A_{\Gm}A^{\dag}_{\Gm})A_{\Lm\sm\Gm}x^*_{\Lm\sm\Gm}\right\|-\left\|A_{\Lm\sm\Gm}^T(I-A_{\Gm}A^{\dag}_{\Gm})e\right\|.
\end{equation}
Substituting (\ref{lhsbound}) and (\ref{rhsbound}) into (\ref{2norm}), we arrive at \req{stablecond_noise}.}

Theorem~\ref{singlefp} simplifies further in the noiseless case.

\begin{corollary}[\textbf{Stable point condition; noiseless case}]\label{stable_noiseless}
Consider problem \req{sparserecovery} with $e\eqdef 0$ and let $\Lm={\rm supp}(x^*)$.  Suppose {\bf A.1}
 holds and suppose there exists an $\au$-stable point of G-IHT on some $\Gm$ such that $\Gm\neq\Lm$. Then
\begin{equation}\label{fpcondNEW}
\left\|A_{\Gm}^{\dag}A_{\Lm\sm\Gm}x^*_{\Lm\sm\Gm}\right\|\geq\au\left\|A_{\Lm\sm\Gm}^T(I-A_{\Gm}A^{\dag}_{\Gm})A_{\Lm\sm\Gm}x^*_{\Lm\sm\Gm}\right\|,
\end{equation}
where $A_{\Gm}^{\dag}$ is defined in \req{pseudoinv}.
\end{corollary}
\proof{The result follows immediately by setting $e\eqdef 0$ in (\ref{stablecond_noise}).}

Clearly, Corollary \ref{stable_noiseless} implies that if the reverse inequality in \req{fpcondNEW} holds for
all support sets $\Gm \neq \Lm$,  then $x^*$ is the only $\au$-stable point of G-IHT.

\subsection{A convergence condition}
\label{convergencesection}

This section gives conditions for IHT algorithms to convergence to stable points.
Recalling \req{c:l0problem} and \req{ihtIT}, we introduce the notation
\eqn{Gmgm}{g^m\eqdef \nabla \Psi(x^m) \tim{and} \Gm^m\eqdef {\rm supp}(x^m), \tim{for all $m$.}}
Some useful properties of the G-IHT iterates are given
in the next lemma.
\begin{lemma}\label{iterprop}
Apply the G-IHT algorithm to solve \req{c:l0problem}. Then the G-IHT iterates satisfy for all $m\geq 0$,
\begin{equation}\label{gradresult1}
\|x^{m+1}-x^m\|^2+2\al^m (g^m)^T(x^{m+1}-x^m)\le 0
\end{equation}
and
\begin{equation}\label{taylorgen}
\Psi(x^{m+1})-\Psi(x^m)= (g^m)^T(x^{m+1}-x^m)+\frac{1}{2}\left\|A(x^{m+1}-x^m)\right\|^2.
\end{equation}
\end{lemma}
\proof{Since the hard thresholding operation in \req{ihtIT}
 can be viewed as a projection onto the constraint of \req{c:l0problem}, we may rewrite the G-IHT iteration \req{ihtIT} as
$$
x^{m+1}=\displaystyle\arg\min_{\|z\|_0\le k}\|z-\left\{x^m-\al^m g^m\right\}\|^2.
$$
This further gives that
$$
\|x^{m+1}-(x^m-\al^m g^m)\|^2\le\|x^m-(x^m-\al^m g^m)\|^2=(\al^m)^2\|g^m\|^2,
$$
which expands to give
$\|x^{m+1}-x^m\|^2+2\al^m
(g^m)^T(x^{m+1}-x^m)+(\al^m)^2\|g^m\|^2\le(\al^m)^2\|g^m\|^2$,
and so  \req{gradresult1} holds. Since $\Psi$ in \req{c:l0problem} is quadratic,
we have no remainder in the following second-order Taylor expansion
$$
\begin{array}{lcl}
\Psi(x^{m+1})-\Psi(x^m)&=&\left[\nabla\Psi(x^m)\right]^T(x^{m+1}-x^m)+\half (x^{m+1}-x^m)^T\left[\nabla^2\Psi\right](x^{m+1}-x^m)\\
&=&(g^m)^T(x^{m+1}-x^m)+\half (x^{m+1}-x^m)^T A^T
  A(x^{m+1}-x^m),
\end{array}
$$
and so \req{taylorgen} follows.}


 A sufficient condition for G-IHT convergence is given next.
\begin{lemma}\label{convlem}
Let {\bf A.1} hold and the G-IHT iterates satisfy
\begin{equation}\label{convsuff}
\|x^{m+1}-x^m\|^2\le d\left[\Psi(x^m)-\Psi(x^{m+1})\right], \tim{for
  all $m\geq 0$,}
\end{equation}
for some $d>0$ independent of $m$ and
$\Psi$ defined in \req{c:l0problem}. Assume that there exist
$\overline{\alpha}\geq \au>0$
such
\eqn{lbalpham}{\overline{\alpha}\geq \alpha^m\geq \au \tim{for all $m\geq 0$.}}
Then  $x^m\rightarrow\xb$ as
$m\rightarrow \infty$, where $\xb$ is an $\au$-stable point of G-IHT.
\end{lemma}

\proof{We deduce from (\ref{convsuff}) that
$$
\sum_{m=0}^{\infty} \|x^{m+1}-x^m\|^2 \leq d\sum_{m=0}^{\infty}\left[\Psi(x^m)-\Psi(x^{m+1})\right]\leq d\Psi(x^0),
$$
where to obtain the last inequality, we used  $\Psi(x^m)\geq
0$. Thus convergent series properties provide
\begin{equation}\label{getclosergen}
\|x^{m+1}-x^m\|\longrightarrow 0 \quad \mbox{as}\quad m\longrightarrow\infty.
\end{equation}
From \req{ihtIT} and \req{Gmgm}, we deduce
$$
x_{\Gm^{m+1}}^{m+1}=x_{\Gm^{m+1}}^m-\al^m  g^m_{\Gm^{m+1}} \quad
{\rm and}\quad
x_{(\Gm^{m+1})^C}^{m+1}=0.
$$
Thus restricting (\ref{getclosergen}) to
$\Gm^{m+1}$ and using \req{lbalpham} provide
\begin{equation}\label{gradplusone}
\|g_{\Gm^{m+1}}^m\|\longrightarrow 0 \quad \mbox{as}\quad m\longrightarrow\infty,
\end{equation}
while restricting (\ref{getclosergen})  to $\Gm^m\sm\Gm^{m+1}$ yields
\begin{equation}\label{rejectedx}
\|x_{\Gm^m\sm\Gm^{m+1}}^m\|\longrightarrow 0.
\end{equation}
For $m\geq 0$, let $y^{m}$ denote the minimum-norm solution on $\Gm^m$, namely,
\eqn{minnormsol}{
y^{m}_{\Gm^m}\eqdef
A_{\Gm^m}^{\dag}b\tim{and}y^{m}_{(\Gm^m)^C}\eqdef 0,
}
which is well-defined due to {\bf A.1}. Then \req{minnormsol} and
$x^m_{({\Gm^m})^C}=0$ provide
$$\begin{array}{lcl}
\|y^{m+1}-x^m\|&\le&\|y^{m+1}_{\Gm^{m+1}}-x^m_{\Gm^{m+1}}\|+\|x^m_{(\Gm^{m+1})^C}\|
=\|A_{\Gm^{m+1}}^{\dag}b-x^m_{\Gm^{m+1}}\|+\|x_{\Gm^m\sm\Gm^{m+1}}^m\|\\
&=&\|(A_{\Gm^{m+1}}^T A_{\Gm^{m+1}})^{-1}A_{\Gm^{m+1}}^T(b-A_{\Gm^{m+1}}x^m_{\Gm^{m+1}})\|+\|x_{\Gm^m\sm\Gm^{m+1}}^m\|\\
&=&\|(A_{\Gm^{m+1}}^T A_{\Gm^{m+1}})^{-1}g^m_{\Gm^{m+1}}\|+\|x_{\Gm^m\sm\Gm^{m+1}}^m\|
\longrightarrow 0, \tim{as $m\longrightarrow \infty$,}
\end{array}$$
where the limit follows from (\ref{gradplusone}), (\ref{rejectedx}),
{\bf A.1}
and the fact that there are finitely many distinct support sets
$\Gm^m$, $m\geq 0$. This and \req{getclosergen} further give
\eqn{ymxmlim}{
\|y^m-x^m\|\longrightarrow 0 \tim{as $m\longrightarrow \infty$,}
}
and so  for any $\epsilon>0$, there exists $m_0\geq 0$ such that
\begin{equation}\label{c:closeminnorm}
\|y^{m}-x^m\|\leq \epsilon, \tim{for all $m\geq m_0$.}
\end{equation}
We denote the index set of changing minimal-norm solutions by
$$\calS\eqdef \left\{m\geq m_0:\, y^{m+1}\neq y^{m}\right\},$$
and we will show that $\calS$ is finite. Define
\begin{equation}\label{c:farminnorm}
\epsilon\eqdef \fourth\min_{m\in
  \calS}\|y^{m+1}-y^{m}\|.
\end{equation}
Note that $\epsilon>0$ since there are finitely many distinct support
sets $\Gm^m$, $m\geq 0$.
Then, the triangle inequality,
(\ref{c:closeminnorm}) and (\ref{c:farminnorm})  yield
$$
\|x^{m+1}-x^m\|\geq
\|y^{m+1}-y^m\|-\|y^{m+1}-x^{m+1}\|-\|y^m-x^m\|\geq 4\epsilon
-\epsilon-\epsilon >\epsilon, \tim{for all $m\in\calS$.}
$$
This and (\ref{getclosergen}) imply that $\calS$ must be finite and so
there exists $m_1\geq m_0$ such that $y^{m+1}=y^m=\xb$ for all $m\geq
m_1$, where $\xb_{\Gm}=A_{\Gm}^{\dag}b$ and $\xb_{\Gm^C}=0$, for some
$\Gm$ with $|\Gm|=k$. This and \req{ymxmlim} give
\eqn{convxm}{
x^m\longrightarrow\xb, \tim{as $m\longrightarrow \infty$.}}
Clearly,
 (\ref{stable1}) holds for the limit point $\xb$ of the iterates $\{x^m\}$.
To complete the proof, it remains to establish \req{stable2}.
The thresholding operation that defines $x^{m+1}$ in G-IHT gives that
\begin{equation}\label{optimalchoice}
\min_{i\in\Gm^{m+1}}|x_i^{m+1}|\geq
  \displaystyle\max_{j\in(\Gm^{m+1})^C}|\{x^m-\al^m g^m\}_j|,\tim{for all $m\geq 0$,}
\end{equation}
and \req{lbalpham} implies that there exists a convergent subsequence of stepsizes,
\eqn{alphamconvsubseq}{\alpha^{m_r}\longrightarrow \tilde{\alpha}\geq \underline{\alpha} \tim{as $r\longrightarrow \infty$.}}
Letting $\epsilon\eqdef\half \min_{i\in{\rm supp}(\xb)} \xb_i$,
\req{convxm} implies that $\|x^m-\xb\|\leq \epsilon$, and so
\eqn{suppinclusion}{
{\rm supp}(\xb)\subseteq \Gm^m, \tim{for all
$m$ sufficiently large.}}
Firstly, assume that ${\rm supp}(\xb)=\Gm$. Then, since $|\Gm|=|\Gm^m|=k$, \req{suppinclusion} implies
that $\Gm^m=\Gm$ for all $m$ sufficiently large, which together with \req{optimalchoice}, provides
\begin{equation}\label{optimalchoice1}
\min_{i\in\Gm}|x_i^{m+1}|\geq
  \displaystyle\max_{j\in \Gm^C}|\{x^m-\al^m g^m\}_j|,\tim{for all $m$ sufficiently large.}
\end{equation}
Passing to the limit in \req{optimalchoice1} on the subsequence $m_r$ for which
\req{alphamconvsubseq} holds,  using \req{convxm}, $\xb_{\Gm^C}=0$ and the right-hand side of \req{lbalpham}
imply \req{stable2} holds in this case. It remains to consider the case when ${\rm supp}(\xb)\subset \Gm$. Then
$\min_{i\in\Gm}|\xb_i|=0$ and so \req{convxm} further provides
\eqn{m+1lim}{\min_{i\in \Gm^{m+1}} |x_i^{m+1}|\longrightarrow 0 \tim{as $m\longrightarrow \infty$.}}
Now \req{suppinclusion} and again \req{convxm} provide
\eqn{m+2lim}{x^m_{\Gm^{m+1}}\longrightarrow 0 \tim{as $m\longrightarrow \infty$.}}
Passing to the limit in  \req{optimalchoice} on the subsequence $m_r$ for which
\req{alphamconvsubseq} holds, and using \req{m+1lim} and \req{m+2lim}, we obtain
that $g^m_{(\Gm^{m+1})^C}\longrightarrow 0$ as $m\longrightarrow \infty$. This and \req{gradplusone}
now give that $g^m= A^T(Ax^m-b)\longrightarrow 0$, which
due to \req{convxm}, implies that $A^T(b-A\xb)=0$ and so \req{stable2} trivially holds in this case.}

In order to ensure \req{convsuff} and \req{lbalpham}, we make use of
the well-known Restricted Isometry Property (RIP) constants of the matrix $A$,
defined as follows.

\begin{definition}\label{RIP} \cite{Candes, lqphase}
Define $L_s$ and $U_s$, the {\it lower} and {\it upper} RIP constants
of $A$ of order $s$, to be, respectively,
\eqn{RIC}{
L_s=1-\displaystyle\min_{1\le\|y\|_0\le
  s}\displaystyle\frac{\|Ay\|^2}{\|y\|^2} \quad {\tim{and}}\quad
U_s=\displaystyle\max_{1\le\|y\|_0\le s}\displaystyle\frac{\|Ay\|^2}{\|y\|^2}-1.}
\end{definition}

Note that {\bf A.1} is equivalent to the requirement that
$L_{2k}<1$.  This and all other RIP conditions we use here are substantially weaker than
those employed in existing worst-case analyses for IHT algorithms.


In order to ensure \req{convsuff} and \req{lbalpham} -- using RIP constants or
otherwise -- we must specify the choice of stepsize $\alpha^m$ in G-IHT.
(This is by contrast to the stable point condition for
which only lower bounds on the stepsizes $\alpha^m$ matter.)
Hence we now return to the constant-stepsize IHT and variable-stepsize N-IHT variants defined in Section \ref{algorithms}.

\subsubsection{A convergence condition for the IHT algorithm}

 In~\cite{thresh}, Blumensath and Davies prove convergence of IHT
iterates to a fixed point that is also a local minimizer of
(\ref{c:l0problem}) (that may or may not be the original signal
$x^*$) under the assumption that $\alpha\|A\|_2<1$. Similarly, Beck
and Eldar
\cite[Theorem 3.2]{beck} show IHT iterates converge to an
$L$-stationary point, an equivalent notion to that of a fixed
point, under a commensurate condition on the stepsize, namely,
$\alpha\|A^TA\|<1$. Largely following the method of proof in
\cite{thresh}, we now show that the requirement on the IHT stepsize
in both these analyses can be weakened to a condition involving the
RIP constant $U_{2k}$ of $A$.

\begin{theorem}\label{conv1}
Suppose that {\bf A.1} holds, and that the IHT stepsize
is chosen to satisfy
\eqn{RIC:IHTstepsize}{\alpha<\frac{1}{1+U_{2k}}.}
Then the IHT iterates $\{x^m\}$ converge to
an $\al$-stable point $\xb$ of IHT.
\end{theorem}

\proof{Let $m\geq 0$. Since the support size of the change to the iterates $x^{m+1}-x^m$ is at most $2k$,
the upper RIP of $A$ in \req{RIC} with $s=2k$   provides
$\|A(x^{m+1}-x^m)\|^2\le(1+U_{2k})\|x^{m+1}-x^m\|^2$.
Using this bound, and  (\ref{gradresult1}) with the choice \req{IHTstep}, in (\ref{taylorgen}), we obtain
$$
\Psi(x^{m+1})-\Psi(x^m)\le-\frac{1}{2\al}\|x^{m+1}-x^m\|^2+\frac{1}{2}(1+U_{2k})\|x^{m+1}-x^m\|^2
=\displaystyle\frac{\al(1+U_{2k})-1}{2\al}\|x^{m+1}-x^m\|^2,
$$
which due to \req{RIC:IHTstepsize}, implies that \req{convsuff} holds with  $d\eqdef 2\al/[1-\al(1+U_{2k})]$.
Due to \req{IHTstep}, \req{lbalpham} trivially holds with $\overline{\alpha}=\au=\al$. Thus  Lemma~\ref{convlem}
applies, and so the IHT iterates  $x^m$ converge   to an $\al$-stable point of IHT.}

\subsubsection{A convergence condition for the N-IHT algorithm}
\label{NIHTconvergencesection}

Using the notation \req{RIC} and {\bf A.1}, we obtain that the N-IHT stepsize $\alpha^m$ satisfies
\begin{equation}\label{NIHTbound}
\frac{1}{1+U_k}\le\al^m\le\frac{1}{1-L_k}, \tim{whenever $\alpha^m$ satisfies \req{N-IHTstep},}
\end{equation}
and using also \req{ihtIT}, that
\begin{equation}\label{NIHTbound1}
\frac{1}{\kappa(1+U_{2k})}\le\al^m\le\frac{1-c}{1-L_{2k}}, \tim{otherwise (i.e., whenever $\alpha^m$ is shrunk according to Steps 1.1--1.2).}
\end{equation}
As the RIP constants of $A$ are monotonically increasing with $k$ and $\kappa, \,c \in (0,1)$, \req{NIHTbound} and \req{NIHTbound1} imply
\eqn{finalLB}{\frac{1}{\kappa(1+U_{2k})}\leq \alpha^m \leq \frac{1-c}{1-L_{2k}},\tim{for all $m\geq 0$.}}
\begin{theorem}\label{conv3}
Suppose {\bf A.1} holds. Then
the N-IHT iterates $\{x^m\}$ converge to a $[\kappa(1+U_{2k})]^{-1}$-stable point $\xb$ of N-IHT.
\end{theorem}

\proof{Firstly, we consider the case when $\alpha^m$ satisfies \req{N-IHTstep}. Then \req{Gmgm} implies
$\Gamma^{m+1}=\Gamma^m$, and \req{ihtIT} implies
\begin{equation}\label{map2}
x_{\Gm^m}^{m+1}=x_{\Gm^m}^m-\al^m g^m_{\Gm^m}.
\end{equation}
Using \req{map2}, \req{N-IHTstep} becomes
\begin{equation}\label{exactrule}
\al^m=\frac{\|g_{\Gm^m}^m\|^2}{\|A_{\Gm^m}g_{\Gm^m}^m\|^2}=\frac{\|x^{m+1}-x^m\|^2}{\|A(x^{m+1}-x^m)\|^2}.
\end{equation}
Using that $x^{m+1}-x^m$ is supported on $\Gm^m$,  expressing $g^m_{\Gm^m}$ from \req{map2} and substituting
 into \req{taylorgen}, we deduce that
\eqn{niht1}{
\begin{array}{lcl}
\Psi(x^{m+1})-\Psi(x^m)&=&
-\frac{1}{\al^m}(x_{\Gm^m}^{m+1}-x_{\Gm^m}^m)^T(x_{\Gm}^{m+1}-x_{\Gm}^m)+\frac{1}{2}\|A(x^{m+1}-x^m)\|^2\\[1ex]
&=&-\frac{1}{\al^m}\|x^{m+1}-x^m\|^2+ \frac{1}{2\al^m}\|x^{m+1}-x^m\|^2=-\frac{1}{2\al^m}\|x^{m+1}-x^m\|^2,
\end{array}}
where to obtain the second equality, we also used \req{exactrule}.
Alternatively, when $\al^m$ is computed by shrinkage,  we deduce
that
$$
\|A(x^{m+1}-x^m)\|^2\leq \frac{1-c}{2\al^m} \|x^{m+1}-x^m\|^2.
$$
Substituting this and \req{gradresult1} into \req{taylorgen}, we obtain
\eqn{niht2}{
\Psi(x^{m+1})-\Psi(x^m)\leq -\frac{1}{2\alpha^m}\|x^{m+1}-x^m\|^2 + \frac{1-c}{2\al^m} \|x^{m+1}-x^m\|^2= -\frac{c}{2\al^m} \|x^{m+1}-x^m\|^2.}
Thus \req{niht1}, \req{niht2} and $c\in (0,1)$ imply that for all $m\geq 0$,
$$
\|x^{m+1}-x^m\|^2 \leq  \frac{2\alpha^m}{c} [\Psi(x^m)-\Psi(x^{m+1})] \leq \frac{2(1-c)}{c(1-L_{2k})}  [\Psi(x^m)-\Psi(x^{m+1})],
$$
due to \req{finalLB}. Hence \req{convsuff} holds
with $d\eqdef 2 (1-c)/[c(1-L_{2k})]$, and so does \req{lbalpham} due to \req{finalLB}. Lemma \ref{convlem} applies and together with \req{finalLB}
provides the required conclusion.}

Note that due to \req{finalLB}, the shrinkage strategy, rather than the exact linesearch, determines the value of
$\au$ in Theorem \ref{conv3}, which is crucial for our phase transition bounds. However, we cannot guarantee that
the less-conservative exact linesearch strategy is taken asymptotically.

\subsection{Deterministic recovery conditions}

In the case of zero measurement noise, combining the two parts of our analysis in Sections~\ref{stablesection} and~\ref{convergencesection} respectively leads immediately to recovery conditions for both IHT and N-IHT.

\begin{theorem}\label{IHT_det_recov}
Consider problem \req{sparserecovery} with $e\eqdef 0$ and let $\Lm={\rm supp}(x^*)$.  Suppose that {\bf A.1}
 holds, that the stepsize $\alpha$ of IHT satisfies (\ref{RIC:IHTstepsize}), and  that
\begin{equation}\label{fpcond1}
\left\|A_{\Gm}^{\dag}A_{\Lm\sm\Gm}x^*_{\Lm\sm\Gm}\right\|<\al\left\|A_{\Lm\sm\Gm}^T(I-A_{\Gm}A^{\dag}_{\Gm})A_{\Lm\sm\Gm}x^*_{\Lm\sm\Gm}\right\|
\end{equation}
for all $\Gm\neq\Lm$ such that $|\Gm|=k$, where $A_{\Gm}^{\dag}$ is defined in \req{pseudoinv}. Then the IHT iterates $\{x^m\}$ converge to its only fixed point, namely, the
original signal $x^*$.
\end{theorem}

\proof{Under Assumption {\bf A.1}, Corollary~\ref{stable_noiseless} and (\ref{fpcond1}) imply that there exists no $\al$-stable point on any $\Gm\neq\Lm$ such that $|\Gm|=k$. It follows that any $\al$-stable point is supported on $\Lm$, and therefore by Lemma~\ref{necfp}, it must coincide with $x^*$. Also under Assumption {\bf A.1}, it follows from (\ref{RIC:IHTstepsize}) and Theorem~\ref{conv1} that IHT converges to an $\al$-stable point, and hence to $x^*$. Since a fixed point of IHT with stepsize $\alpha$ is the same as an $\alpha-$stable point, we
conclude the proof.}

\begin{theorem}\label{NIHT_det_recov}
Consider problem \req{sparserecovery} with $e\eqdef 0$ and let $\Lm={\rm supp}(x^*)$.  Suppose that {\bf A.1}
 holds and  that
\begin{equation}\label{fpcond}
\left\|A_{\Gm}^{\dag}A_{\Lm\sm\Gm}x^*_{\Lm\sm\Gm}\right\|<[\kappa(1+U_{2k})]^{-1}\left\|A_{\Lm\sm\Gm}^T(I-A_{\Gm}A^{\dag}_{\Gm})A_{\Lm\sm\Gm}x^*_{\Lm\sm\Gm}\right\|
\end{equation}
for all $\Gm\neq\Lm$ such that $|\Gm|=k$, where $A_{\Gm}^{\dag}$ is defined in \req{pseudoinv}. Then the N-IHT iterates $\{x^m\}$ converge to the original signal $x^*$.
\end{theorem}

\proof{Under Assumption {\bf A.1}, Corollary~\ref{stable_noiseless} and (\ref{fpcond}) imply that there exists no $[\kappa(1+U_{2k})]^{-1}$-stable point on any $\Gm\neq\Lm$ such that $|\Gm|=k$. It follows that any $[\kappa(1+U_{2k})]^{-1}$-stable point is supported on $\Lm$, and therefore by Lemma~\ref{necfp} it must be $x^*$. Also under Assumption {\bf A.1}, Theorem~\ref{conv3} implies that we also have convergence to an $[\kappa(1+U_{2k})]^{-1}$-stable point, which concludes the proof.}

While Theorems~\ref{IHT_det_recov} and~\ref{NIHT_det_recov} give conditions guaranteeing recovery, what is less clear is when one might expect these conditions to be satisfied. We provide answers to this question in the rest of the paper, quantifying when these conditions are satisfied in the case of Gaussian matrices. Furthermore, we also extend our analysis for Gaussian matrices to the case of measurements contaminated by Gaussian noise.

\numsection{Probabilistic quantification of the deterministic analysis}
\label{phasetrans}

{\bf Brief roadmap for Sections \ref{phasetrans} and \ref{phase_proofs}.}
For the remainder of the paper, we focus our attention on quantifying the deterministic recovery conditions of Section~\ref{general} in the case of Gaussian matrices and discussing our results. In the case of IHT, two conditions must be satisfied to ensure recovery: an RIP-based convergence condition (\ref{RIC:IHTstepsize}) and the stable point condition (\ref{stablecond_noise}). For N-IHT, the two notions combine to give a single condition (this condition was given in Theorem~\ref{NIHT_det_recov} in the case of zero measurement noise; a corresponding condition will be obtained in the case of nonzero noise in Section~\ref{NIHT_proofs}). The quantification of these conditions for Gaussian matrices will be performed in the proportional-growth asymptotic of Definition \ref{propdimdef}. As we explain below, there is a need to quantify each condition for a given support, as well as the union/intersection of these conditions over all possible support sets.

Such a quantification has already been done for the RIP constants of Gaussian matrices. Namely, 
it was shown in~\cite{lqphase} that bounds on RIP constants of Gaussian matrices can be obtained in the proportional-growth asymptotic; a subsequent improvement on these bounds was obtained in~\cite{BT}.
\begin{lemma}[\textbf{Gaussian RIP bounds~\cite[Theorem 2.3]{BT}}]\label{RIP_bounds}
Suppose $A\sim\mathcal{N}_{n,N}(0,1/n)$ has RIP constants $L_{k}$ and $U_{k}$ as defined in Definition~\ref{RIP}, and let the implicit but computable expressions $\LL(\dd,\rr)$ and $\UU(\dd,\rr)$ be defined as in \cite[Definition 2.2]{BT}. Then, for any fixed $\e$, in the proportional-growth asymptotic,
$$\PP[L_k<\LL(\dd,\rr)+\e]\ra 1\;\;\;\;\mbox{and}\;\;\;\;\PP[U_k<\UU(\dd,\rr)+\e]\ra 1,$$
exponentially in $n$.
\end{lemma}
Using the bounds in Lemma \ref{RIP_bounds}, quantifying the RIP-based  convergence condition (\ref{RIC:IHTstepsize}) for IHT and the corresponding one for N-IHT is straightforward and we are left with quantifying the stable point condition.
In this section, we obtain analogous bounds to those in Lemma \ref{RIP_bounds}, in the proportional-dimensional asymptotic, for the stable point condition (\ref{stablecond_noise}) with Gaussian matrices. These asymptotic bounds are then combined with the RIP bounds of Lemma~\ref{RIP_bounds} in Section~\ref{phase_proofs} to determine a region of phase-space in which recovery is asymptotically guaranteed. The boundary of this region can be viewed as a lower bound on the average-case recovery phase transition for IHT algorithms (see Section~\ref{discussion} for further discussion). 

{\bf Roadmap for the results in Section \ref{phasetrans}.}
The stable point condition was itself analysed using the RIP and the asymptotic bounds of Lemma~\ref{RIP_bounds} in~\cite{thesis}. However, we take a different approach in this paper, motivated by the observation that the stable point condition has no dependence upon the iterates of the algorithm, but depends only upon the original signal, the measurement matrix and the measurement noise. This opens up a possibility that rarely presents itself in the recovery analysis of CS algorithms: we can exploit the reasonable assumption that these three quantities are independent\footnote{A central tenet of CS is the design of nonadaptive measurement schemes, i.e. measurement matrices which are independent of the signal.}.

Our asymptotic bounds for the stable point condition are obtained by first deducing the precise distribution of each term in the stable point condition (or bounds thereon) in terms of the $\chi^2$ and $\FF$ distributions. These distributional results are derived in Section~\ref{stableanalysis}, culminating in Lemma~\ref{dist}. 

Recall that the stable point condition (\ref{stablecond_noise}) takes the form of an inequality that is required to hold over all supports $\Gm\neq\Lm$, where $\Lm=\mbox{supp}(\xs)$. Our asymptotic bounds therefore take the form of tail bounds for a combinatorial number of certain $\chi^2$ and $\FF$ distributions. In Section~\ref{largedev}, we use union bounds to derive three tail bound functions (upper and lower bounds for $\chi^2$ and an upper bound for $\FF$) which are defined implicitly as the solution to equations involving $\dd$ and $\rr$. These results rely on the asymptotic behaviour of the distribution functions for the $\chi^2$ and $\FF$ distributions, which is analysed in Appendix A.

The $\chi^2$ tail bounds have a nice interpretation as bounds on \textit{independent RIP constants}. We close Section~\ref{largedev} by explaining this connection and numerically illustrating that the independence assumption leads to a tightening of RIP bounds.

\subsection{Distribution results for the stable point condition}\label{stableanalysis}

The aim of this section is to derive distribution results in the context of Gaussian measurement matrices for each of the terms in the stable point condition (\ref{stablecond_noise}) of Theorem~\ref{singlefp}. We first give some definitions of Gaussian and Gaussian-related matrix variate distributions, along with some fundamental results concerning their Rayleigh quotients when applied to independent vectors.

We consider a particular kind of matrix variate Gaussian distribution in which all entries are i.i.d. Gaussian random variables, and a few other related distributions.

\begin{definition}[\textbf{Matrix variate Gaussian distribution\cite{handbook}}]\label{gaussian_def}
We say that an $s\times t$ matrix $B$ follows the \textit{matrix variate Gaussian distribution} $B\sim\mN_{s,t}(\mu,\sigma^2)$, if each entry of $B$ independently follows the (univariate) Gaussian distribution $B_{ij}\sim\mN(\mu,\sg^2)$.
\end{definition}


\begin{definition}[\textbf{Matrix variate Wishart distribution~\cite{handbook}}]\label{wishart_def}
Let $B\sim\mathcal{N}_{s,t}(\mu,\sg^2)$ such that $s\geq t$. Then we say that $B^T B$ follows a \textit{matrix variate Wishart distribution} $\mW_t(s;\mu,\sg^2)$ with $s$ degrees of freedom, mean $\mu$ and variance $\sg^2$.
\end{definition}
\begin{definition}[\textbf{$\chi^2$ and $\FF$ distributions~\cite[pp.940,946]{handbook}}]\label{chisq_F_def}
Given a positive integer $s$, let $Z_i\sim\mN(0,1)$ be independent random variables for $1\le i\le s$. Then we say $P=Z_1^2+Z_2^2+\ldots+Z_s^2$ follows a \textit{chi-squared} distribution with $s$ degrees of freedom, and we write $P\sim\chi^2_s$. Furthermore, given positive integers $s$ and $t$, if $P\sim\frac{1}{s}\chi^2_s$ and $Q\sim\frac{1}{t}\chi^2_t$ are independent random variables, we say that $P/Q$ follows the \textit{$\FF$-distribution}, and we write $P/Q\sim\mathcal{F}(s,t)$.
\end{definition}

Crucial to our argument will be the well-known result that the central matrix variate Gaussian distribution defined in Definition~\ref{gaussian_def} is invariant under transformation by an independent orthonormal matrix.

\begin{lemma}[\textbf{Orthogonal invariance~\cite{edelman}}]\label{orth_inv}
Let $B\sim\mathcal{N}_{s,t}(0,\sg^2)$ and let $Z_1\in\RR^{s\times s}$ and $Z_2\in\RR^{t\times t}$ be orthonormal and independent of $B$. Then 
\begin{equation}\label{orthinv1}
Z_1B\sim\mathcal{N}_{s,t}(0,\sg^2),\;\;\;\;\mbox{independently of $Z_1$,}
\end{equation}
and
\begin{equation}\label{orthinv2}
BZ_2\sim\mathcal{N}_{s,t}(0,\sg^2),\;\;\;\;\mbox{independently of $Z_2$.}
\end{equation}
\end{lemma}

Useful results concerning the distributions of Rayleigh quotients related to Gaussian and Wishart matrices are given in the next lemma.

\begin{lemma}[\textbf{Distributions of Rayleigh quotients}]\label{mvdlem}
Let $B\sim\mathcal{N}_{s,t}(0,\sg^2)$ with $s\geq t$. Let $z\in\RR^t$ be
independent of $B$, and such that $\PP(z\neq 0)=1$. Then
\begin{equation}\label{rayleigh1}
\frac{z^T B^T Bz}{z^T z}\sim\sg^2\chi^2_s\;\;\mbox{and is
  independent of $z$;}
\end{equation}
\begin{equation}\label{rayleigh2}
\frac{z^T z}{z^T(B^T
  B)^{-1}z}\sim\sg^2\chi^2_{s-t+1}\;\;\mbox{and is independent of
  $z$;}
\end{equation}
\begin{equation}\label{rayleigh3}
\frac{z^T(B^T B)^2z}{z^T z}\;\;\mbox{has the same
  distribution as $\left\{(B^T B)^2\right\}_{11}$.}
\end{equation}
\end{lemma}

\proof{Let $B\sim\mathcal{N}_{s,t}(0,\sg^2)$ with $s\geq t$. \cite[Theorem 3.3.12]{mvd} gives a more general result than (\ref{rayleigh1}) for when the entries of $B$ are not necessarily independent. The present result follows by setting $\Sigma=\sg^2 I$ for the covariance matrix. Similarly, (\ref{rayleigh2}) follows by setting $\Sigma=\sg^2 I$ in \cite[Corollary 3.3.14.1]{mvd}. To prove (\ref{rayleigh3}), let $S=B^T B$ so that $S\sim\mathcal{W}_t(s;0,\sg^2)$ and let $Z\in\RR^{t\times t}$ be any orthonormal matrix which is independent of $B$. Lemma~\ref{orth_inv} yields $BZ\sim\mathcal{N}_{s,t}(0,\sg^2)$ independently of $Z$, and, writing $T:=Z^T SZ$, we therefore have
\begin{equation}\label{Tdist}
T=Z^T SZ=Z^T B^T BZ=(BZ)^T BZ\sim\mathcal{W}_t(s;0,\sg^2),
\end{equation}
independently of $Z$. In particular, let us fix the first column of $Z$ as $z$ normalized so that
$$Z=\left.\left[\displaystyle\frac{z}{\|z\|}\right|Z_2\right],$$
which leads to
$$\begin{array}{l}
z^TS^2z=z^T(ZTZ^T)^2 z=z^TZTZ^TZTZ^T z=z^TZT^2 Z^T z\\
=z^T
\left.\left[\displaystyle\frac{z}{\|z\|}\right|Z_2\right]T^2
\left[\begin{array}{c}
\frac{z^T}{\|z\|}\\
\hline
Z_2^T
\end{array}\right]
z
=\left[\;\|z\|\;|\;0\;\right]T^2\left[\begin{array}{c}
\|z\|\\
\hline
0
\end{array}\right]
=(T^2)_{11}\|z\|^2.
\end{array}$$
Dividing by $\|z\|^2$ and using (\ref{Tdist}) then gives the desired result.}

We now make the assumption that the measurement matrix in \req{sparserecovery}
is drawn from the (central) matrix variate Gaussian distribution with appropriate normalization.

\ass{gaussianass}{\quad The measurement matrix $A$ has i.i.d. $\mathcal{N}(0,1/n)$ entries, so that $A\sim\mathcal{N}_{n,N}(0,1/n)$. Furthermore, $A$ is independent of $\xs$.}

Given Assumption {\bf A.2} and the standard compressed sensing regime with $2k\leq n$,
we can dispense with Assumption {\bf A.1} \cite[Theorem 3.2.1]{mvd}\footnote{\cite[Theorem 3.2.1]{mvd} states that
$B^T B$ is positive definite with probability $1$ when  $B\sim\mathcal{N}_{s,t}(0,\sg^2)$ with $s\geq t$.}.



We also impose the additional assumption that measurement noise is itself Gaussian and independent of both the original signal and the measurement matrix.

\ass{noiseass}{\quad The noise vector $e$ has i.i.d. Gaussian entries $e_i\sim N(0,\sigma^2/n)$, independently of $A$ and $\xs$.}

Note that, under Assumption {\bf A.3}, $\EE\|e\|^2=\sg^2$, so that $\|e\|\approx\sg$.

We now give the main result of this section, in which we derive precise distributions for various expressions which make up the stable point condition (\ref{stablecond_noise}) of Theorem~\ref{singlefp}, in terms of the $\chi^2$ and $\FF$ distributions.

\begin{lemma}[\textbf{Distribution results for the stable point condition}]\label{dist}
Suppose Assumptions {\bf A.2} and {\bf A.3} hold, and
let $\Gm$ and $\Lm$ be index sets of cardinality $k$, where $k<n$, such that $\Gm\neq\Lm$. Then
\begin{equation}\label{lhs}
\frac{\|A_{\Gm}^{\dag}A_{\Lm\sm\Gm}\xs_{\Lm\sm\Gm}\|}{\|\xs_{\Lm\sm\Gm}\|}=\sqrt{F_{\Gm}},\;\;\;\mbox{where}\;\;\;F_{\Gm}\sim\frac{k}{n-k+1}\FF(k,n-k+1);
\end{equation}
\begin{equation}\label{rhs}
\frac{\|A_{\Lm\sm\Gm}^T(I-A_{\Gm}A^{\dag}_{\Gm})A_{\Lm\sm\Gm}\xs_{\Lm\sm\Gm}\|}{\|\xs_{\Lm\sm\Gm}\|}\geq{\left(\frac{n-k}{n}\right)}\cdot R_{\Gm},\;\;\;\mbox{where}\;\;\;R_{\Gm}\sim\frac{1}{n-k}\chi^2_{n-k};
\end{equation}
\begin{equation}\label{lhsnoise}
\|A_{\Gm}^{\dag}e\|\le\sg\cdot\sqrt{G_{\Gm}},\;\;\;\mbox{where}\;\;\;G_{\Gm}\sim\frac{k}{n-k+1}\mathcal{F}(k,n-k+1);
\end{equation}
\begin{equation}\label{rhsnoise}
\|A_{\Lm\sm\Gm}^T(I-A_{\Gm}A_{\Gm}^{\dag})e\|\le\sg\sqrt{\frac{k(n-k)}{n^2}\cdot(S_{\Gm})(T_{\Gm})},\;\;\;\mbox{where}\;\;\;S_{\Gm}\sim\frac{1}{n-k}\chi^2_{n-k},\;\;\;T_{\Gm}\sim\frac{1}{k}\chi^2_k.
\end{equation}
\end{lemma}

\textbf{Proof of (\ref{lhs}):} Let $A_{\Gm}$ have the singular value
decomposition
\begin{equation}\label{AgSVD}
A_{\Gm}:=U[D\;|\;0]V^T=U_1 DV^T,
\end{equation}
where $D\in\RR^{k\times k}$ is diagonal, and where $V\in\RR^{k\times k}$ and
$U=[U_1\;|\;U_2]\in\RR^{n\times n}$ are orthonormal, with $U_1\in\RR^{n\times k}$. By Assumption {\bf A.2}, $A_{\Gm}^{\dag}$ is well-defined and we have the standard result
\begin{equation}\label{dagger}
A_{\Gm}^{\dag}=VD^{-1}U_1^T,
\end{equation}
and since $(A_{\Gm}^TA_{\Gm})^{-1}=VD^{-2}V^T,$ it follows by
rearrangement that
\begin{equation}\label{invwish}
D^{-2}=V^T(A_{\Gm}^T A_{\Gm})^{-1}V.
\end{equation}
Using (\ref{dagger}) and (\ref{invwish}), we have
\begin{eqnarray}\label{rewrite}
\displaystyle\|A_{\Gm}^{\dag}A_{\Lm\sm\Gm}\xs_{\Lm\sm\Gm}\|^2&=&(\xs_{\Lm\sm\Gm})^T
A_{\Lm\sm\Gm}^T(A_{\Gm}^{\dag})^T(A_{\Gm}^{\dag})A_{\Lm\sm\Gm}\xs_{\Lm\sm\Gm}\nonumber\\
&=&(\xs_{\Lm\sm\Gm})^T
A_{\Lm\sm\Gm}^TU_1 D^{-1}V^T VD^{-1}U_1^T A_{\Lm\sm\Gm}\xs_{\Lm\sm\Gm}\nonumber\\
&=&(\xs_{\Lm\sm\Gm})^T
A_{\Lm\sm\Gm}^TU_1 D^{-2}U_1^T A_{\Lm\sm\Gm}\xs_{\Lm\sm\Gm}\nonumber\\
&=&(\xs_{\Lm\sm\Gm})^T
A_{\Lm\sm\Gm}^TU_1 V^T(A_{\Gm}^T A_{\Gm})^{-1}VU_1^T
A_{\Lm\sm\Gm}\xs_{\Lm\sm\Gm}.
\end{eqnarray}
By Lemma~\ref{orth_inv}, we have $U^T A_{\Lm\sm\Gm}\sim\mathcal{N}_{n,r}(0,1/n)$, independently of $U$, where $r:=|\Lm\sm\Gm|$. Since $U_1^T A_{\Lm\sm\Gm}$ is a submatrix of $U^T A_{\Lm\sm\Gm}$, it follows that $U_1^T A_{\Lm\sm\Gm}\sim\mathcal{N}_{k,r}(0,1/n)$, independently of $U$. Writing $C:=VU_1^T A_{\Lm\sm\Gm}\in\RR^{k\times r}$, we also have by Lemma~\ref{orth_inv} that $C\sim\mathcal{N}_{k,r}(0,1/n)$, independently of both $U$ and $V$, and
therefore independently of $A_{\Gm}$. Substituting for $C$ in (\ref{rewrite}), we have
\begin{eqnarray}\displaystyle
\frac{\|A_{\Gm}^{\dag}A_{\Lm\sm\Gm}\xs_{\Lm\sm\Gm}\|^2}{\|\xs_{\Lm\sm\Gm}\|^2}&=&\displaystyle\frac{(\xs_{\Lm\sm\Gm})^T
  C^T (A_{\Gm}^T A_{\Gm})^{-1} C\xs_{\Lm\sm\Gm}}{(\xs_{\Lm\sm\Gm})^T
  \xs_{\Lm\sm\Gm}}\nonumber\\
&=&\displaystyle\frac{(\xs_{\Lm\sm\Gm})^T
  C^T (A_{\Gm}^T A_{\Gm})^{-1} C\xs_{\Lm\sm\Gm}}{(\xs_{\Lm\sm\Gm})^T C^T
  C\xs_{\Lm\sm\Gm}}\;\cdot\;\displaystyle\frac{(\xs_{\Lm\sm\Gm})^T C^T C\xs_{\Lm\sm\Gm}}{(\xs_{\Lm\sm\Gm})^T
  \xs_{\Lm\sm\Gm}},\label{split}
\end{eqnarray}
where $\xs$, $C$ and $A_{\Gm}$ are all independent. Now it follows from Lemma~\ref{mvdlem} that
\begin{equation}\label{chisqsplit}
\frac{(\xs_{\Lm\sm\Gm})^T C^T C\xs_{\Lm\sm\Gm}}{(\xs_{\Lm\sm\Gm})^T
  \xs_{\Lm\sm\Gm}}\sim\frac{1}{n}\;\chi^2_k\;\;\;\;\mbox{and}\;\;\;\;\frac{(\xs_{\Lm\sm\Gm})^T C^T
  C\xs_{\Lm\sm\Gm}}{(\xs_{\Lm\sm\Gm})^T
  C^T (A_{\Gm}^T A_{\Gm})^{-1}
  C\xs_{\Lm\sm\Gm}}\sim\frac{1}{n}\;\chi^2_{n-k+1},
\end{equation}
where both distributions are independent of each other. Combining (\ref{split}) and (\ref{chisqsplit}) leads us to conclude
$$\frac{\|A_{\Gm}^{\dag}A_{\Lm\sm\Gm}\xs_{\Lm\sm\Gm}\|^2}{\|\xs_{\Lm\sm\Gm}\|^2}\sim\frac{\chi^2_k}{\chi^2_{n-k+1}}=\frac{k}{n-k+1}\mathcal{F}(k,n-k+1),$$
where in the last step we use the fact that the two distributions are independent, which proves (\ref{lhs}).

\textbf{Proof of (\ref{rhs}):} Using (\ref{AgSVD}) and (\ref{dagger}), we have
$$A_{\Gm}A_{\Gm}^{\dag}=U_1 DV^T VD^{-1}U_1^T=U_1 U_1^T=U\left[\begin{array}{cc}
I&0\\
0&0\end{array}\right]U^T,$$
and writing $I=UU^T$,
\begin{equation}\label{U2result}
I-A_{\Gm}A_{\Gm}^{\dag}=U\left\{\left[\begin{array}{cc}I&0\\0&I\end{array}\right]-\left[\begin{array}{cc}I&0\\0&0\end{array}\right]\right\}U^T=U\left[\begin{array}{cc}0&0\\0&I\end{array}\right]U^T=U_2 U_2^T,
\end{equation}
which in turn gives
\begin{equation}\label{beforeD}A_{\Lm\sm\Gm}^T(I-A_{\Gm}A^{\dag}_{\Gm})A_{\Lm\sm\Gm}=A_{\Lm\sm\Gm}^TU_2
U_2^T A_{\Lm\sm\Gm}.\end{equation}
Writing $F:=U_2^T A_{\Lm\sm\Gm}$, we have $U^T A_{\Lm\sm\Gm}\sim\mathcal{N}_{n,r}(0,1/n)$ by Lemma~\ref{orth_inv}, and since $U_2^T A_{\Lm\sm\Gm}\in\RR^{(n-k)\times r}$ is a submatrix of $U^T A_{\Lm\sm\Gm}$, it follows that
\begin{equation}\label{Fdistn}
F\sim\mathcal{N}_{(n-k),r}(0,1/n).
\end{equation}
Substituting for $F$ in (\ref{beforeD}) gives
\begin{equation}\label{withD}
A_{\Lm\sm\Gm}^T(I-A_{\Gm}A^{\dag}_{\Gm})A_{\Lm\sm\Gm}=F^T F.
\end{equation}
Now, writing $M:=F^T F$, and using (\ref{withD}) and (\ref{rayleigh3}) of Lemma~\ref{mvdlem}, we deduce
\begin{equation}\label{Msq}
\displaystyle\frac{\|A_{\Lm\sm\Gm}^T(I-A_{\Gm}A^{\dag}_{\Gm})A_{\Lm\sm\Gm}\xs_{\Lm\sm\Gm}\|^2}{\|\xs_{\Lm\sm\Gm}\|^2}=\displaystyle\frac{\|F^T
  Fx_{\Lm\sm\Gm}\|^2}{\|\xs_{\Lm\sm\Gm}\|^2}=\displaystyle\frac{(\xs_{\Lm\sm\Gm})^T(F^T F)^2 \xs_{\Lm\sm\Gm}}{(\xs_{\Lm\sm\Gm})^T \xs_{\Lm\sm\Gm}}\sim(M^2)_{11}.
\end{equation}
To obtain a lower bound in terms of the chi-squared distribution, note that
\begin{equation}\label{M}
(M^2)_{11}=\sum_{i=1}^r M_{i1}^2=M_{11}^2 + \sum_{i=2}^r
M_{i1}^2\geq M_{11}^2.
\end{equation}
Meanwhile it follows from (\ref{Fdistn}) and (\ref{chisq_F_def}) that
$$M_{11}=\sum_{i=1}^{n-k} F_{i1}^2\sim\frac{1}{n}\chi^2_{n-k},$$
which combines with (\ref{Msq}) and (\ref{M}) to give (\ref{rhs}).

\textbf{Proof of (\ref{lhsnoise}):} By (\ref{dagger}), we have
\begin{equation}\label{error1}
A^{\dag}_{\Gm}e=VD^{-1}U_1^T e=VD^{-1}p,
\end{equation}
where $p:=U_1^T e\in\RR^{n-k}$. Using Assumption {\bf A.3}, we may view $e$ as a one-column Gaussian matrix, such that $e\sim\mathcal{N}_{n,1}(0,\sg^2/n)$, it follows from Lemma~\ref{orth_inv} that
\begin{equation}\label{pdist}
p\sim\mathcal{N}_{k,1}(0,\sg^2/n),
\end{equation}
independently of $U$ and therefore independently of $A_{\Gm}$. Substituting (\ref{invwish}) into (\ref{error1}) then gives
\begin{equation}\label{quadform}
\|A_{\Gm}^{\dag}e\|^2=\|VD^{-1}p\|^2=\|D^{-1}p\|^2=p^T D^{-2}p=p^T V^T(A_{\Gm}^T A_{\Gm})^{-1}Vp=q^T(A_{\Gm}^T A_{\Gm})^{-1}q,
\end{equation}
where $q:=Vp\in\RR^k$. It now follows from (\ref{pdist}) and Lemma~\ref{orth_inv} that $q\sim\mathcal{N}_{k,1}(0,\sg^2/n)$, independently of $V$ and therefore independently of $A_{\Gm}$, and consequently that
\begin{equation}\label{normal2norm}
q^T q\sim\sg^2\chi^2_k.
\end{equation}
By (\ref{rayleigh2}) of Lemma~\ref{mvdlem},
\begin{equation}\label{restrdist}
\frac{q^T q}{q^T(A_{\Gm}^T
  A_{\Gm})^{-1}q}\sim\frac{1}{n}\chi^2_{n-k+1}.
\end{equation}
Since $q$ and $A_{\Gm}$ are independent, we may combine (\ref{quadform}), (\ref{normal2norm}) and (\ref{restrdist}) to give
\begin{equation}\label{edashdist}
\|A^{\dag}_{\Gm}e\|\sim\sg\sqrt{G_{\Gm}},\;\;\;\mbox{where}\;\;\;G_{\Gm}\sim\frac{k}{n-k+1}\mathcal{F}(k,n-k+1),
\end{equation}
and (\ref{lhsbound}) now follows.

\textbf{Proof of (\ref{rhsnoise}):} Using (\ref{U2result}), we have
\begin{equation}\label{orthprojnoise}
A_{\Lm\sm\Gm}^T(I-A_{\Gm}A_{\Gm}^{\dag})e=A_{\Lm\sm\Gm}^T U_2 U_2^Te=A_{\Lm\sm\Gm}^T U_2 f=B^T f,
\end{equation}
where $B:=U_2^T A_{\Lm\sm\Gm}\sim\mathcal{N}_{n-k,r}(0,1/n)$ by Lemma~\ref{orth_inv}, and where
\begin{equation}\label{fdist}
f:=U_2^Te\sim\mathcal{N}_{n-k,1}(0,\sg^2/n)
\end{equation}
by Lemma~\ref{orth_inv}. Now let $B$ have singular value
decomposition
\begin{equation}\label{BSVD}
W[F\;|\;0]Y^T=W_1FY^T,
\end{equation}
where $F\in\RR^{r\times r}$ is diagonal, and where $Y\in\RR^{r\times r}$ and
$W=[W_1\;|\;W_2]\in\RR^{(n-k)\times(n-k)}$ are orthonormal, noting
that $W_1\in\RR^{(n-k)\times r}$. We have
\begin{equation}\label{gdist}
g:=W_1^T f\sim\mathcal{N}_{r,1}(0,\sg^2/n)
\end{equation}
by (\ref{fdist}) and Lemma~\ref{orth_inv}, and we may apply (\ref{orthprojnoise}) to give
\begin{equation}\label{quadformorth}
\begin{array}{l}
\|A_{\Lm\sm\Gm}^T(I-A_{\Gm}A_{\Gm}^{\dag})e\|^2\le\|B^T f\|^2=\|YFW_1^T f\|^2\\
=\|Fg\|^2=g^T F^2 g=g^T Y^T(B^T B)Yg=h^T(B^T B)h,
\end{array}
\end{equation}
where $h:=Yg\in\RR^k$. Since $h\sim\mathcal{N}_{r,1}(0,\sg^2/n)$ by (\ref{gdist}) and Lemma~\ref{orth_inv}, it follows that
\begin{equation}\label{normal2normorth}
h^T h\sim\sg^2\chi^2_r\le\sg^2\chi^2_k,
\end{equation}
since a $\chi^2_r$ random variate may be viewed as a truncation of its extension to a $\chi^2_k$ random variate. By (\ref{rayleigh1}) of Lemma~\ref{mvdlem},
\begin{equation}\label{restrdistorth}
\frac{h^T B^T
  Bh}{h^T h}\sim\frac{1}{n}\chi^2_{n-k}.
\end{equation}
Combining (\ref{quadformorth}), (\ref{normal2normorth}) and (\ref{restrdistorth}) then proves (\ref{rhsnoise}).\hfill$\Box$\\

In order for the (converse of the) stable point condition \req{stablecond_noise} to provide a recovery result regarding the proximity of all stable points to the
underlying signal, we need to quantify the quantities in Lemma \ref{dist} on all possible fixed points on support sets $\Gamma$ of cardinality $k$. Similarly,
the convergence conditions in Section \ref{convergencesection} involve  RIP constants which again involve looking over combinatorially many supports.
Thus we need to derive union bounds for the relevant distributions involved in the stable point and convergence conditions.

\subsection{Large deviation results involving Gaussian matrices}\label{largedev}

In this section, we derive large deviations results for quantities relating to Gaussian matrices within the proportional-growth asymptotic that is defined on page 2. We define three tail bound functions.

\begin{definition}[\textbf{$\chi^2$ tail bounds}]\label{Idef}
Let $\dd\in(0,1]$, $\rr\in(0,1)$ and $\lm\in(0,1]$. Let the tail bound function $\IU(\dd,\rr,\lm)$  be the unique solution to
\begin{equation}\label{udef}
\nu-\ln(1+\nu)=\frac{2H(\dd\rr)}{\lm}\;\;\;\;\mbox{for}\;\;\;\;\nu>0,
\end{equation}
and let the tail bound function $\IL(\dd,\rr,\lm)$ be the unique solution to
\begin{equation}\label{ldef}
-\nu-\ln(1-\nu)=\frac{2H(\dd\rr)}{\lm}\;\;\;\;\mbox{for}\;\;\;\;\nu\in(0,1),
\end{equation}
where $H(\cdot)$ is the  Shannon entropy with base $e$ logarithms \cite{lqphase}, namely,
\begin{equation}\label{shannon_def}
H(p):=-p\ln(p)-(1-p)\ln(1-p).
\end{equation}
\end{definition}

That $\IU$ is well-defined follows since the left-hand side of (\ref{udef}) is zero at $\nu=0$, tends to infinity as $\nu\ra\infty$, and is strictly increasing on $\nu>0$. Similarly, $\IL$ is well-defined since the left-hand side of (\ref{ldef}) is zero at $\nu=0$, tends to infinity as $\nu\ra 1$, and is strictly increasing on $\nu\in(0,1)$.

\begin{definition}[\textbf{$\FF$ tail bound}]\label{IFdef}
Let $\dd\in(0,1]$ and $\rr\in(0,1/2]$. Let the tail bound function $\IFF(\dd,\rr)$ be the unique solution in $f$ to
\begin{equation}\label{Fdef}
\ln(1+f)-\rr\ln f=2H(\dd\rr)+H(\rr)\;\;\;\;\mbox{for}\;\;\;\;f>\frac{\rr}{1-\rr},
\end{equation}
where $H(\cdot)$ is defined in (\ref{shannon_def}).
\end{definition}

That $\IFF$ is well-defined follows since the left-hand side of (\ref{Fdef}) is equal to $H(\rr)$ at $f=\rr/(1-\rr)$, tends to infinity as $f\ra\infty$, and is strictly increasing on $f>\rr/(1-\rr)$.

Defining $S_n$ as
\begin{equation}\label{Sn_def}
S_n\eqdef\left\{1,\ldots,{N\choose k}\right\},
\end{equation}
we have the following large deviation bound for a combinatorial number of $\chi^2$ distributions.

\begin{lemma}[\textbf{Large deviations result for $\chi^2$}]\label{chisq}
Let $l\in\{1,\ldots,n\}$ and let the random variables $X_l^i\sim\displaystyle\frac{1}{l}\chi^2_l$ for all $i\in S_n$, and let $\e>0$. In the proportional-growth asymptotic, let $l/n\ra\lm\in(0,1]$. Then
\begin{equation}\label{chisqresult1}
\PP\left\{\cup_{i\in S_n}[X_l^i\geq 1+\IU(\dd,\rr,\lm)+\e]\right\}\rightarrow
0
\end{equation}
and
\begin{equation}\label{chisqresult2}
\PP\left\{\cup_{i\in S_n}[X_l^i\le 1-\IL(\dd,\rr,\lm)-\e]\right\}\rightarrow
0,
\end{equation}
exponentially in $n$, where $\IU(\dd,\rr,\lm)$ and $\IL(\dd,\rr,\lm)$ are defined in (\ref{udef}) and (\ref{ldef}) respectively.
\end{lemma}

The proof of Lemma \ref{chisq} is delegated to Appendix A. It employs asymptotic results derived by Temme \cite{temme82} for
the incomplete gamma function which is related to the $\chi^2$ distribution.

\begin{lemma}[\textbf{Large deviations result for $\FF$}]\label{Fdist}
Let the random variables $X_n^i\sim\frac{k}{n-k+1}\;\mathcal{F}(k,n-k+1)$ for all $i\in S_n$, and let $\e>0$. In the proportional-growth asymptotic,
\begin{equation}\label{Fresult}
\PP\left\{\cup_{i\in S_n}[X_n^i\geq \IFF(\dd,\rr)+\e]\right\}\rightarrow 0,
\end{equation}
exponentially in $n$, where $\IFF(\dd,\rr)$ is defined in (\ref{Fdef}).
\end{lemma}

The proof of Lemma \ref{Fdist} is delegated to Appendix A. It employs asymptotic results derived by Temme \cite{temme82} for
the incomplete beta function which is related to the $\FF$ distribution.\\
\\
{\bf Comparison of Lemmas~\ref{RIP_bounds} and~\ref{chisq}; RIP versus Independent RIP constants.}\quad
Suppose $A\sim\mathcal{N}_{n,N}(0,1/n)$, let $\Gm$ be an index set of cardinality $k$ and fix $\e>0$.
Then in the conditions of Lemma~\ref{RIP_bounds}, in the proportional-growth asymptotic, for any $y\in\RR^k$,
\begin{equation}\label{rayleigh_RIP}
1-\LL(\dd,\rr)-\e<\frac{\|A_{\Gm}y\|^2}{\|y\|^2}<1+\UU(\dd,\rr)+\e.
\end{equation}
However, if $y$ is independent of $A$, we may set $\lm=1$ in Lemma~\ref{chisq}, giving in the proportional-growth asymptotic,
\begin{equation}\label{rayleigh_ind}
1-\IL(\dd,\rr,1)-\e<\frac{\|A_{\Gm}y\|^2}{\|y\|^2}<1+\IU(\dd,\rr,1)+\e.
\end{equation}
Comparing (\ref{rayleigh_RIP}) and (\ref{rayleigh_ind}), we see that $\IU(\dd,\rr,1)$ and $\IL(\dd,\rr,1)$ may be viewed as upper bounds on `{\it independent RIP}' constants for Gaussian matrices.

\begin{figure}[ht]
\centering
\subfigure[]{\includegraphics[width=2.8in,height=2in]{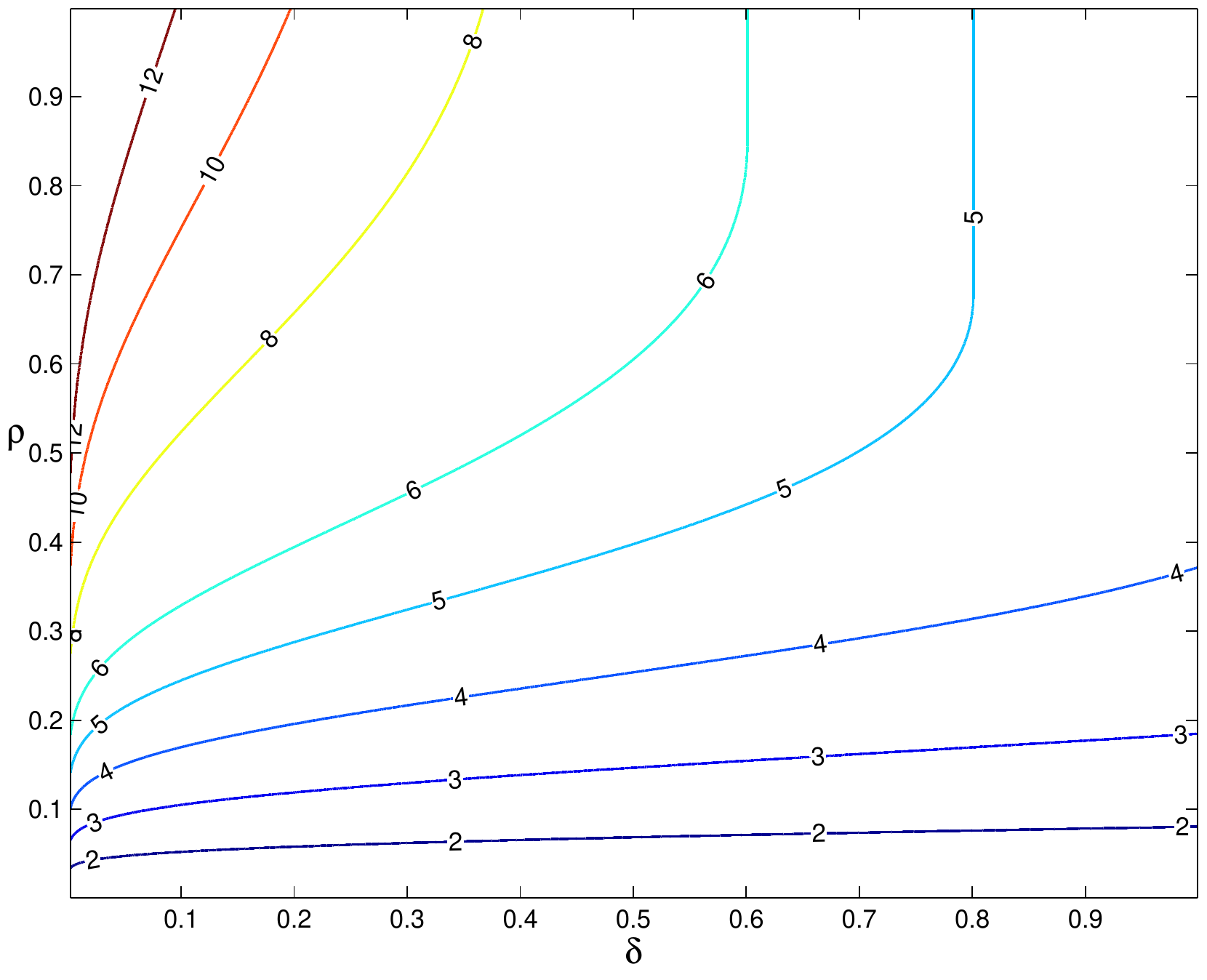}}
\subfigure[]{\includegraphics[width=2.8in,height=2in]{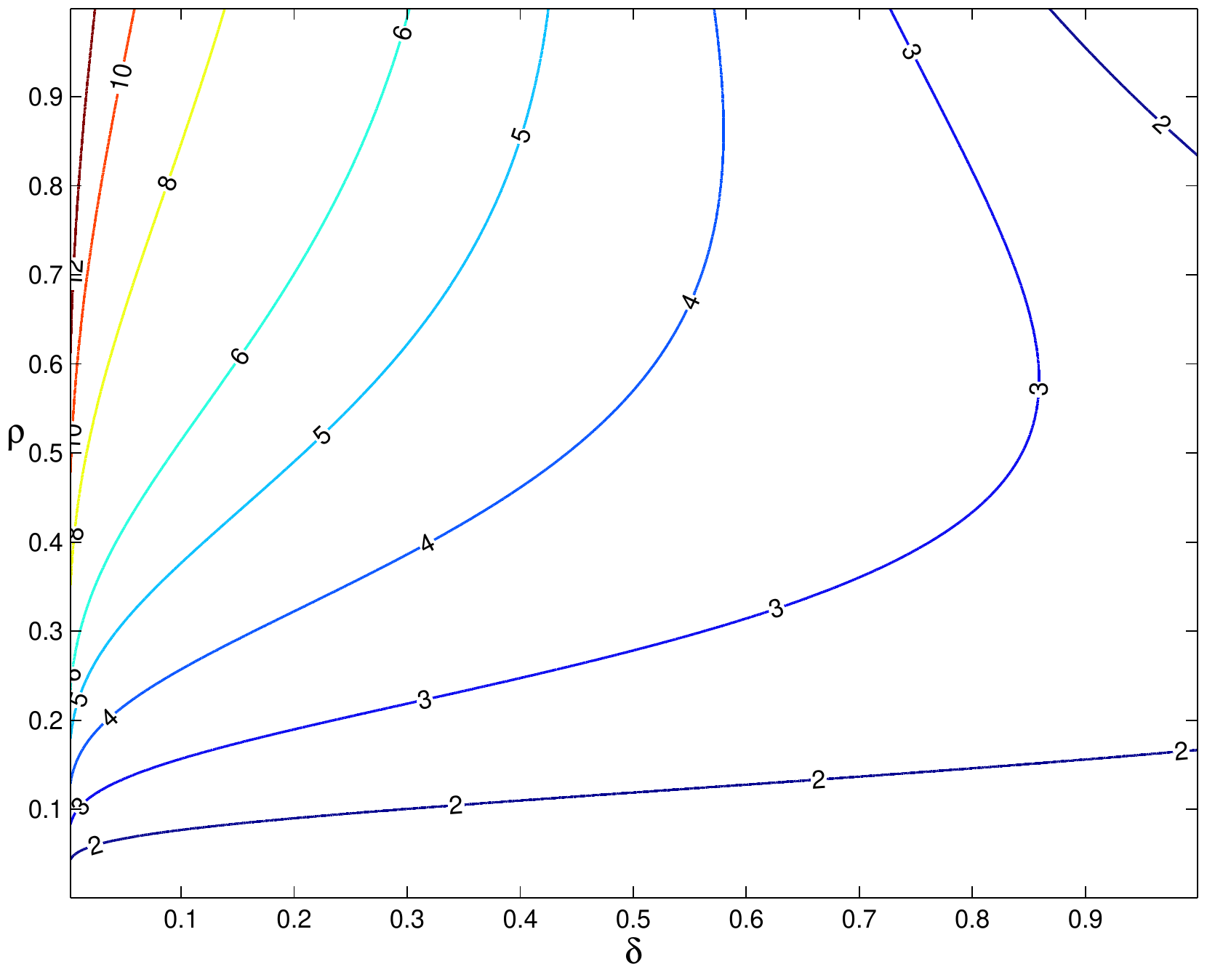}}
\subfigure[]{\includegraphics[width=2.8in,height=2in]{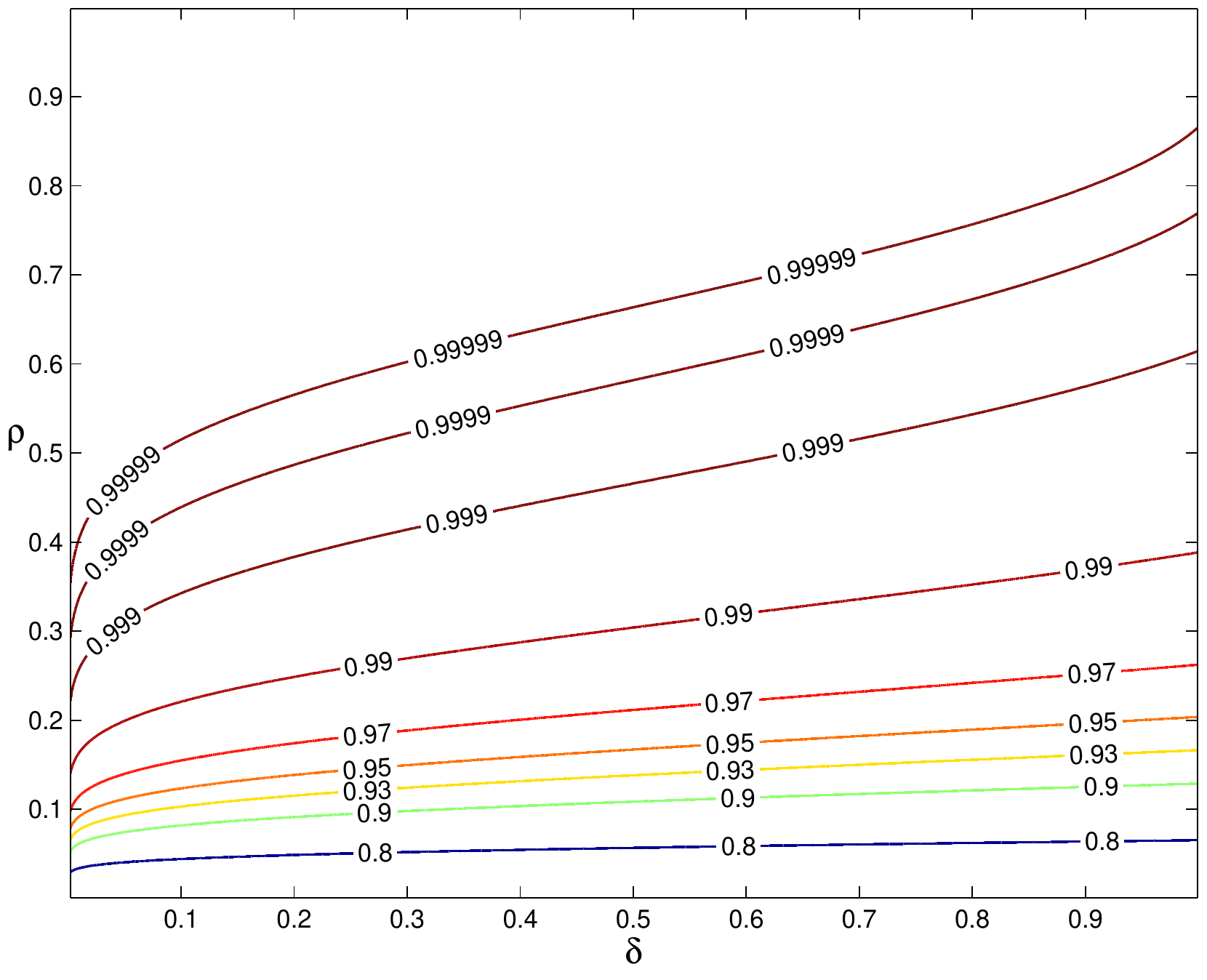}}
\subfigure[]{\includegraphics[width=2.8in,height=2in]{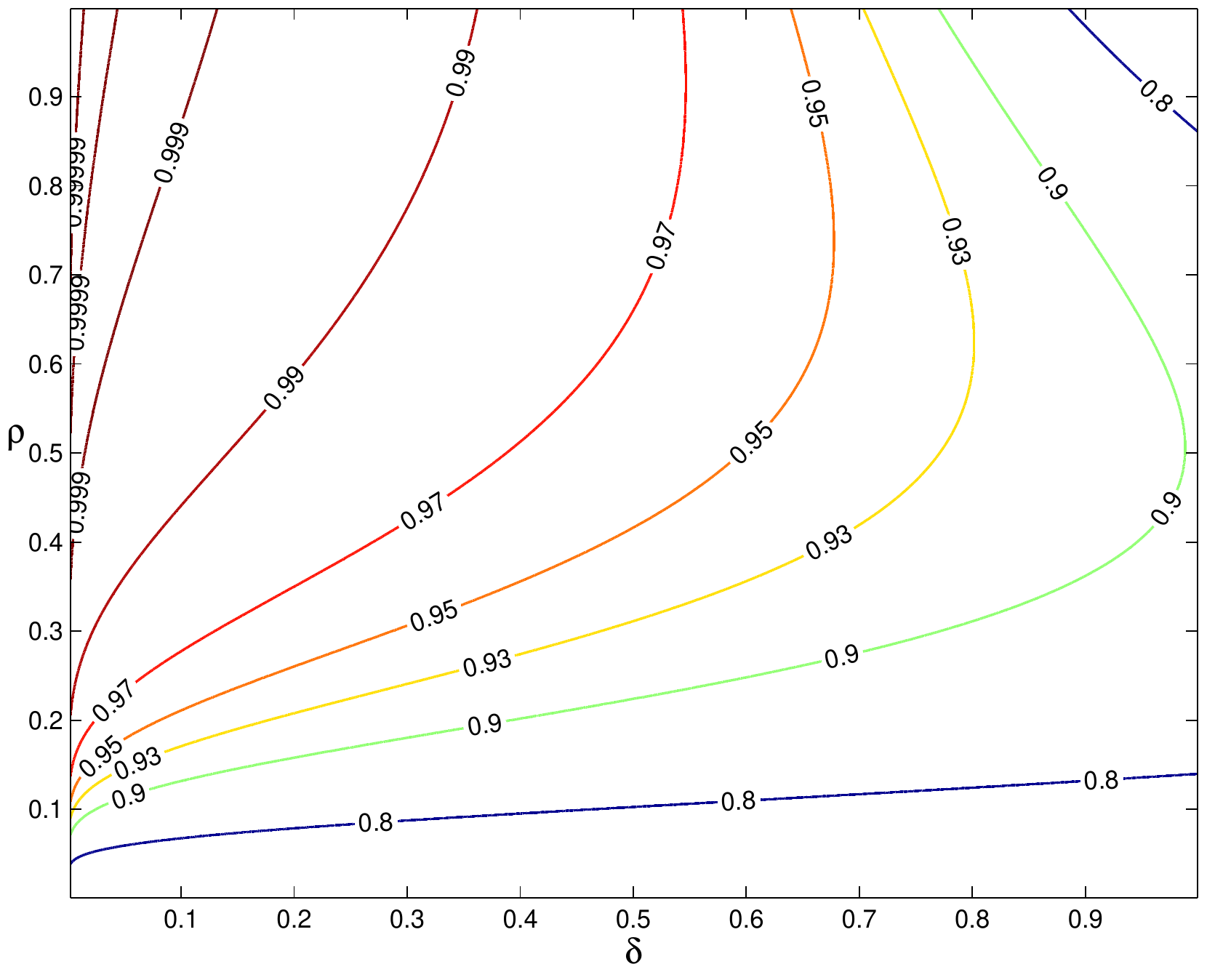}}
\caption{A comparison of standard RIP bounds \cite{BT} and `independent RIP' bounds for Gaussian matrices:
(a) $\UU(\dd,\rr)$ (b) $\IU(\dd,\rr,1)$ (c) $\LL(\dd,\rr)$ (d) $\IL(\dd,\rr,1)$.}
\label{bounds_plots}
\end{figure}

Figure~\ref{bounds_plots} gives plots of the `independent RIP' bounds for Gaussian matrices $\IU(\dd,\rr,1)$ and $\IL(\dd,\rr,1)$ derived in this paper, along with plots of the RIP bounds for Gaussian matrices $\UU(\dd,\rr)$ and $\LL(\dd,\rr)$ in \cite{BT}. One observes empirically the inequalities
$$\IU(\dd,\rr,1)<\UU(\dd,\rr)\;\;\;\;\mbox{and}\;\;\;\;\IL(\dd,\rr,1)<\LL(\dd,\rr).$$
A simple interpretation is that the additional information that the matrix and vector are independent allows us to tighten the bounds in (\ref{rayleigh_RIP}) to obtain (\ref{rayleigh_ind}). This consideration accounts for a large part of the quantitative improvement that is obtained in this paper over existing recovery results for IHT algorithms which rely solely upon the RIP.
Of course, our improved analysis is only possible because our proposed stable point condition (\ref{stablecond_noise}) can exploit the assumption of
matrix-vector independence.

\numsection{Novel recovery analysis for IHT algorithms}\label{phase_proofs}

{\bf Roadmap for the results in Section \ref{phase_proofs}.}
Revisiting the roadmap at the start of Section~\ref{phasetrans}, we find that we now have 
all the necessary asymptotic bounds for quantifying both the RIP-based convergence conditions and the stable point conditions
required for ensuring recovery using IHT and N-IHT. Putting these ingredients together, we are now ready to present our main
quantitative recovery results for IHT and N-IHT when Gaussian measurement matrices are employed. 
We begin the next section by 
defining the phase transition bounds and noise stability factors which feature in the statement of our main recovery results. We then
state our main recovery results for IHT and N-IHT, respectively, before illustrating and discussing their significance in Section~\ref{discussion}. 
Proofs for the IHT results can be found in Section~\ref{IHT_proofs} (with a roadmap for the
line of argument given at the start of the respective section), while the proofs for N-IHT (which follow very similar lines) are delegated to Appendix B.


\subsection{Statement of main recovery results}\label{statement}

\subsubsection{Results for IHT}\label{IHTresults_Section}

We first give definitions of the lower bound on the phase transition and noise stability factor featuring in the main result.
The function $\rh^{IHT}(\dd)$ is a lower bound on the phase transition for recovery using IHT (see Section~\ref{discussion} for further explanation). The function $\xi(\dd,\rr)$ represents a stability factor in our results, bounding the approximation error of the output of IHT as a multiple of the noise level $\Sg$. Both functions are numerically computable.


\begin{definition}[\textbf{Phase transition lower bound for IHT}]\label{rhIHT}
Given $\dd\in(0,1]$, define the phase transition lower bound $\rh^{IHT}(\dd)$ to be the unique solution to
\begin{equation}\label{rhoIHTdef}
\frac{\sqrt{\IFF(\dd,\rr)}}{(1-\rr)\left[1-\IL(\dd,\rr,1-\rr)\right]}=\frac{1}{1+\UU(\dd,2\rr)}\;\;\;\;\mbox{for}\;\;\;\;\rr\in(0,1/2],
\end{equation}
where $\IFF$ is defined in (\ref{Fdef}), $\IL$ is defined in (\ref{ldef}), $\UU$ is defined in~\cite[Definition 2.2]{BT}.\footnote{A proof that $\rh^{IHT}(\dd)$ is well-defined can be found in \cite[Section 5.2]{thesis}.}
\end{definition}


\begin{definition}[\textbf{Stability factor for IHT}]\label{massive_def}
Given $\dd\in(0,1]$, $\rr\in(0,1/2]$ and $\al>0$, provided
\begin{equation}\label{SSPcond}
\al>\frac{\sqrt{\IFF(\dd,\rr)}}{(1-\rr)[1-\IL(\dd,\rr,1-\rr)]},
\end{equation}
define the stability factor $\xi(\dd,\rr)$ to be
\begin{equation}\label{xidef}
\xi(\dd,\rr)\eqdef\sqrt{\IFF(\dd,\rr)\left[1+a(\dd,\rr)\right]^2+\left[a(\dd,\rr)\right]^2},
\end{equation}
where\footnote{Note that (\ref{SSPcond}) ensures that the denominator in (\ref{adefn}) is strictly positive and that $a(\dd,\rr)$ is therefore well-defined.}
\begin{equation}\label{adefn}
a(\dd,\rr)\eqdef\frac{\sqrt{\IFF(\dd,\rr)}+\al\sqrt{\rr(1-\rr)[1+\IU(\dd,\rr,1-\rr)][1+\IU(\dd,\rr,\rr)]}}{\al(1-\rr)[1-\IL(\dd,\rr,1-\rr)]-\sqrt{\IFF(\dd,\rr)}},
\end{equation}
and where $\IFF$ is defined in (\ref{Fdef}), $\IU$ is defined in (\ref{udef}), and where $\IL$ is defined in (\ref{ldef}).
\end{definition}


We have the following recovery result for IHT.

\begin{theorem}[\textbf{Recovery result for IHT; noise case}]\label{recov1noise}
Suppose Assumptions {\bf A.2} and {\bf A.3} hold, suppose that
\begin{equation}\label{IHTstablecond}
\rr<\rh^{IHT}(\dd),
\end{equation}
where $\rh^{IHT}(\dd)$ is defined in (\ref{rhoIHTdef}), and  that the IHT stepsize $\al$ satisfies
\begin{equation}\label{alphabound1}
\frac{\sqrt{\IFF(\dd,\rr)}}{(1-\rr)\left[1-\IL(\dd,\rr,1-\rr)\right]}<\al<\frac{1}{1+\UU(\dd,2\rr)}.
\end{equation}
Then, in the proportional-growth asymptotic\footnote{In other words, we consider instances of the Gaussian random variables $A$ and $e$ for a sequence of triples
$(k,n,N)$ where $n\rightarrow\infty$, where $n$ is the number of measurements, $N$, the signal dimension and $k$, the sparsity of the underlying signal.}, IHT converges to $\xb$ that is close to $\xs$ in the sense that
\begin{equation}\label{error}
\|\xb-\xs\|\le\xi(\dd,\rr)\cdot\Sg
\end{equation}
holds with probability tending to $1$ exponentially in $n$, where $\xi(\dd,\rr)$ is defined in (\ref{xidef}).
\end{theorem}
Comparing our result with previous RIP-based recovery results for IHT, Theorem~\ref{recov1noise} proves a phase transition bound that is equally valid over a continuous stepsize range. In contrast, the recovery results in~\cite{threshRIP,gargkandekhar,foucart,greedy} either require a specific fixed stepsize or degrade with the choice of stepsize.

In the absence of noise, the same condition guarantees exact recovery of the original signal $\xs$.

\begin{corollary}[\textbf{Recovery result for IHT; noiseless case}]\label{recov1noiseless}
Suppose Assumption {\bf A.2} holds, as well as (\ref{IHTstablecond}), and that $\al$ satisfies (\ref{alphabound1})
and  the noise $e\eqdef 0$. Then, in the proportional-growth asymptotic, IHT converges to $\xs$ with probability tending to $1$ exponentially in $n$.
\end{corollary}

In the case of IHT applied to problems with zero noise, the above result has a surprising corollary: a condition can be given which guarantees that, with overwhelming probability, the underlying $k$-sparse signal $\xs$ is the algorithm's only fixed point. In other words, within some portion of phase space, there is only one possible solution to which the IHT algorithm can converge, namely the underlying signal $\xs$. This is remarkable since IHT is a gradient projection algorithm for the nonconvex problem (\ref{c:l0problem}) which can be shown to have a combinatorially large number of local minimizers. The conclusion is that the properties of Gaussian matrices ensure that, within this region of phase space, the IHT algorithm will never `get stuck' at an unwanted local minimizer, thus exhibiting a behaviour one would usually only expect if a convex problem was being solved. The result follows.

\begin{corollary}[\textbf{Single fixed point condition; noiseless case}]\label{theta1cor}
Suppose Assumption {\bf A.2} holds, as well as (\ref{SSPcond}), 
and that $e=0$. Then, in the proportional-growth asymptotic, $\xs$ is the only fixed point of IHT with stepsize $\al$, with probability tending to $1$ exponentially in $n$.
\end{corollary}

\subsubsection{Results for N-IHT}\label{NIHT_proofs}

Again, we first define two numerically computable functions, namely, the lower bound $\rh^{N-IHT}(\dd)$ on the phase transition for recovery using N-IHT and
the stability factor  $\xi(\dd,\rr)$.

\begin{definition}[\textbf{Phase transition lower bound for N-IHT}]\label{rhIHT1}
Given $\dd\in(0,1]$, define the phase transition lower bound $\rh^{N-IHT}(\dd)$ to be the unique solution to
\begin{equation}\label{rhoNIHTdef}
\frac{\sqrt{\IFF(\dd,\rr)}}{(1-\rr)\left[1-\IL(\dd,\rr,1-\rr)\right]}=\frac{1}{\kappa[1+\UU(\dd,2\rr)]}\;\;\;\;\mbox{for}\;\;\;\;\rr\in(0,1/2],
\end{equation}
where $\IFF$ is defined in (\ref{Fdef}), $\IL$ is defined in (\ref{ldef}), $\UU$ is defined in~\cite[Definition 2.2]{BT}, and $\kappa$ is an N-IHT algorithm parameter.\footnote{A proof that $\rh^{N-IHT}(\dd)$ is well-defined can be found in \cite[Section 5.2]{thesis}.}
\end{definition}

\begin{definition}[\textbf{Stability factor for N-IHT}]\label{massive_def2}
Given $\dd\in(0,1]$ and $\rr\in(0,1/2]$, provided $\rr<\rh^{N-IHT}(\dd)$ holds, where $\rh^{N-IHT}(\dd)$ is defined in (\ref{rhoNIHTdef}), 
define the stability factor $\xi(\dd,\rr)$ to be
\begin{equation}\label{xidef2}
\xi(\dd,\rr)\eqdef\sqrt{\IFF(\dd,\rr)\left[a(\dd,\rr)\right]^2+\left[a(\dd,\rr)\right]^2},
\end{equation}
where\footnote{Note that (\ref{NIHTstablecond}) ensures that the denominator in (\ref{adefn2}) is strictly positive and that $a(\dd,\rr)$ is therefore well-defined. The reader may verify by comparison with Definitions~\ref{massive_def} that the $\al$ terms have been replaced by the term $\{\kappa[1+\UU(\dd,2\rr)]\}^{-1}$.}  
\begin{equation}\label{adefn2}
a(\dd,\rr)\eqdef\frac{\sqrt{\IFF(\dd,\rr)}+\{\kappa[1+\UU(\dd,2\rr)]\}^{-1}\sqrt{\rr(1-\rr)[1+\IU(\dd,\rr,1-\rr)][1+\IU(\dd,\rr,\rr)]}}{(1-\rr)\{\kappa[1+\UU(\dd,2\rr)]\}^{-1}[1-\IL(\dd,\rr,1-\rr)]-\sqrt{\IFF(\dd,\rr)}}.
\end{equation} and where 
$\IFF$ is defined in (\ref{Fdef}), $\IU$ is defined in (\ref{udef}), $\IL$ is defined in (\ref{ldef}), $\UU$ is defined in~\cite[Definition 2.2]{BT}.
\end{definition}

We have the following recovery result for N-IHT.

\begin{theorem}[\textbf{Recovery result for N-IHT; noise case}]\label{recov2noise}
Suppose Assumptions {\bf A.2} and {\bf A.3} hold, as well as
\begin{equation}\label{NIHTstablecond}
\rr<\rh^{N-IHT}(\dd),
\end{equation}
where $\rh^{N-IHT}$ is defined in (\ref{rhoNIHTdef}). Then, in the proportional-growth asymptotic, N-IHT converges to $\xb$ such that
\begin{equation}\label{error2}
\|\xb-\xs\|\le\xi(\dd,\rr)\cdot\Sg,
\end{equation}
with probability tending to $1$ exponentially in $n$, where $\xi(\dd,\rr)$ is defined in \req{xidef2}.
\end{theorem}

In the case of zero noise, Theorem~\ref{recov2noise} also simplifies to an exact recovery result.

\begin{corollary}[\textbf{Recovery result for N-IHT; noiseless case}]\label{recov2noiseless}
Suppose Assumption {\bf A.2} holds, as well as (\ref{NIHTstablecond}), and that the noise
$e\eqdef 0$. Then, in the proportional-growth asymptotic, N-IHT converges to $\xs$ with probability tending to $1$ exponentially in $n$.
\end{corollary}

\subsection{Illustration and discussion of results}\label{discussion}

\begin{figure}[h]
\centering
\subfigure[]{\includegraphics[width=3in,height=2.4in]{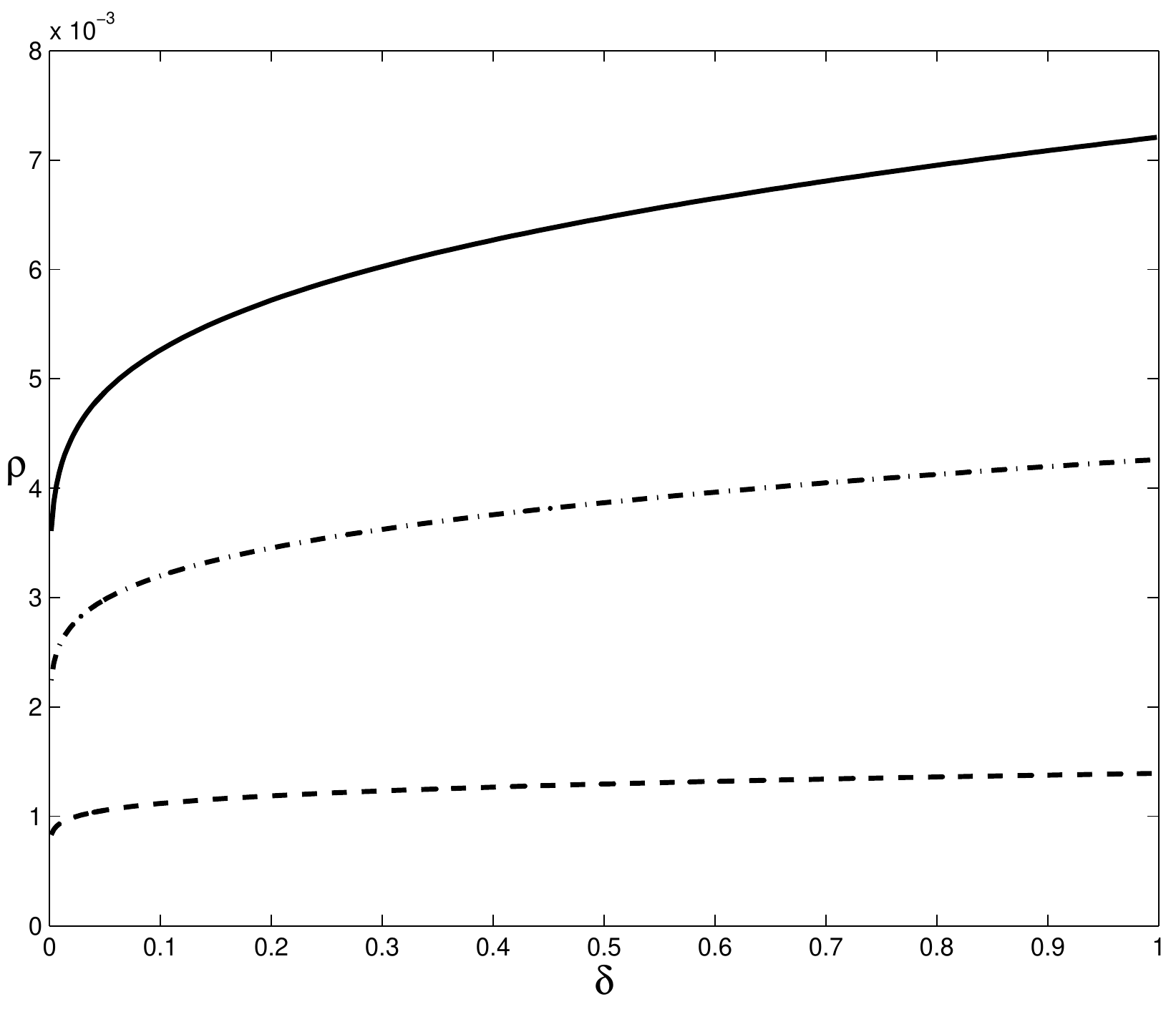}}
\subfigure[]{\includegraphics[width=3in,height=2.4in]{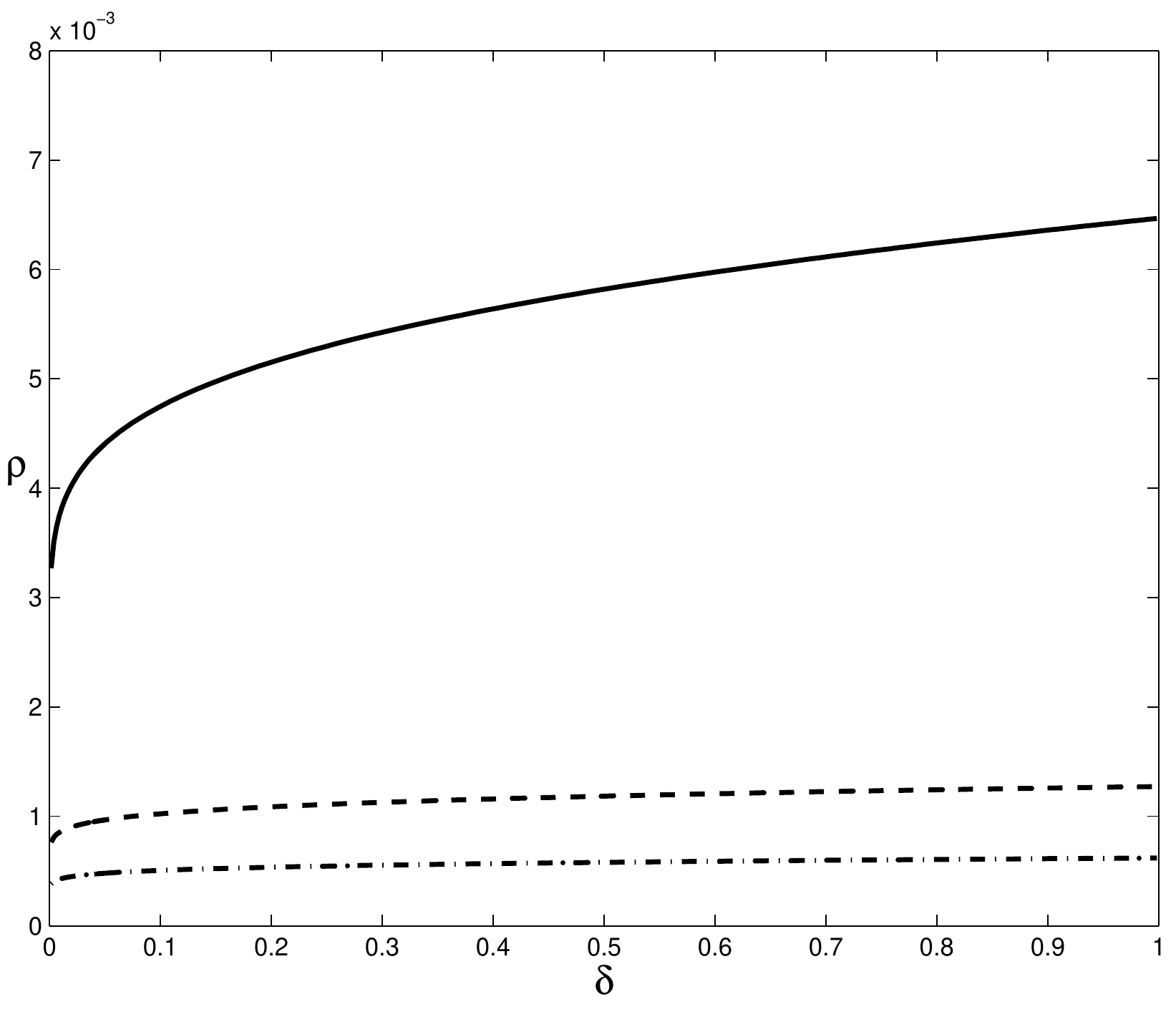}}
\caption{Our average-case phase transition bounds for IHT algorithms (unbroken) compared with the best-known RIP-based phase transition bounds based on our stable point analysis~\cite{thesis} (dashed) and the analysis in~\cite{HTP} (dash-dot): (a) IHT (b) N-IHT.}
\label{phaseplots}
\end{figure}

\begin{figure}[h]
\centering
\subfigure[]{\includegraphics[width=3in,height=2.4in]{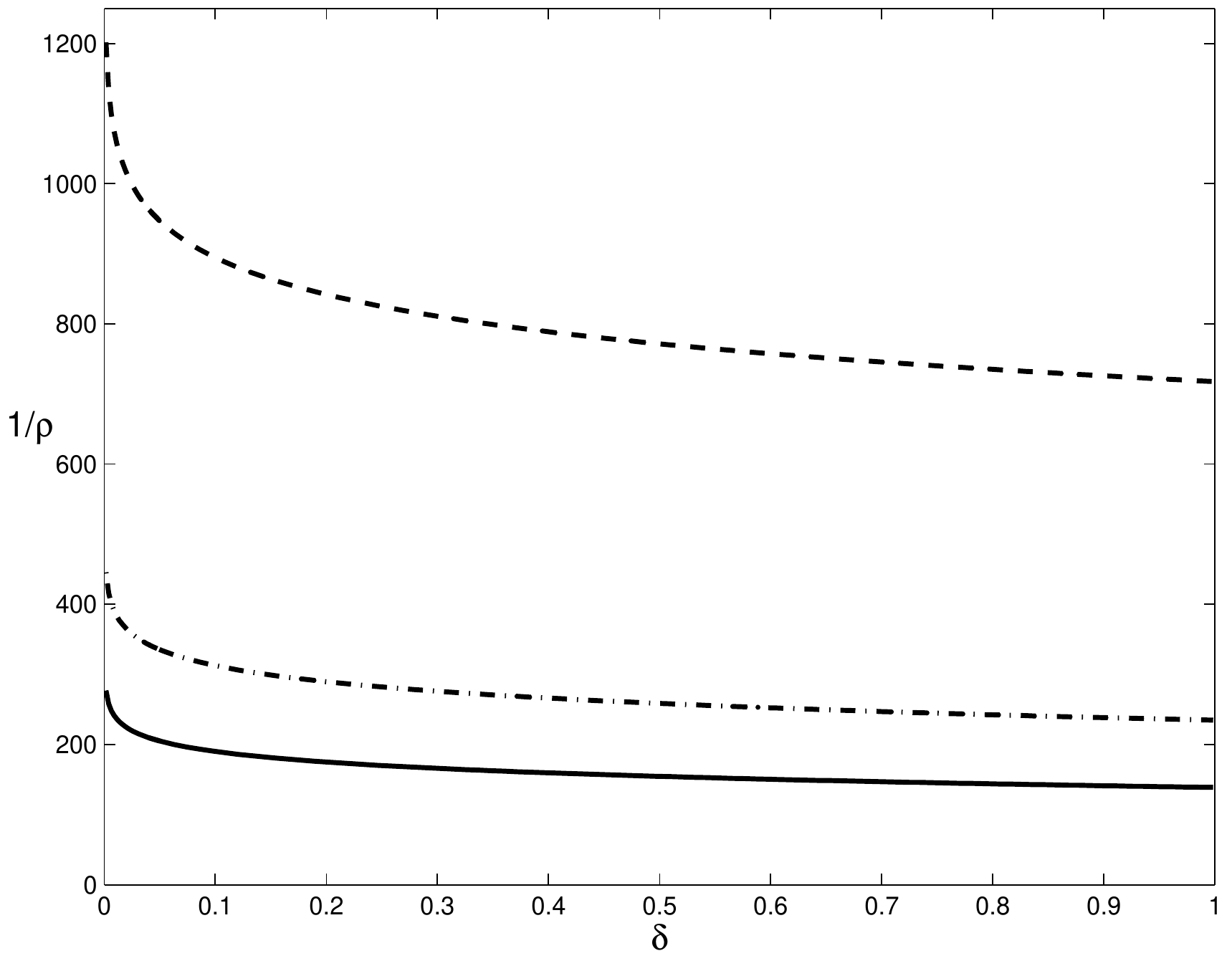}}
\subfigure[]{\includegraphics[width=3in,height=2.4in]{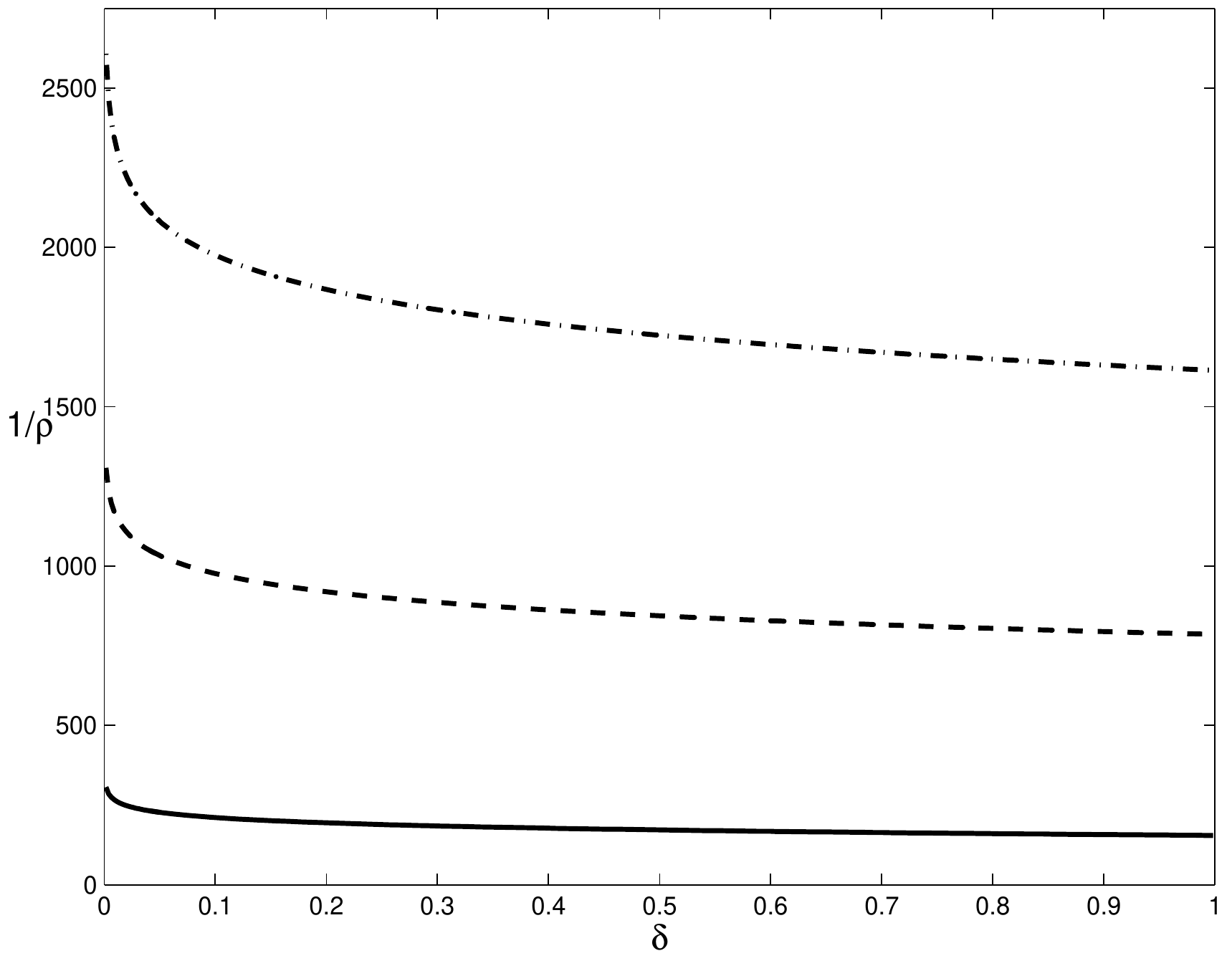}}
\caption{Inverse of the phase transition bounds in Figure~\ref{phaseplots}: (a) IHT (b) N-IHT.}
\label{invphaseplots}
\end{figure}

{\bf Noiseless case.}\quad The recovery phase transition bounds given in Definition~\ref{rhIHT} for IHT and N-IHT (with $\kappa=1.1$) respectively are displayed in Figure~\ref{phaseplots}. Exact recovery in the case of zero noise is guaranteed asymptotically for $(\dd,\rr)$ pairs falling below the respective curves. The best-known lower bounds on exact recovery phase transitions obtained in~\cite{thesis} from previous RIP analysis are included for comparison: the IHT phase transition bound applies the RIP bounds in~\cite{BT} to Foucart's analysis in~\cite{HTP}, while an extension of the same approach leads to the phase transition bound for N-IHT. An RIP analysis of the stable point approach adopted in this paper was also carried out in~\cite{thesis}, and the resulting phase transition bounds are also displayed in Figure~\ref{phaseplots}. We see a considerable improvement over the phase transition bounds corresponding to previous RIP analysis, with recovery being guaranteed for IHT for values of $\rr$ around $1.7$ times higher than before, and for N-IHT around $10$ times higher than before. Figure~\ref{invphaseplots} displays the inverse of the phase transition bound for each stepsize scheme. Previous RIP analysis requires a lower bound of $n\geq 234k$ measurements to guarantee recovery using IHT, and $n\geq 1617k$ using N-IHT. By comparison, we reduce these lower bounds to $n\geq 138k$ for IHT and $n\geq 154k$ for N-IHT. It should also be added that our result for IHT holds for a continuous stepsize range, while the result based upon~\cite{HTP}, in keeping with all other similar RIP-based results for IHT (see~\cite{thesis}), holds true only if the stepsize is optimized to a particular value.

{\bf Interpretation of recovery results as lower bounds on a weak phase transition.}\quad We have obtained an improvement by switching to a new method of analysis which allows us to leverage the assumption that the measurements are statistically independent of the signal. The latter has allowed us to make a partial transition from worst-case to average-case analysis.

The distinction between worst-case and average-case phase transitions can also be found in the phase transitions of Donoho and Tanner for recovery using $l_1$-minimization~\cite{precise}, where successful recovery by means of $l_1$-minimization is shown to be equivalent to the neighbourliness of projected $l_1$ balls~\cite{neighborliness}. There are both strong and weak version of neighbourliness: strong neighbourliness guarantees recovery of \textit{any} signal by means of $l_1$-minimization, while weaker forms of neighbourliness assume either a randomly-chosen support and/or randomly-chosen sign pattern for the signal. It is appropriate then to see our results as lower bounds on a weak phase transition for IHT, in contrast to an RIP analysis which gives lower bounds on the strong phase transition. The notions of weakness are comparable but not identical: in the case of $l_1$-minimization, some dependency between the signal and measurement matrix is permitted: it is only required that the support set and sign pattern of the signal are chosen independently of the matrix. However, independence is the only assumption we place upon the signal, and beyond this there is no further restriction upon the signal's coefficients.

 It is worth pointing out that it is the weak phase transition that is observed empirically for recovery by means of $l_1$-minimization, and the same is also to be expected for IHT algorithms. While we obtain a significant improvement, our lower bound is still pessimistic compared to the weak phase transition observed empirically, though we have succeeded in narrowing the gap between the two. It is no surprise that our results do not give the precise weak phase transition, due to the continued (but limited) use of worst-case techniques, such as the RIP and large deviations analysis. However, the use of the average-case independence assumption to analyse the stable point condition has allowed us to break free in part from the restrictions of worst-case analysis.

\begin{figure}[h]
\centering
\includegraphics[height=3in,width=5in]{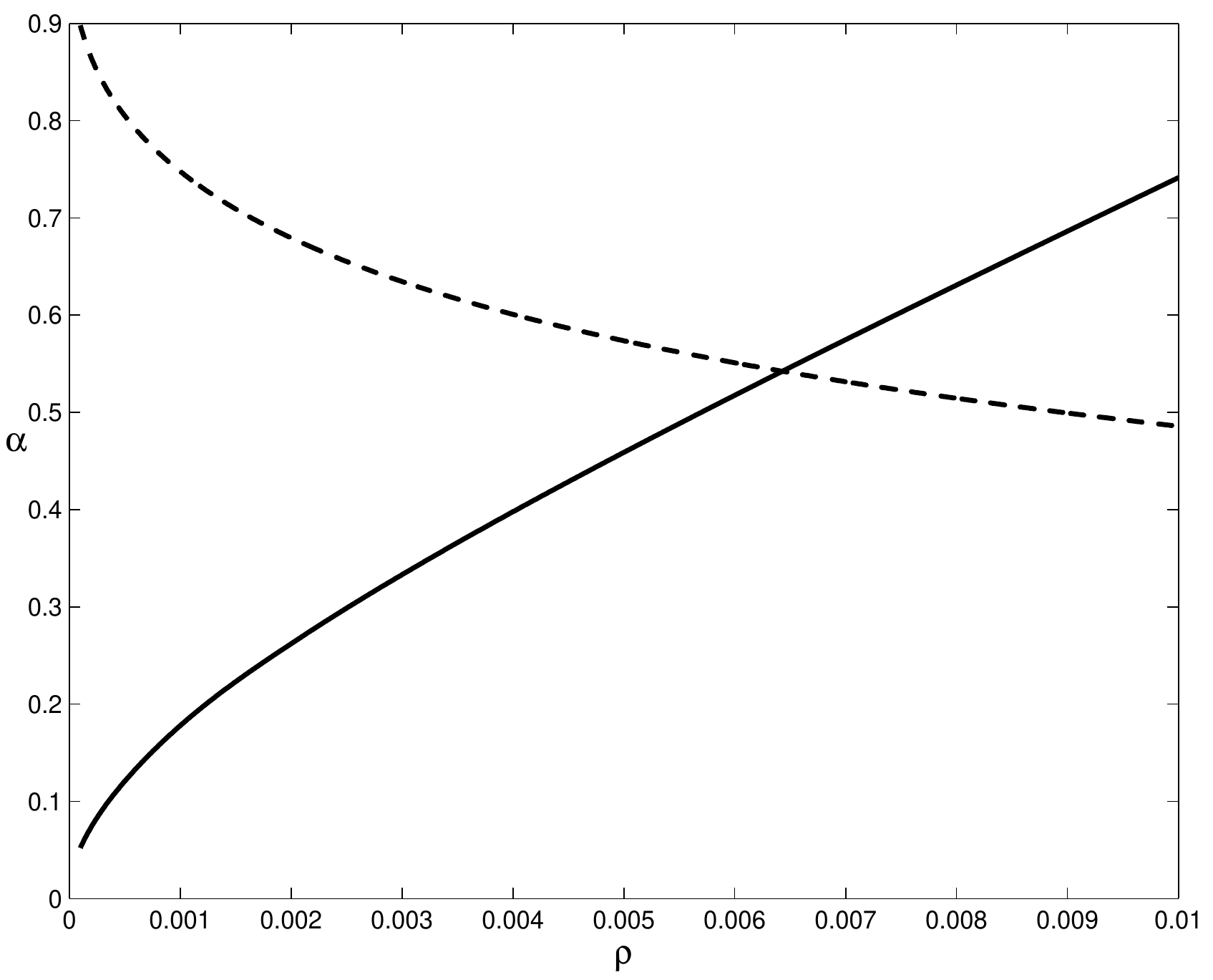}
\caption{Lower bound (unbroken) and upper bound (dashed) on the IHT stepsize for $\dd=0.5$.}
\label{deltahalf}
\end{figure}

{\bf Choice of stepsize for IHT.}\quad Corollary~\ref{recov1noiseless} guarantees exact recovery using IHT provided the stepsize $\al$ falls within the interval given in (\ref{alphabound1}), provided this interval is well-defined. In fact, an inspection of the proof of these two results reveals that the lower bound in (\ref{alphabound1}) arises from the stable point condition, while the upper bound in (\ref{alphabound1}) arises from the convergence condition. Figure~\ref{deltahalf} illustrates these bounds for the case $\dd=0.5$. We see that, as $\rr$ is increased, the admissible stepsize range contracts, until a critical $\rr$-value is reached at which the interval is no longer well-defined.

It has been observed empirically~\cite{optimaltuning} that care must be taken to ensure that the IHT stepsize is neither too small or too large. Our analysis gives theoretical insight into this observation: the stepsize must be small enough to ensure that the algorithm converges, but large enough to ensure that it does not converge to fixed points other than the underlying sparse signal.

\begin{figure}[h]
\centering
\subfigure[]{\includegraphics[width=3in,height=2.4in]{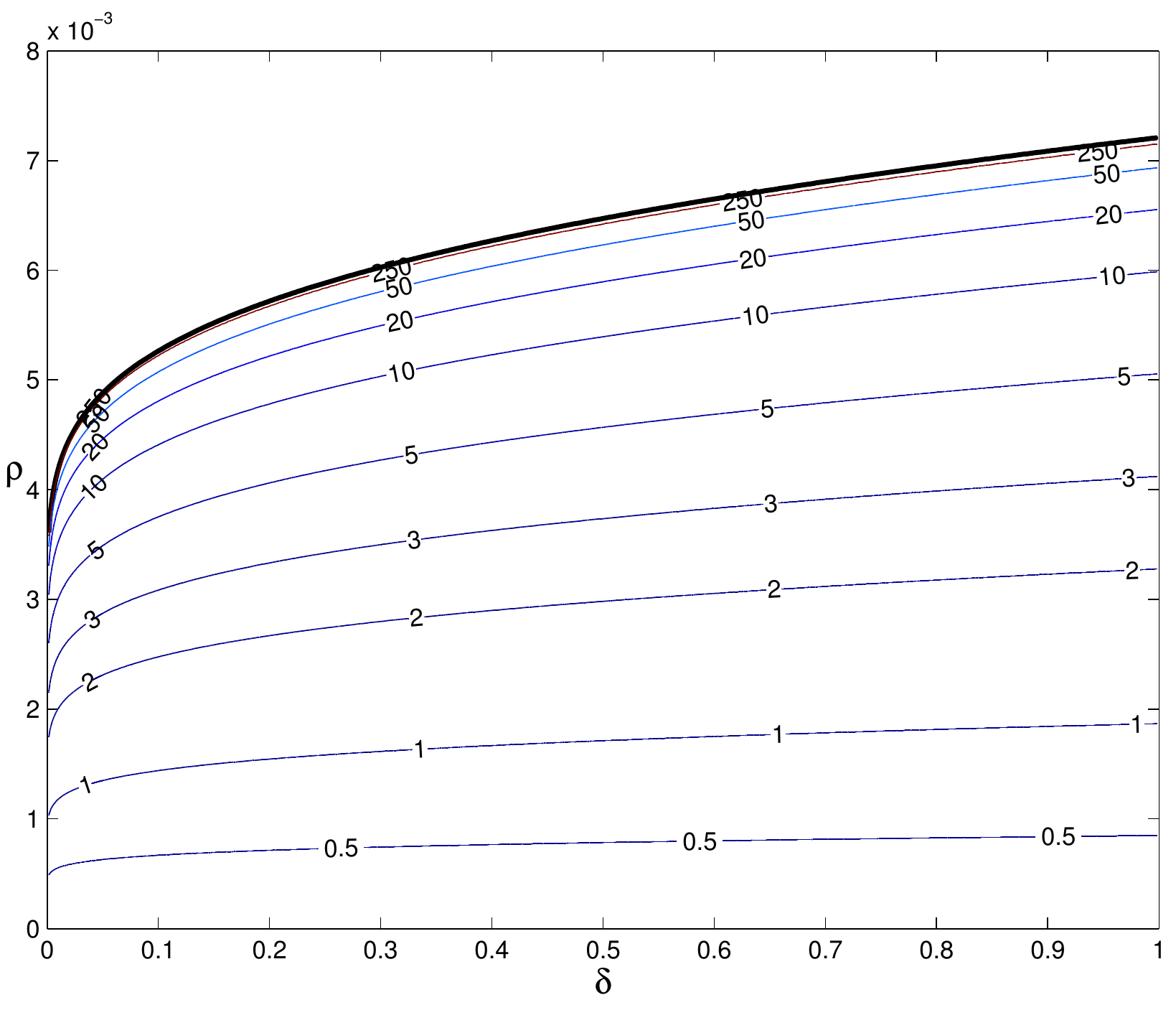}}
\subfigure[]{\includegraphics[width=3in,height=2.4in]{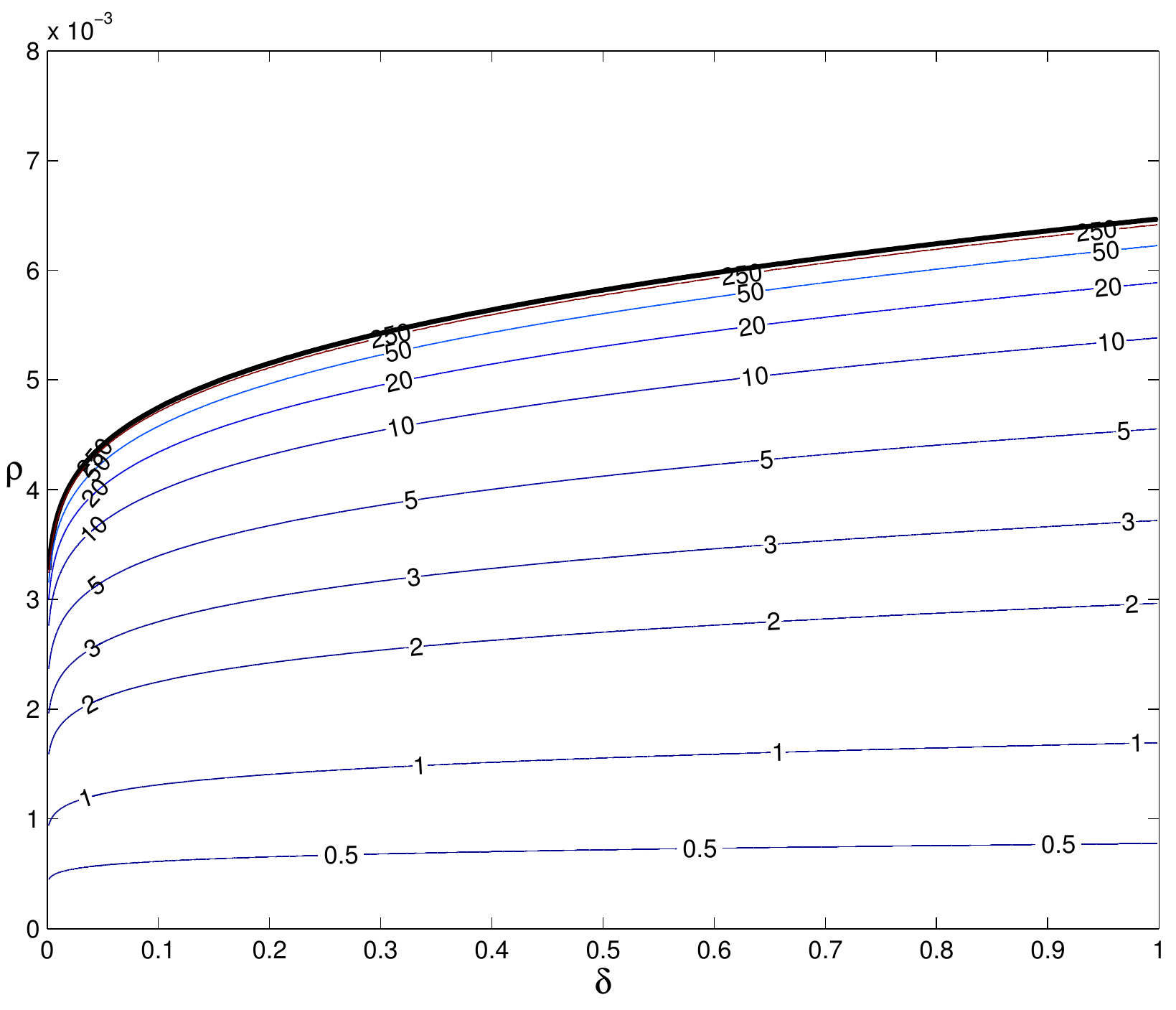}}
\caption{Plot of the stability factor $\xi(\dd,\rr)$ for (a) IHT (b) N-IHT.}
\label{stabilityplots}
\end{figure}

{\bf Extension to noise.}\quad In the case where measurements are contaminated by noise, exact recovery of the original signal is impossible. However, Theorems~\ref{recov1noise} and~\ref{recov2noise} guarantee that, in the same region of phase space defined by the exact recovery phase transition bound, the approximation error of the output of IHT/N-IHT is asymptotically bounded by some known stability factor multiplied by the noise level $\sg$. Figure~\ref{stabilityplots} plots this noise stability factor $\xi(\dd,\rr)$ for each of the two stepsize schemes considered ($\kappa=1.1$ for N-IHT). In keeping with the results in~\cite{lqphase} and~\cite{greedy}, we observe that the stability factor tends to infinity as the transition point is reached.

For both IHT and N-IHT, in the region for which the stability factors derived in this paper are defined, they are everywhere lower than the corresponding stability factors derived from the previous analysis in~\cite{HTP}; see \cite{thesis} for a comparison. It should be pointed out that we have obtained improved stability results by imposing additional restrictions upon the noise, namely that the noise is Gaussian distributed and independent of the signal and measurement matrix. This assumption is in keeping with our aim of using average-case assumptions. Our analysis could, however, be altered to deal with the case of non-independent noise by making more use of the RIP, though this would lead to larger stability constants.

We have also extended our analysis in~\cite{thesis} to the case of signals which are only approximately $k$-sparse, for both IHT and N-IHT, though we omit this extension in the present work for the sake of brevity. In this extension, a stability factor is derived which multiplies the unrecoverable energy of the signal, due to both measurement noise and inaccuracy of the $k$-sparse model.

\subsection{Proof of recovery results for IHT}\label{IHT_proofs}

{\bf Roadmap for the results in this section.}
Here we prove the results stated in Section \ref{IHTresults_Section}. 
Let us outline how our argument will proceed. We first define a support set partition, see (\ref{thetadefn}) that follows. This partition 
is defined in such a way that, provided (\ref{SSPcond}) holds, an analysis of the stable point condition (\ref{stablecond_noise}) shows that there are asymptotically no $\al$-stable points on any $\Gm_i$ such that $i\in\Ta$, and this is proved in Lemma~\ref{theta1lem}. On the other hand, it is also possible to use the large deviations results of Section~\ref{largedev} to bound the error in approximating $\xs$ by any $\al$-stable point on $\Gm_i$ such that $i\in\Tb$, which is achieved by Lemma~\ref{theta2lem}. It follows that, for any $\al>0$, all $\al$-stable points have bounded approximation error. Finally, Lemma~\ref{stableconvlem} builds on the convergence result in Theorem~\ref{conv1} and gives a condition on the stepsize $\al$ which asymptotically guarantees convergence of IHT to some $\al$-stable point. Combining all three results, we have convergence to some $\al$-stable point with guaranteed approximation error, provided the conditions in each lemma hold; combining the conditions leads to the phase transition bound defined in (\ref{rhIHT}).

We begin by defining the above-mentioned support set partition. 

\begin{definition}[\textbf{Support set partition for IHT}]\label{zeta_def}
Suppose $\dd\in(0,1]$, $\rr\in(0,1/2]$ and $\al>0$. Given $\zeta>0$, let us write
\begin{equation}\label{astardef}
\astar\eqdef a(\dd,\rr)+\zeta,
\end{equation}
where $a(\dd,\rr)$ is defined in (\ref{massive_def}), let us write $\{\Gm_i:i\in S_n\}$ for the set of all possible support sets of cardinality $k$, and let us disjointly partition $S_n\eqdef\Ta\cup\Tb$ such that
\begin{equation}\label{thetadefn}
\Ta\eqdef\left\{i\in S_n\;\;:\;\;\|\xs_{\Lm\sm\Gm_i}\|>\Sg\cdot\astar\right\};\;\;\;\;\Tb\eqdef\left\{i\in S_n\;\;:\;\;\|\xs_{\Lm\sm\Gm_i}\|\le\Sg\cdot\astar\right\}.
\end{equation}
\end{definition}

We recall that $\Lm$ is defined to be the support of the original signal $\xs$. Note that the partition $S_n:=\Ta\cup\Tb$ defined in (\ref{thetadefn}) also depends on $\zeta$, though we omit this dependency from our notation for the sake of brevity. Note also that if $\Gm_i=\Lm$, then $\|\xs_{\Lm\sm\Gm_i}\|=0$ and $i\in\Tb$. In other words, the index corresponding to $\Lm$ is contained in $\Tb$.

We first show that, asymptotically, there are no $\al$-stable points on any $\Gm_i$ with $i\in\Ta$, and we write $NSP_{\al}$ for this event.

\begin{lemma}\label{theta1lem}
Choose $\zeta>0$. Suppose Assumptions {\bf A.2} and {\bf A.3} hold, as well as (\ref{SSPcond}). Then, in the proportional-growth asymptotic, there are no $\al$-stable points on any $\Gm_i$ such that $i\in\Ta$, with probability tending to $1$ exponentially in $n$.
\end{lemma}

\proof{For any $\Gm_i$ such that $i\in\Ta$, we have $\Gm_i\neq\Lm$, and we may therefore use Theorem~\ref{singlefp} and Lemma~\ref{dist} with $\Gm:=\Gm_i$ to deduce that a necessary condition for there to be an $\al$-stable point on $\Gm_i$ is
\begin{equation}\label{stablecond2}
\|\xs_{\Lm\sm\Gm_i}\|\cdot\sqrt{F_{\Gm_i}}+\sg\cdot\sqrt{G_{\Gm_i}}\geq\al\left[\left(\frac{n-k}{n}\right)\|\xs_{\Lm\sm\Gm_i}\|\cdot R_{\Gm_i}-\sg\cdot\sqrt{\frac{k(n-k)}{n^2}\cdot S_{\Gm_i}\cdot T_{\Gm_i}}\right],
\end{equation}
where
$$\begin{array}{c}F_{\Gm_i}\sim\displaystyle\frac{k}{n-k+1}\mathcal{F}(k,n-k+1);\;\;\;\;G_{\Gm_i}\sim\displaystyle\frac{k}{n-k+1}\mathcal{F}(k,n-k+1);\\
R_{\Gm_i}\sim\displaystyle\frac{1}{n-k}\chi^2_{n-k};\;\;\;\;S_{\Gm_i}\sim\displaystyle\frac{1}{n-k}\chi^2_{n-k};\;\;\;\;T_{\Gm_i}\sim\displaystyle\frac{1}{k}\chi^2_{k}.\end{array}$$
We also have, by (\ref{thetadefn}),
\begin{equation}\label{Sgbound}
\Sg\le\frac{\|\xs_{\Lm\sm\Gm_i}\|}{\astar}
\end{equation}
for any $\Gm_i$ such that $i\in\Ta$. Since $\Gm_i\neq\Lm$, $\|\xs_{\Lm\sm\Gm}\|>0$, and substitution of (\ref{Sgbound}) into (\ref{stablecond2}), rearrangement and division by $\|\xs_{\Lm\sm\Gm_i}\|$ yields
$$\astar\left[\al\left(\frac{n-k}{n}\right)\cdot R_{\Gm_i}-\sqrt{F_{\Gm_i}}\right]\le\sqrt{G_{\Gm_i}}+\al\sqrt{\frac{k(n-k)}{n^2}\cdot S_{\Gm_i}\cdot T_{\Gm_i}}.$$
Consequently,
\begin{eqnarray}\label{necrhon}
&&\PP(\bNa)\nonumber\\
&=&\PP\left\{\cup_{i\in\Ta}(\exists\;\mbox{an $\al$-stable point
  supported on}\;\Gm_i)\right\}\nonumber\\
&\le&\PP\left\{\bigcup_{i\in\Ta}\left[\astar\left[\al\left(1-\rr_n\right)\cdot R_{\Gm_i}-\sqrt{F_{\Gm_i}}\right]\le\sqrt{G_{\Gm_i}}+\al\sqrt{\rr_n(1-\rr_n)\cdot S_{\Gm_i}\cdot T_{\Gm_i}}\right]\right\}\nonumber,\\
&&
\end{eqnarray}
where we write $\rr_n$ for the sequence of values of the ratio $k/n$. For brevity's sake, let us define
\begin{equation}\label{fdef}
\Phi[\rr,F,G,R,S,T]\eqdef\sqrt{G}+\al\sqrt{\rr(1-\rr)(S)(T)}-\astar\cdot\left[\al(1-\rr)\cdot R-\sqrt{F}\right],
\end{equation}
so that (\ref{necrhon}) may be equivalently written as
\begin{equation}\label{necsplit0}
\PP(\bNa)\le\PP\left\{\cup_{i\in\Ta}\left(\Phi[\rr_n,F_{\Gm_i},G_{\Gm_i},R_{\Gm_i},S_{\Gm_i},T_{\Gm_i}]\geq
  0\right)\right\}.
\end{equation}
Given some $\e>0$, we now define
\begin{equation}\label{astdef}
\begin{array}{c}
F^{\ast}=G^{\ast}\eqdef\IFF(\dd,\rr)+\e;\;\;\;\;\;R^{\ast}\eqdef1-\IL(\dd,\rr,1-\rr)-\e;\\
S^{\ast}\eqdef1+\IU(\dd,\rr,1-\rr)+\e;\;\;\;\;\;T^{\ast}\eqdef1+\IU(\dd,\rr,\rr)+\e.\end{array}
\end{equation}
Using (\ref{astdef}), we deduce from (\ref{necsplit0}) that
\begin{eqnarray}
&&\PP(\bNa)\nonumber\\
&\le&\PP\left\{\cup_{i\in\Ta}\left(\Phi[\rr_n,F_{\Gm_i},G_{\Gm_i},R_{\Gm_i},S_{\Gm_i},T_{\Gm_i}]\geq \Phi[\rr_n,F^{\ast},G^{\ast},R^{\ast},S^{\ast},T^{\ast}]\right)\right\}\label{necsplitnoise1}\\
&+&\PP\left\{\Phi[\rr_n,F^{\ast},G^{\ast},R^{\ast},S^{\ast},T^{\ast}]\geq \Phi[\rr,F^{\ast},G^{\ast},R^{\ast},S^{\ast},T^{\ast}]+\e\right\}\label{necsplitnoise2}\\
&+&\PP\left\{\Phi[\rr,F^{\ast},G^{\ast},R^{\ast},S^{\ast},T^{\ast}]+\e\geq
  0\right\}\label{necsplitnoise3},
\end{eqnarray}
since the event in the right-hand side of (\ref{necsplit0}) lies in the union of the three
events in (\ref{necsplitnoise1}), (\ref{necsplitnoise2}) and
(\ref{necsplitnoise3}). Now (\ref{necsplitnoise3}) is a deterministic event, and $\astar$ has been defined in such a way that, for any $\zeta>0$, provided $\e$ is taken sufficiently small, the event has probability
$0$. This follows from (\ref{SSPcond}), (\ref{adefn}), (\ref{astardef}), and by the continuity of $\Phi$. The event (\ref{necsplitnoise2}) is also deterministic, and by continuity and since $\rr_n\rightarrow\rr$, it follows that there exists some $\tilde{n}$ such that
$$\PP\left\{\Phi[\rr_n,F^{\ast},G^{\ast},R^{\ast},S^{\ast},T^{\ast}]\geq \Phi[\rr,F^{\ast},G^{\ast},R^{\ast},S^{\ast},T^{\ast}]+\e\right\}=0\;\;\;\;\mbox{for all}\;\;n\geq\tilde{n}.$$
Taking limits as $n\rightarrow\infty$, the terms (\ref{necsplitnoise2}) and
(\ref{necsplitnoise3}) are zero, leaving only (\ref{necsplitnoise1}), and we have
\begin{eqnarray}
&&\lim_{n\rightarrow\infty}\PP(\bNa)\nonumber\\
&\le&\lim_{n\rightarrow\infty}\PP\left\{\cup_{i\in\Ta}\left(\Phi[\rr_n,F_{\Gm_i},G_{\Gm_i},R_{\Gm_i},S_{\Gm_i},T_{\Gm_i}]\geq \Phi[\rr_n,F^{\ast},G^{\ast},R^{\ast},S^{\ast},T^{\ast}]\right)\right\}\nonumber\\
&\le&\lim_{n\rightarrow\infty}\PP\left\{\cup_{i\in\Ta}(F_{\Gm_i}\geq F^{\ast})\right\}+\lim_{n\rightarrow\infty}\PP\left\{\cup_{i\in\Ta}(G_{\Gm_i}\geq G^{\ast})\right\}+\lim_{n\rightarrow\infty}\PP\left\{\cup_{i\in\Ta}(R_{\Gm_i}\le R^{\ast})\right\}\nonumber\\
&+&\lim_{n\rightarrow\infty}\PP\left\{\cup_{i\in\Ta}(S_{\Gm_i}\geq S^{\ast})\right\}+\lim_{n\rightarrow\infty}\PP\left\{\cup_{i\in\Ta}(T_{\Gm_i}\geq T^{\ast})\right\},\label{neclemnoise}
\end{eqnarray}
where the last line follows from the monotonicity of $\Phi$ with respect to $F$, $G$, $R$, $S$ and $T$. Since
  $\Ta\subseteq S_n$, we may apply
  Lemmas~\ref{chisq} and~\ref{Fdist} to (\ref{neclemnoise}), and we deduce $\PP(\bNa)\ra 0$ as $n\ra\infty$, exponentially in $n$, as required.}

Next, we show that any $\al$-stable points on $\Gm_i$ with $i\in\Tb$ are `close' to $\xs$.

\begin{lemma}\label{theta2lem}
Suppose Assumptions {\bf A.2} and {\bf A.3} hold, as well as (\ref{SSPcond}). Then there exists $\zeta$ sufficiently small such that, in the proportional-growth asymptotic, any $\al$-stable point $\xb$ on $\Gm_i$ such that $i\in\Tb$ satisfies (\ref{error}) with probability tending to $1$ exponentially in $n$, where $\xi(\dd,\rr)$ is defined in (\ref{xidef}).
\end{lemma}

\proof{If $\sg=0$, the result follows trivially from Lemma~\ref{theta1cor}, so let us assume that $\sg>0$. Suppose $\xb$ is a minimum-norm solution on $\Gm$, so that $\xb_{\Gm}=A_{\Gm}^{\dag}b$. Then, using $A_{\Gm}^{\dag}A_{\Gm}=I$, we have
\begin{eqnarray}
(\xb-\xs)_{\Gm}&=&A_{\Gm}^{\dag}(A_{\Gm}\xs_{\Gm}+A_{\Gm^C}\xs_{\Gm^C}+e)-\xs_{\Gm}\nonumber\\
&=&\xs_{\Gm}+A_{\Gm}^{\dag}(A_{\Lm\sm\Gm}\xs_{\Lm\sm\Gm}+A_{(\Lm\cup\Gm)^C}\xs_{(\Lm\cup\Gm)^C}+e)-\xs_{\Gm}\nonumber\\
&=&A_{\Gm}^{\dag}(A_{\Lm\sm\Gm}\xs_{\Lm\sm\Gm}+\et)+\xs_{\Gm}-\xs_{\Gm}\nonumber\\
&=&A_{\Gm}^{\dag}(A_{\Lm\sm\Gm}\xs_{\Lm\sm\Gm}+\et),\label{errorgam}
\end{eqnarray}
while
\begin{equation}\label{errorgamC}
(\xb-\xs)_{\Gm^C}=-\xs_{\Gm^C}.
\end{equation}
Combining (\ref{errorgam}) and (\ref{errorgamC}) using the triangle inequality, we may bound
\begin{eqnarray}
\|\xb-\xs\|^2&=&\|(\xb-\xs)_{\Gm}\|^2+\|(\xb-\xs)_{\Gm^C}\|^2\nonumber\\
&=&\|A_{\Gm}^{\dag}(A_{\Lm\sm\Gm}\xs_{\Lm\sm\Gm}+\et)\|^2+\|\xs_{\Gm^C}\|^2\nonumber\\
&\le&\left[\|A_{\Gm}^{\dag}A_{\Lm\sm\Gm}\xs_{\Lm\sm\Gm}\|+\|A_{\Gm}^{\dag}\et\|\right]^2+\|\xs_{\Lm\sm\Gm}\|^2\label{errorbound1}
\end{eqnarray}
We may deduce, by (\ref{lhs}) of Lemma~\ref{dist},
\begin{equation}\label{errorbound2}
\|A_{\Gm}^{\dag}A_{\Lm\sm\Gm}\xs_{\Lm\sm\Gm}\|^2=\|\xs_{\Lm\sm\Gm}\|^2\cdot P_{\Gm},\;\;\;\mbox{where}\;\;\;P_{\Gm}\sim\frac{k}{n-k+1}\FF(k,n-k+1),
\end{equation}
and by (\ref{lhsnoise}) of Lemma~\ref{dist},
\begin{equation}\label{errorbound3}
\|A_{\Gm}^{\dag}e\|^2=\sg^2\cdot Q_{\Gm},\;\;\;\mbox{where}\;\;\;Q_{\Gm}\sim\frac{k}{n-k+1}\FF(k,n-k+1).
\end{equation}
Substituting (\ref{errorbound2}) and (\ref{errorbound3})
into (\ref{errorbound1}), we have
\begin{equation}\label{PQ_bound}
\|\xb-\xs\|^2\le\left[\|\xs_{\Lm\sm\Gm}\|\cdot\sqrt{P_{\Gm}}+\sg\cdot\sqrt{Q_{\Gm}}\right]^2+\|\xs_{\Lm\sm\Gm}\|^2,
\end{equation}
and we may use (\ref{thetadefn}) to further deduce
\begin{eqnarray}
\|\xb-\xs\|^2&\le&\sg^2\left[\astar\cdot\sqrt{P_{\Gm}}+\sqrt{Q_{\Gm}}\right]^2+\left[\astar\right]^2\cdot\sg^2\nonumber\\
&=&\sg^2\left\{\left[\astar\cdot\sqrt{P_{\Gm}}+\sqrt{Q_{\Gm}}\right]^2+\left[\astar\right]^2\right\}.\label{errorbound}
\end{eqnarray}
For the sake of brevity, let us define
\begin{equation}\label{psi_def}
\Psi(P,Q):=\sqrt{\left(\astar\cdot\sqrt{P}+\sqrt{Q}\right)^2+\astar^2},
\end{equation}
so that (\ref{errorbound}) may equivalently be written as
\begin{equation}\label{errorbound_brevity}
\|\xb-\xs\|\le\Sg\cdot\Psi\left[P_{\Gm},Q_{\Gm}\right].
\end{equation}
Given $\zeta>0$, let us define
\begin{equation}\label{PQdef}
P^{\ast}=Q^{\ast}:=\IFF(\dd,\rr)+\zeta.
\end{equation}
Now we use (\ref{errorbound_brevity}) to perform a union bound over all $\Gm_i$ such that $i\in\Tb$, writing $\xb_i$ for the minimum-norm solution on $\Gm_i$, giving
\begin{eqnarray}
&&\PP\left\{\exists\;\mbox{some}\;\Gm_i\;\mbox{such
    that}\;i\in\Tb\;\mbox{and}\;\|\xb_i-\xs\|>\Sg\cdot\Psi\left[P^{\ast},Q^{\ast}\right]\right\}\nonumber\\
&=&\PP\left\{\bigcup_{i\in\Tb}\left(\|\xb_i-\xs\|>\Sg\cdot\Psi\left[P^{\ast},Q^{\ast}\right]\right)\right\}\label{2necsplit0}\\
&\le&\PP\left\{\bigcup_{i\in\Tb}\left(\|\xb_i-\xs\|>\Sg\cdot\Psi\left[P_{\Gm_i},Q_{\Gm_i}\right]\right)\right\}\label{2necsplit1}\\
&+&\PP\left\{\bigcup_{i\in\Tb}\left(\Sg\cdot\Psi\left[P_{\Gm_i},Q_{\Gm_i}\right]\geq\Sg\cdot\Psi\left[P^{\ast},Q^{\ast}\right]\right)\right\},\nonumber\\
&&\;\label{2necsplit2}
\end{eqnarray}
since the event in (\ref{2necsplit0}) lies in the union of the two
events in (\ref{2necsplit1}) and (\ref{2necsplit2}). It is an immediate consequence of
(\ref{errorbound_brevity}) that the event in (\ref{2necsplit1}) has probability
$0$. Taking limits of (\ref{2necsplit2}) as $n\rightarrow\infty$, and cancelling $\Sg$, we have
\begin{eqnarray}
&&\lim_{n\rightarrow\infty}\PP\left\{\exists\;\mbox{some}\;\Gm_i\;\mbox{such
    that}\;i\in\Tb\;\mbox{and}\;\|\xb_i-\xs\|>\Sg\cdot\Psi\left[P^{\ast},Q^{\ast}\right]\right\}\nonumber\\
&\le&\lim_{n\rightarrow\infty}\PP\left\{\bigcup_{i\in\Tb}\left(\Psi\left[P_{\Gm_i},Q_{\Gm_i}\right]\geq\Psi\left[P^{\ast},Q^{\ast}\right]\right)\right\}\nonumber\\
&\le&\lim_{n\rightarrow\infty}\PP\left\{\cup_{i\in\Tb}(P_{\Gm_i}\geq P^{\ast})\right\}+\lim_{n\rightarrow\infty}\PP\left\{\cup_{i\in\Tb}(Q_{\Gm_i}\geq Q^{\ast})\right\},\label{2neclem}
\end{eqnarray}
where we used the monotonicity of $\Psi$ with respect to $P$ and $Q$ in the last line. Since
  $\Tb\subseteq S_n$, and using (\ref{errorbound2}) and (\ref{errorbound3}), we may apply
  Lemma~\ref{Fdist} to (\ref{2neclem}), yielding that each of the limits in the right-hand side of (\ref{2neclem}) converges to zero exponentially in $n$, and so finally
$$\lim_{n\rightarrow\infty}\PP\left\{\exists\;\mbox{some}\;\Gm_i\;\mbox{such
    that}\;i\in\Tb\;\mbox{and}\;\|\xb_i-\xs\|>\Sg\cdot\Psi\left[\astar,P^{\ast},Q^{\ast}\right]\right\}=0,$$
exponentially in $n$. Since by Lemma~\ref{necfp}, any stable point is necessarily a minimum-norm solution, and recalling the definition of $\astar$ in (\ref{astardef}), $\Psi(a,P,Q)$ in (\ref{psi_def}), and the definitions of $P^{\ast}$, $Q^{\ast}$ in (\ref{PQdef}), we have
\begin{equation}\label{finalbound}
\lim_{n\rightarrow\infty}\PP\left\{\begin{array}{ll}\exists\;\mbox{some $\al$-stable point $\xb_i$ on}\;\Gm_i\;\mbox{such
    that}\;i\in\Tb\;\mbox{and}\\
\;\|\xb_i-\xs\|>\Sg\sqrt{\IFF(\dd,\rr)\left[1+a(\dd,\rr)+\zeta\right]^2+\left[a(\dd,\rr)+\zeta\right]^2}\end{array}
\right\}=0,
\end{equation}
with convergence exponential in $n$. Finally, by continuity,
$$\begin{array}{l}\|\xb_i-\xs\|>\Sg\sqrt{\IFF(\dd,\rr)\left[1+a(\dd,\rr)\right]^2+1+\left[a(\dd,\rr)\right]^2}\\
\;\;\;\;\;\;\;\;\;\;\;\;\;\Longrightarrow\|\xb_i-\xs\|>\Sg\sqrt{\IFF(\dd,\rr)\left[1+a(\dd,\rr)+\zeta\right]^2+\left[a(\dd,\rr)+\zeta\right]^2},\end{array}$$
for some $\zeta$ suitably small, and the result now follows from the definition of $\xi(\dd,\rr)$ in (\ref{xidef}).}

In the context of IHT, we obtain the following convergence result in the proportional-dimensional asymptotic framework.

\begin{lemma}\label{stableconvlem}
Suppose Assumption {\bf A.2} holds and that the stepsize $\al$ of IHT is chosen to satisfy
\begin{equation}\label{alphaconvbound}
\al<\frac{1}{1+\UU(\dd,2\rr)}.
\end{equation}
Then, in the proportional-growth asymptotic, IHT converges to an $\al$-stable point with probability tending to $1$ exponentially in $n$.
\end{lemma}

\proof{Given (\ref{alphaconvbound}), we may apply Lemma~\ref{RIP_bounds} with $\e$ sufficiently small to deduce $\al(1+U_{2k})<1$, with probability tending to $1$ exponentially in $n$. Under Assumption {\bf A.2}, we may then apply Theorem~\ref{conv1} and deduce convergence of IHT to an $\al$-stable point.}

We now combine Lemmas~\ref{theta1lem},~\ref{theta2lem} and~\ref{stableconvlem} and prove the three main recovery results for IHT.

\textbf{Proof of Theorem~\ref{recov1noise}:} First note that (\ref{IHTstablecond}) implies that the interval in (\ref{alphabound1}) is well-defined. Provided $\al$ is chosen to satisfy (\ref{alphabound1}), (\ref{alphaconvbound}) holds, and under Assumption {\bf A.2}, we may apply Lemma~\ref{stableconvlem} to deduce convergence of IHT to an $\al$-stable point. On the other hand, Lemma~\ref{theta1lem} establishes that there are asymptotically no $\al$-stable points on any $\Gm_i$ such that $i\in\Ta$, while we may apply Lemma~\ref{theta2lem} to deduce that any $\al$-stable points on any $\Gm_i$ such that $i\in\Tb$ satisfy (\ref{error}).

\textbf{Proof of Corollary~\ref{recov1noiseless}:} The result follows by setting $\Sg\eqdef 0$ in Theorem~\ref{recov1noise}.

\textbf{Proof of Corollary~\ref{theta1cor}:} Lemma~\ref{theta1lem} establishes that, if (\ref{SSPcond}) holds, there are asymptotically no $\al$-stable points on any $\Gm_i$ such that $i\in\Ta$. Setting $\Sg\eqdef 0$ in (\ref{thetadefn}), we have $i\in\Tb\Rightarrow\Gm_i=\Lm$. Therefore any $\al$-stable point is supported on $\Lm$, and Lemma~\ref{necfp} implies that it must be $\xs$. However, any fixed point of IHT with stepsize $\al$ is necessarily an $\al$-stable point, and therefore $\xs$ is also the only fixed point of IHT with stepsize $\al$.

\numsection{Conclusions and future directions}

While CS was first developed within the framework of $l_1$-minimization, there is growing evidence that recovery algorithms which do not rely on convex relaxation and the $l_1$-norm can be equally effective in practice~\cite{large_scale}. Two such examples are the gradient-based IHT~\cite{thresh} and N-IHT~\cite{normalized} algorithms, which also have favourable computational efficiency in comparison with other CS approaches. It is important that a CS recovery algorithm is supported by theory which quantitatively determines the degree of undersampling that the algorithm permits. Such results now exist for $l_1$-minimization, where precise phase transitions have been determined within a proportional-growth asymptotic framework in the case of Gaussian matrices~\cite{precise}. By contrast, worst-case recovery guarantees for IHT algorithms using the RIP are pessimistic in comparison with observed empirical behaviour~\cite{greedy}.

To address this issue, we introduced a new method of recovery analysis for IHT algorithms in which we analysed the algorithms' \textit{stable points}, a generalization of the notion of fixed points. By making the realistic assumption of independence between the signal and measurement matrix, we obtained the first recovery guarantees for IHT algorithms and Gaussian measurement matrices which make use of average-case assumptions. In contrast to RIP analysis, which leads to lower bounds on the strong phase transition, we obtained lower bounds on a weak phase transition for recovery using IHT algorithms, which is the notion of practical interest. By breaking free in part from the restrictions of worst-case analysis, we have obtained, to the best of our knowledge, the highest phase transition bounds yet guaranteeing exact recovery of sparse signals by means of IHT and N-IHT. Our results extend to the realistic model of noisy measurements, guaranteeing an improved robustness to these inaccuracies.

The ultimate remaining goal of the work is to fully close the gap between theoretical guarantees and empirical performance for IHT algorithms. At present, the continued use of worst-case methods of analysis such as union bounds over combinatorially many support sets is a hindrance to significant further improvements in phase transition bound. It is an open question whether such a strong requirement is necessary for ensuring signal recovery on average. Though we have obtained quantitative results only for Gaussian matrices here, many other families of random or randomized measurement matrices exhibit similar empirical behaviour and are important to practitioners. Obtaining quantitative guarantees for IHT algorithms applied to such CS measurement schemes is an open avenue of research.

\newappendix{A}
{\Large{\bf Proofs of results in Section \ref{largedev}}}
\vspace*{0.25cm}

We make use of asymptotic results derived by Temme~\cite{temme82} for the incomplete gamma and beta functions, which are related to the $\chi^2$ and $\FF$ distributions respectively. We denote by $P(s,t)$ the lower regularized incomplete gamma function $P(s,t)$~\cite{temme82}, and we let $Q(s,t)=1-P(s,t)$ be the upper regularized incomplete gamma function. We also define the complementary error function erfc$(\omega)$ in the usual way as
$$\mbox{erfc}(\omega)\eqdef\frac{2}{\sqrt{\pi}}\int_{\omega}^{\infty}e^{-u^2}\mathrm{d}u.$$
The result for the gamma function follows.

\begin{lemma}[\textbf{Gamma asymptotic~\cite[Section 3.4 and (2.20)]{temme82}}]\label{temmegamma}
For $0<s<t$,
\begin{equation}\label{Qexpansion}
Q(s,t)=\frac{1}{2}\mbox{erfc}\left(\eta_Q\sqrt{\frac{s}{2}}\right)-R_s(\eta_Q)\;\;\;\;\mbox{where}\;\;\;\;\eta_Q=\sqrt{2\left[\frac{t}{s}-\ln\left(1+\frac{t}{s}\right)\right]},
\end{equation}
and for $s>t>0$,
\begin{equation}\label{Pexpansion}
P(s,t)=\frac{1}{2}\mbox{erfc}\left(-\eta_P\sqrt{\frac{s}{2}}\right)+R_s(\eta_P)\;\;\;\;\mbox{where}\;\;\;\;\eta_P=-\sqrt{2\left[-\frac{t}{s}-\ln\left(1-\frac{t}{s}\right)\right]},
\end{equation}
where $R_s(\cdot)$ is a residual term. Furthermore, if $t/s$ remains fixed so that $\eta$ is held constant,
\begin{equation}\label{gammaresidual}
R_s(\eta)=\mathcal{O}\left(\frac{1}{\sqrt{s}}\right)e^{-\frac{1}{2}s\eta^2}\;\;\;\;\mbox{for $s$ sufficiently large}.
\end{equation}
\end{lemma}

\begin{lemma}\label{chisqlem}
Let $0<l\le n$ and let the random variable $X_l\sim\displaystyle\frac{1}{l}\chi^2_l$. Let $l/n\ra\gm\in(0,1]$ as $n\ra\infty$. Then, for any $\nu>0$,
\begin{equation}\label{loglimit1}
\lim_{n\rightarrow\infty}\frac{1}{n}\ln\PP(X_l^i\geq 1+\nu)=-\frac{\gm}{2}[\nu-\ln(1+\nu)]
\end{equation}
and, for any $\nu\in(0,1)$,
\begin{equation}\label{loglimit2}
\lim_{n\rightarrow\infty}\frac{1}{n}\ln\PP(X_l^i\le 1-\nu)=-\frac{\gm}{2}[-\nu-\ln(1-\nu)].
\end{equation}
\end{lemma}

\proof{We first show (\ref{loglimit1}). We have
\begin{equation}\label{RIGF}
\PP(X_l\geq 1+\nu)=\PP[\chi^2_l\geq l(1+\nu)]=Q\left[\frac{l}{2},\frac{l(1+\nu)}{2}\right],
\end{equation}
where the first step follows from the definition of $X_l$, and the second step follows from the properties of the $\chi^2$ distribution. We can further express the right-hand side of (\ref{RIGF}) by using (\ref{Qexpansion}) with $s=l/2$ and $t=l(1+\nu)/2$, which then gives
\begin{equation}\label{temme82result}
\PP(X_l\geq 1+\nu)=\frac{1}{2}\mbox{erfc}\left(\frac{\eta_Q}{2}\sqrt{l}\right)-R_l(\eta_Q),
\end{equation}
where
\begin{equation}\label{eta}
\eta_Q\eqdef\sqrt{2[\nu-\ln(1+\nu)]}.
\end{equation}
Applying a standard exponential tail bound on the complementary error function $\mbox{erfc}$ to (\ref{temme82result}) then gives
\begin{equation}\label{temme82erfc}
\PP(X_l\geq 1+\nu)\le\frac{1}{2}e^{-\frac{1}{4}l\eta_Q^2}-R_l(\eta_Q),
\end{equation}
to which we can apply (\ref{gammaresidual}) to obtain
$$\PP(X_l\geq 1+\nu)=\mathcal{O}(1)e^{-\frac{1}{4}l\eta_Q^2}\;\;\;\;\mbox{for all $l$ sufficiently large}.$$
Taking logarithms, letting $n\ra\infty$ and recalling that $l/n\ra\gm$, we deduce
$$\displaystyle\lim_{n\rightarrow\infty}\frac{1}{n}\ln\PP(X_l\geq 1+\nu)=\displaystyle\lim_{n\rightarrow\infty}\frac{1}{n}\ln\mathcal{O}(1)+\displaystyle\lim_{n\rightarrow\infty}\frac{1}{n}\cdot\left(-\frac{1}{4}l\eta_Q^2\right)=-\displaystyle\frac{\gm}{4}\eta_Q^2,$$
which together with (\ref{eta}) yields (\ref{loglimit1}). The proof for the lower tail is similar, since the distribution function of $X_l$ is given by
$$\PP(X_l\le 1-\nu)=\PP[\chi^2_l\le l(1-\nu)]=P\left[\frac{l}{2},\frac{l(1-\nu)}{2}\right],$$
which further becomes, due to (\ref{Pexpansion}) with $s=l/2$ and $t=l(1-\nu)/2$,
$$\PP(X_l\le 1-\nu)=\frac{1}{2}\mbox{erfc}\left(-\frac{\eta_P}{2}\sqrt{l}\right)+R_l(\eta_P),$$
where $\eta_P=-\sqrt{2[-\nu-\ln(1-\nu)]}$. The bound (\ref{loglimit2}) now follows similarly to (\ref{loglimit1}).}

We will need the following lemma which gives the limit of a binomial coefficient in the proportional-growth asymptotic.

\begin{lemma}[\textbf{Combinatorial limit}]\label{comblimlem}
In the proportional-dimensional asymptotic,
\begin{equation}\label{comblim}
\lim_{n\ra\infty}\frac{1}{n}\ln{N\choose k}=\frac{H(\dd\rr)}{\dd},
\end{equation}
where $H(\cdot)$ is defined in (\ref{shannon_def}).
\end{lemma}

\proof{In the proportional-dimensional asymptotic,
$$\displaystyle\lim_{n\ra\infty}\frac{1}{n}\ln{N\choose k}=\displaystyle\lim_{n\ra\infty}\frac{N}{n}\cdot\frac{1}{N}\ln{N\choose k}=\frac{1}{\dd}\cdot H(\dd\rr),$$
where the last step follows from Stirling's formula.}

{\bf Proof of Lemma \ref{chisq} (Large deviation result for $\chi^2$).}
Union bounding $\PP\left(X_l^i\geq 1+\nu\right)$ over all $i\in S_n$ gives
\begin{equation}\label{combsum}
\PP\left\{\cup_{i\in S_n}(X_l^i\geq 1+\nu)\right\}\le\sum_{i\in S_n}\PP\left(X_l^i\geq 1+\nu\right)=|S_n|\cdot\PP(X_l^1\geq 1+\nu).
\end{equation}
Taking logarithms and limits of the right-hand side of (\ref{combsum}), using (\ref{loglimit1}) and (\ref{comblim}), we have
$$\displaystyle\lim_{n\ra\infty}\frac{1}{n}\ln\left[|S_n|\cdot\PP(X_l^1\geq 1+\nu)\right]=H(\dd\rr)-\displaystyle\frac{\lm}{2}[\nu-\ln(1+\nu)],$$
and so (\ref{combsum}) implies that, for any $\eta>0$,
\begin{equation}\label{loglimit}
\frac{1}{n}\ln\PP\left\{\cup_{i\in S_n}(X_l^i\geq 1+\nu)\right\}\le H(\dd\rr)-\displaystyle\frac{\lm}{2}[\nu-\ln(1+\nu)]+\eta,
\end{equation}
for all $n$ sufficiently large. By the definition of $\IU(\dd,\rr,\lm)$ in (\ref{udef}), and since $[\nu-\ln(1+\nu)]$ is strictly increasing on $\nu>0$, then, for any $\e>0$, setting $\nu:=\nu^{\ast}=\IU(\dd,\rr,\lm)+\e$ and choosing $\eta$ sufficiently small in (\ref{loglimit}) ensures
$$\frac{1}{n}\ln\PP\left\{\cup_{i\in S_n}(X_l^i\geq 1+\nu^{\ast})\right\}\le-c_Q\;\;\;\;\mbox{for all $n$ sufficiently large},$$
where $c_Q$ is some positive constant, from which it follows that
$$\PP\left\{\cup_{i\in S_n}(X_l^i\geq 1+\nu^{\ast})\right\}\le e^{-c_Q\cdot n}\;\;\;\;\mbox{for all $n$ sufficiently large},$$
and (\ref{chisqresult1}) follows. Combining the same union bound argument with the lower tail result of Lemma~\ref{chisqlem} shows that, if we take $\nu^{\ast}=\IL(\dd,\rr,\lm)+\e$ for some $\e>0$, then
$$\frac{1}{n}\ln\PP\left\{\cup_{i\in S_n}(X_l^i\le 1-\nu^{\ast})\right\}\le-c_P\;\;\;\;\mbox{for all $n$ sufficiently large},$$
where $c_P$ is some positive constant, and (\ref{chisqresult2}) follows similarly to (\ref{chisqresult1}).\hfill$\Box$

For the $\FF$-distribution, we need an asymptotic result concerning the regularized incomplete beta function~\cite{temme82}, which we denote by $I_{\beta}(d_1,d_2)$.

\begin{lemma}[\textbf{Beta asymptotic~\cite[Section 3.3.2 and (2.20)]{temme82}}]\label{temmebeta}
For $d_1>d_2>0$,
\begin{equation}\label{Iexpansion}
I_{\beta}(d_1,d_2)=\frac{1}{2}\mbox{erfc}\left(\eta_I\sqrt{\frac{d_1+d_2}{2}}\right)+S_n(\eta_I)
\end{equation}
where
\begin{equation}\label{etaIdef}
-\frac{1}{2}\eta_I^2=\left(\frac{d_1}{d_1+d_2}\right)\ln\left[\frac{\beta(d_1+d_2)}{d_1}\right]+\left(\frac{d_2}{d_1+d_2}\right)\ln\left[\frac{(1-\beta)(d_1+d_2)}{d_2}\right],
\end{equation}
where
\begin{equation}\label{etasgn}
\mbox{sgn}(\eta_I)=\mbox{sgn}\left(\beta-\frac{d_1}{d_1+d_2}\right),
\end{equation}
and where $S_n(\cdot)$ is a residual term. Furthermore,
\begin{equation}\label{betaresidual}
S_n(\eta_I)=\mathcal{O}\left(\frac{1}{\sqrt{s}}\right)e^{-\frac{1}{2}s\eta_I^2}\;\;\;\;\mbox{for $l$ sufficiently large},
\end{equation}
uniformly in $\eta_I$ on compactly-supported subsets of $\RR$.
\end{lemma}

\begin{lemma}\label{Fdistlem}
Let the random variable $X_n\sim\displaystyle\frac{k}{n-k+1}\;\mathcal{F}(k,n-k+1).$ Provided
\begin{equation}\label{f_restrict}
f>\frac{\rr}{1-\rr},
\end{equation}
in the proportional-growth asymptotic,
\begin{equation}\label{loglimit3}
\lim_{n\rightarrow\infty}\frac{1}{n}\ln\PP(X_n\geq f)=\frac{1}{2}\left[\ln(1+f)-\rr\ln f-H(\rr)\right].
\end{equation}
\end{lemma}

\proof{We have
\begin{equation}\label{RIBF}
\PP[\mathcal{F}(d_1,d_2)\geq\beta]=I_{\left(\frac{d_2}{d_1\beta+d_2}\right)}\left(\frac{d_2}{2},\frac{d_1}{2}\right),
\end{equation}
where the first step follows from the definition of $X_n$, and the second step follows from the properties of the $\FF$-distribution. Now $n\geq 2k$, and therefore $\frac{n-k+1}{2}>\frac{k}{2}$, and so we may apply (\ref{Iexpansion}) with $d_1=k$, $d_2=n-k+1$ and $\beta=\left(\frac{n-k+1}{k}\right)f$ to the right-hand side of (\ref{RIBF}) to obtain
\begin{equation}\label{temme82result2}
\PP[\mathcal{F}(d_1,d_2)\geq\beta]=\frac{1}{2}\mbox{erfc}\left(-\frac{\eta_I}{2}\sqrt{n+1}\right)+S_n(\eta_I),
\end{equation}
where
\begin{equation}\label{eta2}
-\frac{1}{2}\eta_I^2=\left(\frac{n-k+1}{n+1}\right)\ln\left[\frac{n+1}{(n-k+1)(1+f)}\right]+\left(\frac{k}{n+1}\right)\ln\left[\frac{(n+1)f}{k(1+f)}\right],
\end{equation}
and where
\begin{equation}\label{signeta}
\mbox{sgn}(\eta_I)=\mbox{sgn}\left(\frac{1}{1+f}-\frac{n-k+1}{n+1}\right).
\end{equation}
By (\ref{f_restrict}), $f>\rr/(1-\rr)$, which may be combined with the observation that
$$\frac{1}{1+f}-\frac{n-k+1}{n+1}<0\iff f>\frac{k}{n-k+1},$$
to deduce that $\eta_I<0$ for $(k,n)$ sufficiently large, and therefore that
\begin{eqnarray}
\bar{\eta_I}^2\eqdef\lim_{n\rightarrow\infty}\eta_I^2&=&2\left\{(1-\rr)\ln[(1-\rr)(1+f)]+\rr\ln\left[\frac{\rr(1+f)}{f}\right]\right\}\nonumber\\
&=&2\left[(1-\rr)\ln(1-\rr)+(1-\rr)\ln(1+f)+\rr\ln\rr+\rr\ln(1+f)-\rr\ln f\right]\nonumber\\
&=&2\left[\ln(1+f)-\rr\ln f-H(\rr)\right].\label{etalimit}
\end{eqnarray}
Combining (\ref{temme82result2}) with a standard exponential tail bound on the complementary error function $\mbox{erfc}$ gives
\begin{equation}\label{temme82erfc2}
\PP(X_n\geq f)\le\frac{1}{2}e^{-\frac{1}{4}(n+1)\eta_I^2}+S_n(\eta_I),
\end{equation}
to which we can apply (\ref{betaresidual}) to obtain
$$\PP(X_n\geq f)\le\mathcal{O}(1)e^{-\frac{1}{4}(n+1)\eta_I^2}\;\;\;\;\mbox{for all $n$ sufficiently large}.$$
Taking logarithms and letting $n\ra\infty$, we deduce
$$\displaystyle\lim_{n\rightarrow\infty}\frac{1}{n}\ln\PP(X_n\geq f)\le\displaystyle\lim_{n\rightarrow\infty}\frac{1}{n}\ln\mathcal{O}(1)+\displaystyle\lim_{n\rightarrow\infty}\frac{1}{n}\cdot-\frac{1}{4}k\eta_I^2=-\displaystyle\frac{\rr}{4}\bar{\eta_I}^2,$$
which together with (\ref{etalimit}) proves (\ref{loglimit3}).}

{\bf Proof of Lemma \ref{Fdist} (Large deviation result for $\FF$).}
Union bounding $\PP(X_n^i\geq 1+f)$ over all $i\in S_n$ gives
\begin{equation}\label{combsum2}
\PP\left\{\bigcup_{i\in S_n}(X_n^i\geq f)\right\}\le\sum_{i\in S_n}\PP\left(X_n^i\geq f\right)=|S_n|\cdot\PP(X_n^1\geq f),
\end{equation}
Taking logarithms and limits of the right-hand side of (\ref{combsum2}), using (\ref{loglimit3}) and (\ref{comblim}), we have
$$\displaystyle\lim_{n\ra\infty}\frac{1}{n}\ln\left[|S_n|\cdot\PP(X_n^1\geq f)\right]=\displaystyle H(\dd\rr)-\displaystyle\frac{1}{2}\left[\ln(1+f)-\rr\ln f-H(\rr)\right],$$
which combines with (\ref{combsum2}) to imply that, for any $\eta>0$,
\begin{equation}\label{Feqn}
\frac{1}{n}\ln\PP\left\{\cup_{i\in S_n}(X_n^i\geq f)\right\}\le\displaystyle H(\dd\rr)-\displaystyle\frac{1}{2}\left[\ln(1+f)-\rr\ln f-H(\rr)\right]+\eta,
\end{equation}
for all $n$ sufficiently large. By the definition of $\IFF(\dd,\rr)$ in (\ref{Fdef}), and since the left-hand side of (\ref{Fdef}) on $f>\displaystyle\frac{\rr}{1-\rr}$ is strictly increasing in $f$, then, for any $\e>0$, setting $f:=f^{\ast}=\IFF(\dd,\rr)+\e$ and choosing $\eta$ sufficiently small in (\ref{Feqn}) ensures
$$\frac{1}{n}\ln\PP\left\{\cup_{i\in S_n}(X_n^i\geq f^{\ast})\right\}\le-c_I\;\;\;\;\mbox{for all $n$ sufficiently large},$$
where $c_I$ is some positive constant, from which it follows that
$$\PP\left\{\cup_{i\in S_n}(X_n^i\geq f^{\ast})\right\}\le e^{-c_I\cdot n}\;\;\;\;\mbox{for all $n$ sufficiently large},$$
and (\ref{Fresult}) now follows.\hfill$\Box$

\newappendix{B}
{\Large{\bf Proof of recovery results for N-IHT}}
\vspace*{0.25cm}

{\bf Roadmap for the results in this section.} Here we prove the results stated in Section \ref{NIHT_proofs}.
In the case of N-IHT, it is possible to prove convergence to an $\au(\dd,\rr;\e)$-stable point, where
\begin{equation}\label{au_eps}
\au(\dd,\rr;\e)\eqdef\{\kappa[1+\UU(\dd,2\rr)+\e]\}^{-1},
\end{equation}
for some $\e>0$. 

The proof of Theorem~\ref{recov2noise} for N-IHT takes broadly the same approach as for the corresponding result for IHT in Section~\ref{IHT_proofs}. However, in order to finally eliminate the dependence upon $\e$ in $\au(\dd,\rr;\e)$, the results corresponding to Lemmas~\ref{theta1lem} and~\ref{stableconvlem} for IHT need to be combined together. This is accomplished by Lemma~\ref{theta1lemNIHT}, which establishes that, provided (\ref{NIHTstablecond}) holds and $\e$ is taken sufficiently small, N-IHT converges to an $\au(\dd,\rr;\e)$-stable point on some $\Gm_i$ such that $i\in\Tb$ (the N-IHT support set partition is given in \req{thetadefn2} below). Lemma~\ref{theta2lemNIHT} corresponds to Lemma~\ref{theta2lem} for IHT, giving bounds on the approximation error of an $\au(\dd,\rr;\e)$-stable point on some $\Gm_i$ such that $i\in\Tb$, for any $\e>0$. Combining the two lemmas leads us to conclude that N-IHT converges to some limit point with bounded approximation error. We write $NSP_{\au}$ for the event that there is no $\au(\dd,\rr;\e)$-stable point on any $\Gm_i$ such that $i\in\Ta$.

We next introduce the support set partition definition relevant for N-IHT.

\begin{definition}[\textbf{Support set partition for N-IHT}]\label{zeta_def2}
Suppose $\dd\in(0,1]$ and $\rr\in(0,1/2]$. Given $\zeta>0$, let us write
\begin{equation}\label{astardef2}
\astar:=a(\dd,\rr)+\zeta,
\end{equation}
where $a(\dd,\rr)$ is defined in (\ref{adefn2}), let us write $\{\Gm_i:i\in S_n\}$ for the set of all possible support sets of cardinality $k$, and let us disjointly partition $S_n\eqdef\Ta\cup\Tb$ such that
\begin{equation}\label{thetadefn2}
\Ta\eqdef\left\{i\in S_n\;\;:\;\;\|\xs_{\Lm\sm\Gm_i}\|>\Sg\cdot\astar\right\};\;\;\;\;\Tb\eqdef\left\{i\in S_n\;\;:\;\;\|\xs_{\Lm\sm\Gm_i}\|\le\Sg\cdot\astar\right\}.
\end{equation}
\end{definition}

\begin{lemma}\label{theta1lemNIHT}
Choose $\zeta>0$. Suppose Assumptions {\bf A.2} and {\bf A.3} hold, and suppose that (\ref{NIHTstablecond}) holds.
Then there exists $\e$ such that, in the proportional-growth asymptotic, N-IHT converges to an $\au(\dd,\rr;\e)$-stable point on some $\Gm_i$ such that $i\in\Tb$, with probability tending to $1$ exponentially in $n$.
\end{lemma}

\proof{Under Assumption {\bf A.2}, we have by Theorem~\ref{conv3} convergence of N-IHT to a $[\kappa(1+U_{2k})]^{-1}$-stable point. By Definition~\ref{stable}, for any $\al_1<\al_2$, the set of $\al_1$-stable points includes the set of $\al_2$-stable points, and this observation combines with Lemma~\ref{RIP_bounds} to imply convergence to an $\au(\dd,\rr;\e)$-stable point, where $\au(\dd,\rr;\e)$ is defined in (\ref{au_eps}), with probability tending to $1$ exponentially in $n$. We now rehearse the argument of Lemma~\ref{theta1lem} to show that, provided $\e$ is taken sufficiently small, this stable point must be on $\Gm_i$ such that $i\in\Tb$. For any $\Gm_i$ such that $i\in\Ta$, we have $\Gm_i\neq\Lm$, and we may therefore use Theorem~\ref{singlefp} and Lemma~\ref{dist} with $\Gm:=\Gm_i$ to deduce that, given some $\e>0$, a necessary condition for there to be an $\au(\dd,\rr;\e)$-stable point on $\Gm_i$ is
\begin{equation}\label{stablecond2_NIHT}
\begin{array}{l}
\|\xs_{\Lm\sm\Gm_i}\|\cdot\sqrt{F_{\Gm_i}}+\sti\cdot\sqrt{G_{\Gm_i}}\\
\;\;\;\;\;\;\;\;\geq\au(\dd,\rr;\e)\left[\left(\frac{n-k}{n}\right)\|\xs_{\Lm\sm\Gm_i}\|\cdot R_{\Gm_i}-\sti\cdot\sqrt{\frac{k(n-k)}{n^2}\cdot S_{\Gm_i}\cdot T_{\Gm_i}}\right],\end{array}
\end{equation}
where
$$\begin{array}{c}F_{\Gm_i}\sim\displaystyle\frac{k}{n-k+1}\mathcal{F}(k,n-k+1);\;\;\;\;G_{\Gm_i}\sim\displaystyle\frac{k}{n-k+1}\mathcal{F}(k,n-k+1);\\
R_{\Gm_i}\sim\displaystyle\frac{1}{n-k}\chi^2_{n-k};\;\;\;\;S_{\Gm_i}\sim\displaystyle\frac{1}{n-k}\chi^2_{n-k};\;\;\;\;T_{\Gm_i}\sim\displaystyle\frac{1}{k}\chi^2_{k}.\end{array}$$
We also have, by (\ref{thetadefn2}),
\begin{equation}\label{Sgbound2}
\Sg\le\frac{\|\xs_{\Lm\sm\Gm_i}\|}{\astar}
\end{equation}
for any $\Gm_i$ such that $i\in\Ta$. Since $\Gm_i\neq\Lm$, $\|\xs_{\Lm\sm\Gm}\|>0$, and substitution of (\ref{Sgbound2}) into (\ref{stablecond2_NIHT}), rearrangement and division by $\|\xs_{\Lm\sm\Gm_i}\|$ yields
$$\astar\left[\au(\dd,\rr;\e)\left(\frac{n-k}{n}\right)\cdot R_{\Gm_i}-\sqrt{F_{\Gm_i}}\right]\le\sqrt{G_{\Gm_i}}+\au(\dd,\rr;\e)\sqrt{\frac{k(n-k)}{n^2}\cdot S_{\Gm_i}\cdot T_{\Gm_i}},$$
and consequently
\begin{eqnarray}\label{necsplit0_NIHT}
\PP(\bNN)&=&\PP\left\{\cup_{i\in\Ta}(\exists\;\mbox{an $\au(\dd,\rr;\e)$-stable point
  supported on}\;\Gm_i)\right\}\nonumber\\
&\le&\PP\left\{\cup_{i\in\Ta}\left(\Phi[\rr_n,F_{\Gm_i},G_{\Gm_i},R_{\Gm_i},S_{\Gm_i},T_{\Gm_i}]\geq
  0\right)\right\},
\end{eqnarray}
where we write $\rr_n$ for the sequence of values of the ratio $k/n$, and where
\begin{equation}\label{fdef2}
\Phi[\rr,F,G,R,S,T]\eqdef\sqrt{G}+\au(\dd,\rr;\e)\sqrt{\rr(1-\rr)(S)(T)}-\astar\cdot\left[\au(\dd,\rr;\e)(1-\rr)\cdot R-\sqrt{F}\right].
\end{equation}
We now define
\begin{equation}\label{astdef2}
\begin{array}{c}
F^{\ast}=G^{\ast}\eqdef\IFF(\dd,\rr)+\e;\;\;\;\;\;R^{\ast}\eqdef 1-\IL(\dd,\rr,1-\rr)-\e;\\
S^{\ast}\eqdef 1+\IU(\dd,\rr,1-\rr)+\e;\;\;\;\;\;T^{\ast}\eqdef 1+\IU(\dd,\rr,\rr)+\e.\end{array}
\end{equation}
Using (\ref{astdef2}), we deduce from (\ref{necsplit0_NIHT}) that
\begin{eqnarray}
&&\PP(\bNN)\nonumber\\
&\le&\PP\left\{\cup_{i\in\Ta}\left(\Phi[\rr_n,F_{\Gm_i},G_{\Gm_i},R_{\Gm_i},S_{\Gm_i},T_{\Gm_i}]\geq \Phi[\rr_n,F^{\ast},G^{\ast},R^{\ast},S^{\ast},T^{\ast}]\right)\right\}\label{necsplitnoise1_NIHT}\\
&+&\PP\left\{\Phi[\rr_n,F^{\ast},G^{\ast},R^{\ast},S^{\ast},T^{\ast}]\geq \Phi[\rr,F^{\ast},G^{\ast},R^{\ast},S^{\ast},T^{\ast}]+\e\right\}\label{necsplitnoise2_NIHT}\\
&+&\PP\left\{\Phi[\rr,F^{\ast},G^{\ast},R^{\ast},S^{\ast},T^{\ast}]+\e\geq
0\right\}\label{necsplitnoise3_NIHT},
\end{eqnarray}
since the event in (\ref{necsplit0_NIHT}) lies in the union of the three
events in (\ref{necsplitnoise1_NIHT}), (\ref{necsplitnoise2_NIHT}) and
(\ref{necsplitnoise3_NIHT}). Now (\ref{necsplitnoise3_NIHT}) is a deterministic event, and $\astar$ has been defined in such a way that, for any $\zeta>0$, provided $\e$ is taken sufficiently small, the event has probability
$0$. This follows from (\ref{NIHTstablecond}), (\ref{adefn2}), (\ref{astardef2}), and by the continuity of $\Phi$. The event (\ref{necsplitnoise2_NIHT}) is also deterministic, and by continuity and since $\rr_n\rightarrow\rr$, it follows that there exists some $\tilde{n}$ such that
$$\PP\left\{\Phi[\rr_n,F^{\ast},G^{\ast},R^{\ast},S^{\ast},T^{\ast}]\geq \Phi[\rr,F^{\ast},G^{\ast},R^{\ast},S^{\ast},T^{\ast}]+\e\right\}=0\;\;\;\;\mbox{for all}\;\;n\geq\tilde{n}.$$
Taking limits as $n\rightarrow\infty$, the terms (\ref{necsplitnoise2_NIHT}) and
(\ref{necsplitnoise3_NIHT}) are zero, leaving only (\ref{necsplitnoise1_NIHT}), and we have
\begin{eqnarray}
&&\lim_{n\rightarrow\infty}\PP(\bNN)\nonumber\\
&\le&\lim_{n\rightarrow\infty}\PP\left\{\cup_{i\in\Ta}\left(\Phi[\rr_n,F_{\Gm_i},G_{\Gm_i},R_{\Gm_i},S_{\Gm_i},T_{\Gm_i}]\geq \Phi[\rr_n,F^{\ast},G^{\ast},R^{\ast},S^{\ast},T^{\ast}]\right)\right\}\nonumber\\
&\le&\lim_{n\rightarrow\infty}\PP\left\{\cup_{i\in\Ta}(F_{\Gm_i}\geq F^{\ast})\right\}+\lim_{n\rightarrow\infty}\PP\left\{\cup_{i\in\Ta}(G_{\Gm_i}\geq G^{\ast})\right\}+\lim_{n\rightarrow\infty}\PP\left\{\cup_{i\in\Ta}(R_{\Gm_i}\le R^{\ast})\right\}\nonumber\\
&+&\lim_{n\rightarrow\infty}\PP\left\{\cup_{i\in\Ta}(S_{\Gm_i}\geq S^{\ast})\right\}+\lim_{n\rightarrow\infty}\PP\left\{\cup_{i\in\Ta}(T_{\Gm_i}\geq T^{\ast})\right\},\label{neclemnoise_NIHT}
\end{eqnarray}
where the last line follows from the monotonicity of $\Phi$ with respect to $F$, $G$, $R$, $S$ and $T$. Since
  $\Ta\subseteq S_n$, we may apply
  Lemmas~\ref{chisq} and~\ref{Fdist} to (\ref{neclemnoise_NIHT}), and we deduce $\PP(\bNN)\ra 0$ as $n\ra\infty$, exponentially in $n$, as required.}

\begin{lemma}\label{theta2lemNIHT}
Suppose Assumptions {\bf A.2} and {\bf A.3} hold, and suppose that (\ref{NIHTstablecond}) holds. Given any $\e>0$, there exists $\zeta$ sufficiently small such that, in the proportional-growth asymptotic, any $\au(\dd,\rr;\e)$-stable point on $\Gm_i$ such that $i\in\Tb$ satisfies (\ref{error2}), with probability tending to $1$ exponentially in $n$.
\end{lemma}

\proof{Suppose $\xb$ is a minimum-norm solution on $\Gm$, so that $\xb_{\Gm}=A_{\Gm}^{\dag}b$. Then we may follow the argument of Lemma~\ref{theta2lem} to deduce (\ref{PQ_bound}), where
\begin{equation}\label{dist_def_NIHT}
P_{\Gm}\sim\frac{k}{n-k+1}\mathcal{F}(k,n-k+1);\;\;\;\;Q_{\Gm}\sim\frac{k}{n-k+1}\mathcal{F}(k,n-k+1).
\end{equation}
Combining (\ref{PQ_bound}) with (\ref{thetadefn2}), we may further deduce
\begin{eqnarray}
\|\xb-\xs\|^2&\le&\Sg^2\left[\astar\cdot\sqrt{P_{\Gm}}+\sqrt{Q_{\Gm}}\right]^2+\left[\astar\right]^2\cdot\Sg^2\nonumber\\
&=&\Sg^2\left\{\left[\astar\cdot\sqrt{P_{\Gm}}+\sqrt{Q_{\Gm}}\right]^2+\left[\astar\right]^2\right\}.\label{errorbound_NIHT}
\end{eqnarray}
For the sake of brevity, let us define
\begin{equation}\label{psi_def2}
\Psi[P,Q]\eqdef\sqrt{\left(\astar\cdot\sqrt{P}+\sqrt{Q}\right)^2+\astar^2},
\end{equation}
so that (\ref{errorbound_NIHT}) may equivalently be written as
\begin{equation}\label{errorbound_brevity2}
\|\xb-\xs\|\le\Sg\cdot\Psi\left[P_{\Gm},Q_{\Gm}\right].
\end{equation}
First suppose that $\sg>0$. Given $\zeta>0$, let us define
\begin{equation}\label{PQdef2}
P^{\ast}=Q^{\ast}\eqdef\IFF(\dd,\rr)+\zeta.
\end{equation}
Now we use (\ref{errorbound_brevity2}) to perform a union bound over all $\Gm_i$ such that $i\in\Tb$, writing $\xb_i$ for the minimum-norm solution on $\Gm_i$, giving
\begin{eqnarray}
&&\PP\left\{\exists\;\mbox{some}\;\Gm_i\;\mbox{such
    that}\;i\in\Tb\;\mbox{and}\;\|\xb_i-\xs\|>\Sg\cdot\Psi\left[P^{\ast},Q^{\ast}\right]\right\}\nonumber\\
&=&\PP\left\{\bigcup_{i\in\Tb}\left(\|\xb_i-\xs\|>\Sg\cdot\Psi\left[P^{\ast},Q^{\ast}\right]\right)\right\}\label{2necsplit0_NIHT}\\
&\le&\PP\left\{\bigcup_{i\in\Tb}\left(\|\xb_i-\xs\|>\Sg\cdot\Psi\left[P_{\Gm_i},Q_{\Gm_i}\right]\right)\right\}\label{2necsplit1_NIHT}\\
&+&\PP\left\{\bigcup_{i\in\Tb}\left(\Sg\cdot\Psi\left[P_{\Gm_i},Q_{\Gm_i}\right]\geq\Sg\cdot\Psi\left[P^{\ast},Q^{\ast}\right]\right)\right\},\nonumber\\
&&\;\label{2necsplit2_NIHT}
\end{eqnarray}
since the event in (\ref{2necsplit0_NIHT}) lies in the union of the two
events in (\ref{2necsplit1_NIHT}) and (\ref{2necsplit2_NIHT}). It is an immediate consequence of
(\ref{errorbound_brevity2}) that the event in (\ref{2necsplit1_NIHT}) has probability
$0$. Taking limits of (\ref{2necsplit2_NIHT}) as $n\rightarrow\infty$, and cancelling $\Sg$, we have
\begin{eqnarray}
&&\lim_{n\rightarrow\infty}\PP\left\{\exists\;\mbox{some}\;\Gm_i\;\mbox{such
    that}\;i\in\Tb\;\mbox{and}\;\|\xb_i-\xs\|>\Sg\cdot\Psi\left[P^{\ast},Q^{\ast}\right]\right\}\nonumber\\
&\le&\lim_{n\rightarrow\infty}\PP\left\{\bigcup_{i\in\Tb}\left(\Psi\left[P_{\Gm_i},Q_{\Gm_i}\right]\geq\Psi\left[P^{\ast},Q^{\ast}\right]\right)\right\}\nonumber\\
&\le&\lim_{n\rightarrow\infty}\PP\left\{\cup_{i\in\Tb}(P_{\Gm_i}\geq P^{\ast})\right\}+\lim_{n\rightarrow\infty}\PP\left\{\cup_{i\in\Tb}(Q_{\Gm_i}\geq Q^{\ast})\right\},\label{2neclem_NIHT}
\end{eqnarray}
where we used the monotonicity of $\Psi$ with respect to $P$ and $Q$ in the last line. Since
  $\Tb\subseteq S_n$, and using (\ref{dist_def_NIHT}), we may apply
  Lemma~\ref{Fdist} to (\ref{2neclem_NIHT}), yielding that each of the limits in the right-hand side of (\ref{2neclem_NIHT}) converges to zero exponentially in $n$, and so finally
$$\lim_{n\rightarrow\infty}\PP\left\{\exists\;\mbox{some}\;\Gm_i\;\mbox{such
    that}\;i\in\Tb\;\mbox{and}\;\|\xb_i-\xs\|>\Sg\cdot\Psi\left[P^{\ast},Q^{\ast}\right]\right\}=0,$$
with convergence at a rate exponential in $n$ also by
    Lemma~\ref{Fdist}. The same result also holds when $\sg=0$ by (\ref{errorbound_NIHT}). Since by Lemma~\ref{necfp}, any stable point is necessarily a minimum-norm solution, and recalling the definition of $\Psi(P,Q)$ in (\ref{psi_def}), and the definitions of $P^{\ast}$, $Q^{\ast}$ in (\ref{PQdef2}), we have
\begin{equation}\label{finalbound2}
\lim_{n\rightarrow\infty}\PP\left\{\begin{array}{ll}\exists\;\mbox{some $\au$-stable point $\xb_i$ on}\;\Gm_i\;\mbox{such
    that}\;i\in\Tb\;\mbox{and}\\
\;\|\xb_i-\xs\|>\Sg\sqrt{\IFF(\dd,\rr)\left[1+a(\dd,\rr)+\zeta\right]^2+\left[a(\dd,\rr)+\zeta\right]^2}\end{array}
\right\}=0,
\end{equation}
with convergence exponential in $n$. Finally, by continuity,
$$\begin{array}{l}\|\xb_i-\xs\|>\Sg\sqrt{\IFF(\dd,\rr)\left[1+a(\dd,\rr)\right]^2+1+\left[a(\dd,\rr)\right]^2}\\
\;\;\;\;\;\;\;\;\;\;\;\;\;\Longrightarrow\|\xb_i-\xs\|>\Sg\sqrt{\IFF(\dd,\rr)\left[1+a(\dd,\rr)+\zeta\right]^2+\left[a(\dd,\rr)+\zeta\right]^2},\end{array}$$
for some $\zeta$ suitably small, and the result now follows from the definition of $\xi(\dd,\rr)$ in (\ref{xidef2}).}

It is now straightforward to prove the two main results for N-IHT.

{\bf Proof of Theorem~\ref{recov2noise}:} By Lemma~\ref{theta1lemNIHT}, there exists $\e>0$ such that N-IHT converges to an $\au(\dd,\rr;\e)$-stable point on some $\Gm_i$ such that $i\in\Tb$, and for this choice of $\e$, we can apply Lemma~\ref{theta2lemNIHT} to deduce the result.

{\bf Proof of Corollary~\ref{recov2noiseless}:} The result follows by setting $\Sg\eqdef 0$ in Theorem~\ref{recov2noise}.


\end{document}